\documentclass[]{amsart}

\usepackage[utf8]{inputenc}
\usepackage[T1]{fontenc}

\usepackage{lineno,hyperref}
\usepackage{amsmath}
\usepackage{latexsym}
\usepackage{amsthm}
\usepackage{amssymb, ifthen, enumerate}
\usepackage{longtable}

\usepackage{color}
\usepackage{ulem}

\newtheorem{lemma}{Lemma}[section]
\newtheorem{thm}[lemma]{Theorem}
\newtheorem{prop}[lemma]{Proposition}
\newtheorem{cor}[lemma]{Corollary}
\newtheorem{rem}[lemma]{Remark}
\newtheorem{hyp}[lemma]{Hypothesis}
\newtheorem{definition}[lemma]{Definition}

\renewcommand{\qed}{\hspace*{\fill}\rule{1.5ex}{1.5ex}}
\newcommand{\A}{\mathrm{Alt}}     %script A fuer alt. Gruppe
\newcommand{\Alt}{\mathrm{Alt}}     %noch mal script A fuer alt. Gruppe
\newcommand{\Sym}{\mathrm{Sym}}     %script S fuer symm Gruppe
\newcommand{\ug}{\leqslant}		%UG
\newcommand{\Menge}[2]{\left\{ #1 \;\middle|\; #2 \right\}}

\newcommand{\aut}{\textrm{Aut}}

\newcommand{\cal}{\mathcal}

\newcommand{\GL}{\textrm{GL}}
\newcommand{\SL}{\textrm{SL}}

\renewcommand{\tilde}{\widetilde}

\newcommand{\PSL}{\textrm{PSL}}

\newcommand{\PSU}{\textrm{PSU}}

\newcommand{\SU}{\textrm{SU}}
\newcommand{\SUn}{\textrm{SU}}
\newcommand{\Spl}{\textrm{Sp}}

\newcommand{\GF}{\textrm{GF}}
\newcommand{\Sz}{\textrm{Sz}}
\newcommand{\GLi}{\textrm{GL}}
\newcommand{\POm}{P\Omega}
\newcommand{\SLi}{\textrm{SL}}
\newcommand{\PSp}{\textrm{PSp}}

\newcommand{\M}{\textrm{M}}

\newcommand{\syl}{\textrm{Syl}}
\newcommand{\Syl}{\textrm{Syl}}

\newcommand{\FO}{\textrm{fix}_{\Omega}}

\newcommand{\FL}{\textrm{fix}_{\Lambda}}

\newcommand{\Mat}{M}
\newcommand{\Jan}{J}

\newcommand{\HS}{HS}
\newcommand{\McL}{McL}
\newcommand{\He}{He}

\def \<{\langle }
\def \>{\rangle }
\def \N{\mathbb{N}}

\makeatletter % we need to use kernel commands

\newcommand{\abs}[1]{\left| #1 \right|}

\begin{document}

\begin{center}
\Large{\textbf{Finite simple permutation groups acting with fixity $4$}}

\vspace{0.2cm} \small{Barbara Baumeister, Paula H\"ahndel, Kay Magaard, Patrick Salfeld, and Rebecca Waldecker}
\end{center}

\normalsize
\vspace{2cm}

Abstract:

Motivated by the theory of Riemann surfaces and specifically the significance of Weierstrass points, we
classify all finite simple groups that have a faithful transitive action with fixity 4. We also give details
about all possible such actions.
\vspace{1cm}


%%%%%%%%%

\section{Introduction}

This paper continues the analysis of permutation groups where every nontrivial element has a low number of fixed points (see \cite{MW}).
We call the maximum number of fixed points of nontrivial group elements that occur in a faithful and transitive group action its \textbf{fixity}.
Our motivation to study permutation groups that act with low fixity comes from questions about Riemann surfaces:
If $X$ is a compact Riemann surface of genus at least $2$, then its automorphism group is finite. Proving this uses the existence of a finite, nonempty set of analytically distinguished points of $X$, the so-called Weierstrass points.
The link between fixed points of nontrivial automorphisms of $X$ and Weierstrass points is the motivation for our work. Schoeneberg has proven (see for example page 264 in \cite{FK}) that if a nontrivial automorphism of $X$ has at least five fixed points on $X$, then all these fixed points are Weierstrass points.
Studying groups of automorphisms of Riemann surfaces that act with fixity at most $4$ is therefore a first step towards understanding situations where Schoeneberg's result does not apply.
We refer the reader to \cite{MW} and \cite{FK} for more background on the theory of Riemann surfaces and our intended applications.

Groups of fixity $1$ are Frobenius groups, where the Frobenius kernel consists of the fixed point free elements and the Frobenius complements are the point stabilizers.
Groups of fixity $2$ and $3$ have been discussed in \cite{MW} and \cite{MW3}.
Therefore, it remains to
analyze those of fixity $4$, and the present article contributes to this analysis.
Before we continue, we would like to mention that there have been many contributions to this area of research already,
often under slightly different hypotheses. For example, there is a line of research about permutation groups
where involutions fix very few points (see\cite{Ben}, \cite{Ronse} and others) or where every nontrivial element fixes zero points or a prescribed number of points (see for example \cite{PS1}).
We refer to some, but not all of this previous work explicitly.

\smallskip
Now let $G$ be a group that acts faithfully, transitively and with fixity $4$ on a finite set $\Omega$.
%It has become apparent that strongly $p$-embedded subgroups play a major role for understanding the structure of $G$ (see Lemma~4 and Theorem~13  of \cite{BMW}). For example, if $H$ is a nontrivial four point stabilizer in $G$ and the order of $H$ is coprime to $6$, then it turns out that $N_G(H)$ is strongly $p$-embedded in $G$ for some prime $p \geq 5$.
%This occurs, in particular, if the point stabilizers in $G$ have order coprime to $6$.

Let $T$ be a Sylow $2$-subgroup of $G$.
In \cite[Theorem~13]{BMW} we prove that if $|G_\alpha|$ is even,
then either the set-wise stabilizer of the union of short $T$-orbits is strongly embedded in $G$  or  $T$ has a quite restricted structure.  The next case to consider is the situation where $|G_\alpha|$
is odd and is divisible by $3$. Let $P$ be a Sylow $3$-subgroup of $G$. Then we prove that  either the set-wise stabilizer $K$ of the union of short $P$-orbits is strongly $3$-embedded in $G$
or the structure of $P$ is very restricted. More precisely our main result about the $3$-structure is as follows:

\begin{thm}\label{3total}
Suppose that $G$ is a finite group acting transitively, faithfully and with fixity $4$ on a set~$\Omega$.
Let $P$ be a Sylow $3$-subgroup of $G$, let $\Delta$ be the union of $P$-orbits of length at most $3$ on $\Omega$ and let $K$ denote the stabilizer of the set $\Delta$ in $G$.
Moreover let $\alpha \in \Omega$. Then one of the following holds:
\begin{enumerate}
\item[(a)]  $G_\alpha$ has even order.
\item[(b)]    there is a subnormal subgroup $N$ of $G$ such that one of the following holds:
\begin{enumerate}
\item[(i)]  $N \unlhd G$, $|G:N|$ divides $3$ and $(|N_\alpha|, 6) = 1$, or
\item[(ii)]  $|G:N|$ divides $9$, $K \cap N$ is strongly $3$-embedded in $N$ and $1 \leq |\Delta| \leq 4$.
\end{enumerate}
\item[(c)]   There are four exceptional possibilities:
\begin{enumerate}
\item[(i)]  $|P| = 3$ and $|N_G(P)|_2 \le 4$,
\item[(ii)] $P$ is  elementary abelian of order $9$,
\item[(iii)] $P$ is extra-special of exponent $3$ and order $27$, or
\item[(iv)] $P \cong C_3 \wr C_3$.
\end{enumerate}
\end{enumerate}
\end{thm}

There are some more details that we discuss in Section 3, and in fact Theorem \ref{3total} is a direct consequence of
Theorem \ref{3and2} at the end of that section.
In the special case of finite simple groups, we combine this result with earlier work of Ronse's (\cite{Ronse}). This gives a stronger result and a
blueprint for the classification of all finite simple groups that act with fixity 4, along with information about the point stabilizers.

\begin{thm}\label{3and2simple}
Suppose that $G$ is a finite simple non-abelian group that acts transitively, faithfully and with fixity $4$ on a set $\Omega$.
Let $P \in \syl_3(G)$, $S \in \syl_2(G)$ and let $f$ denote the maximum number of fixed points of involutions in $G$.

 %let $\Delta$ be the union of $P$-orbits of length at most $3$ on $\Omega$ and let $K$ denote the stabilizer of the set $\Delta$ in $G$. Moreover let $\alpha \in \Omega$ and let Let $F:=\langle t,g \mid t,g \in G, o(t)=2, $\FO(t)=f\rangle$.

Then one of the following holds:
\begin{enumerate}

\item[(1)] $f \ge 1$ and $G$ has a strongly embedded subgroup.

\item[(2)] $1 \le f \le 3$ and $S$ is dihedral or semi-dihedral.

\item[(3)] $f=4$ and $G$ has sectional $2$-rank at most $4$.

\item[(4)] The order of the point stabilizers is odd and divisible by $3$. One of the following occurs:

\begin{enumerate}
\item[(a)] $G$ has a strongly $3$-embedded subgroup.
\item[(b)]  $|P| = 3$ and $|N_G(P)|_2 \le 4$.
\item[(c)] $P$ is elementary abelian of order $9$.
\item[(d)] $P$ is extra-special of exponent $3$ and order $27$.
\item[(e)] $P \cong C_3 \wr C_3$.
\end{enumerate}

\item[(5)]  The point stabilizers have order coprime to $6$.

\end{enumerate}
\end{thm}

As it will turn out later, the case where point stabilizers have order coprime to 6
requires a lot of work and leads to many interesting examples.
After a detailed analysis of all cases from Theorem \ref{3and2simple}, we obtain our main classification result:

\begin{thm}\label{main4fpsimple}
Suppose that $G$ is a finite simple non-abelian group that acts transitively, faithfully and with fixity $4$ on a set $\Omega$.
Let $\alpha \in \Omega$.

Then $G$ is isomorphic to one of the groups in the following table, with point stabilizer $G_\alpha$ as described:
\end{thm}

\begin{tabular}{l | p{8cm}| p{4cm}}

\mbox{}
	Group \(G\)                           & Point stabilizer structure  & Comments \\ \hline
	\(\Alt_6 \cong \PSL_2(9)\)            & \(C_2\), \(\Sym_3\), \(C_3 \times C_3\), \(D_{10}\), \((C_3 \times C_3) : C_2\) & \\
	\(\Alt_7\)                            & \(C_5\), \(\Alt_6\)  &\\

	\(\PSL_2(7) \cong \PSL_3(2)\)         & \(C_2\), \(\Sym_3\) &\\
	\(\PSL_2(8)\)                         & \(C_2\), \(\Sym_3\),
	\(D_{14}\), \(D_{18}\)& \\

	\(\PSL_2(11)\)                        & \(C_3\), \(\Alt_4\)  &     \\
	\(\PSL_2(13)\)                        & \(C_3\), \(\Alt_4\), \(C_{13} : C_3\)& \\
		\(\PSL_2(q)\)                        & $G_\alpha$ cyclic of order $\frac{q-1}{4}$ & $q \ge 17$, $q \equiv 1$ mod 4\\
		& $G_\alpha$ has index 2 in the normalizer of a Sylow subgroup in defining characteristic & $q \ge 17$, $q \equiv 1$ mod 4\\
				                        & $G_\alpha$ cyclic of order $\frac{q+1}{4}$& $q \ge 17$, $q \equiv -1$ mod 4\\
	
	\(\PSU_3(3)\)                         & \(((C_3 \times C_3) :C_3) :C_8\) &            \\
	
	\(\PSU_4(2)\cong \PSp_4(3)\)          & \(C_5\)       &                    \\
		\(\PSU_4(3)\)          & \(C_5\)                   &        \\
			\(\PSp_4(q)\)          & $G_\alpha$ cyclic of order $\frac{q^2+1}{(2,q+1)}$         &      $q \ge 3$            \\

		\(\Sz(q)\)                            & $G_\alpha$ cyclic of order $q+\sqrt{2q}+1$ or
		$q-\sqrt{2q}+1$        & $q \ge 8$ a power of 2 with odd exponent                 \\
			\(\POm^-_8(q)\)          &   $G_\alpha$ cyclic of order $\frac{q^4+1}{(2,q+1)}$   &                \\
			\(^3D_4(q)\)          &  $G_\alpha$ cyclic of order $q^4-q^2+1$        &       \\
			\(^2G_2(q)\)          & $G_\alpha$ Frobenius group of order $q^3\cdot \frac{q-1}{2}$ &\\
			& or cyclic of order $\frac{q-1}{2}$          &     $q \ge 27$ a power of $3$ with odd exponent        \\
	
	\(M_{11}\)                            & \(C_5\), \(C_{11} :C_5\), \(\PSL_2(11)\) &                                                                                 \\
	\(M_{12}\)                            & \(M_{11}\)           &     two conjugacy classes of possible point stabilizers\\
	\(M_{22}\)                            & \(C_5\), \(C_{11}:C_5\)       &    \\

	\(J_1\)                               & \(C_{15}\)         &    \\

\end{tabular}	\label{MainTable}

\vspace{2cm}

This article is organised as follows:

Section 2 is short, with some notation, a technical lemma and an initial result about the 3-structure of the groups that we consider.
Then we analyze the 3-structure more deeply in Section 3, which leads to a natural case distinction towards the classification of finite simple groups that act with fixity 4.
For very small groups, and also for some individual arguments here and there, we use
\texttt{GAP} (\cite{GAP}). The code is presented and explained in an Appendix at the very end, and the results are captured in Table \ref{SmallGroupsTable} in Remark \ref{rem4.1} in Section 4, along with some more technical results for preparation.
Some groups, in particular Lie type groups of small Lie rank, occur several times in Theorem \ref{3and2simple}, which is why we begin our analysis with those, in Section 5.
After that we look at the cases (3), (4) and (5) of Theorem \ref{3and2simple}, in this order.
Section~6 looks at the case where some involution fixes four points.
Concerning this case we should mention previous work by Buekenhout and Rowlinson (see \cite{Ro1}, \cite{BuR1} and \cite{BuR2})
that we are aware of, but did not cite for our classification.
There are two reasons:
The authors refer to unpublished work by Fong, and also, the group
$\PSL_2(9)$ with point stabilizers isomorphic to $\Sym_3$ is missing in their classification.
So we thought it would be more reasonable to give an independent proof, reproducing and correcting their list,
based on the work of Gorenstein and Harada which they also rely on.

Section 7 considers point stabilizers of odd order, but order divisible by 3, and then the bulk of the work is done in Section 8 in the case of point stabilizers of order coprime to 6.
This is where we rely heavily, and very visibly, on the Classification of Finite Simple Groups (CFSG) as presented in the GLS series.

Section 9 collects all the intermediate results and finishes the proof of Theorem ~\ref{main4fpsimple}. The appendix explains how we use \texttt{GAP} and gives the fixed point profiles of the groups in the table in Theorem \ref{main4fpsimple} for future reference.

\smallskip
\textbf{Acknowledgments.}
The project on permutation groups with low fixity started in 2012, initiated by Kay Magaard. It was supported by the DFG for several years and also partly
by the National Science Foundation under Grant No. DMS-1440140 while
he was in residence at the Mathematical Sciences Research Institute in Berkeley, California, during
the Spring semester 2018.
We continue our work on this project after Kay's unexpected passing in 2018.
Finally, we would like to thank Chris Parker for pointing out mistakes, for making valuable suggestions for improvement and for contributing Lemma \ref{NormalIndex3Subgroup}.

%%%%%%%%%%%%%%%%

\section{Preliminaries}

In this paper, by ``group'' we always mean a finite
group, and by ``permutation group'' we always mean a group that acts
faithfully on a set.
Throughout $\Omega$ denotes a finite set, and $G$ denotes a
permutation group on $\Omega$, unless stated otherwise.\\

Let $\omega \in \Omega$ and $g \in G$, and moreover let $\Lambda
\subseteq \Omega$ and $H \le G$.
We use standard notation for orbits and point stabilizers, and moreover we write
$\FL(H):=\{\omega \in \Lambda \mid \omega^h=\omega$ for all $h \in
H\}$ for the \textbf{fixed point set of $H$ in $\Lambda$}.
When we say that $H \le G$ \textbf{is a nontrivial four point stabilizer}, then we mean that there is a set $\Delta \subseteq \Omega$ such that $|\Delta|=4$ and such that $1 \neq H$ is exactly the point-wise stabilizer of $\Delta$ in $G$.

If $g \in G$, then we write
$\FL(g)$ instead of $\FL(\<g\>)$. By
$\pi(G)$ we denote the set of prime divisors of $|G|$.
Whenever $n,m$ are natural numbers and $p$ is a prime number, then $(n,m)$ denotes the largest natural
common divisor of $n$ and $m$ and $n_p$ is the largest power of $p$ dividing $n$.
We use the notation $C_n$ for cyclic groups of order $n$.

\begin{definition}
Let $k$ be a non-negative integer and suppose that the group $G$ acts on the finite set $\Omega$. We say that $G$ has \textbf{fixity $k$ on $\Omega$} if and only
if there is some element of $G^\#$ that fixes exactly $k$ distinct points on $\Omega$ and if no element of $G^\#$ fixes more than $k$ distinct points.
\end{definition}

\begin{hyp}\label{4fix}
Suppose that $G$ is a finite group that acts transitively, faithfully and with fixity $4$ on a set $\Omega$.
\end{hyp}

\begin{lemma}\label{Inv-lem:TB}
Suppose that Hypothesis \ref{4fix} holds and let $\alpha \in \Omega$.

If $1 \neq Y \le G_\alpha$, then $|N_G(Y):N_{G_\alpha}(Y)| \le |\FO(Y)| \le 4$.

\smallskip
Moreover let $H \leq G_\alpha$ be a nontrivial four point stabilizer in $G$ and let $\Delta:=\FO(H)$.

Then
the following hold:

(a) If $p \in \pi(G_\alpha)$ and $p\ge 5$, then $G_\alpha$ contains a Sylow $p$-subgroup of $G$.
If $p \in \pi(H)$ and $p\ge 5$, then $H$ contains a Sylow $p$-subgroup of $G$.
In particular, if the point stabilizers (or nontrivial four point stabilizers) have order coprime to 6,
then they are Hall subgroups of $G$.

(b) Suppose that $1 \neq X \le H$. Then $N_G(X)$ stabilizes $\FO(X)$ and is (therefore) contained in $N_G(H)$.
Moreover $\abs{N_G(H):N_{G_\alpha}(H)}\in\{2,4\}$. If $3 \in \pi(N_G(X))$, then $3 \in \pi(G_\alpha)$.

(c) If $X \le G_\alpha$ fixes exactly three points, then $\abs{N_G(X):N_{G_\alpha}(X)}\leq 3$.

(d)
If $3 \in \pi(H)$, then $H$ contains a Sylow $3$-subgroup of $G$ or $N_G(H)$ has a subgroup that induces $\A_4$ on $\FO(H)$.

(e) If $H$ contains a Sylow subgroup of $G_\alpha$, then $N_G(H)$ acts transitively on $\Delta$.

(f)
Either $|N_G(H):N_{G_\alpha}(H)| = 4$ or
$H$ is a $2$-group, but not a Sylow $2$-subgroup of $G$.
In particular $|G|$ is divisible by $4$.

(g) $H$ is a TI-group.

(h) If $|H|$ is coprime to $6$ and $p\in \pi(H)$, then $N_G(H)$ is strongly
$p$-embedded in $G$.

\end{lemma}

\begin{proof}
For the first assertion we note that $N_G(Y)$ stabilizes $\FO(Y)$, a set of size at most 4 that contains $\alpha$.
Then the $N_G(Y)$-orbit of $\alpha$ is contained in $\FO(Y)$, which gives that
$4 \ge |\FO(Y)| \ge |\alpha^{N_G(Y)}|=|N_G(Y):N_G(Y) \cap G_\alpha|$.

The statement in (a) is proven in Lemma 4 of \cite{BMW}.
Lemma 6 of the same article covers the statements in (b),  (d) and (g) as written.

(c) is a special case of the general statement at the beginning of the lemma.
%we let  $X \le G_\alpha$ be such that $X$ fixes exactly three points $\alpha,\beta,\gamma \in \Omega$. Then $N_G(X)/X$ induces a subgroup of $\Sym_3$ on $\FO(X)$. If this factor group is isomorphic to \(\Sym_3\), then there exists an element in \(N_{G_\alpha}(X)\) that interchanges $\beta$ and $\gamma$ and therefore $\abs{N_G(X)/ N_{G_\alpha}(X)}=3$. Otherwise $N_G(X)/X$ is isomorphic to a proper subgroup of \(\Sym_3\). In both cases $\abs{N_G(X):N_{G_\alpha}(X)}\leq 3$, as stated.

We turn to (e).
Let $r \in \pi(G_\alpha)$ be such that $H$ contains a Sylow $r$-subgroup $R$ of $G_\alpha$. Then
let $\beta \in \Delta$ be such that $\alpha \neq \beta$ and let $g \in G$ be such that $\alpha^g = \beta$ (by transitivity of $G$).
Now $R^g \leq G_\beta$, but also $R \le H \le G_\beta$ and therefore we find $h \in G_\beta$ such that $R^{gh}=R$.
Then $1 \neq R \le H$ and (a) yields that $gh \in N_G(H)$, so we proved that $N_G(H)$ acts transitively on $\Delta$.

For (f) we first suppose that
$H$ is a $\{2,3\}$-group.
If $H$ is a $2$-group, then either (e) is applicable or the other statement holds automatically.
So we suppose that $3 \in \pi(H)$. If $H$ contains a Sylow $3$-subgroup of $G$, then (e) yields the result.
Otherwise we use (d) and then $N_G(H)$ has a subgroup that induces $\Alt(4)$ on $\Delta$. In particular
$N_G(H)$ acts transitively on $\Delta$ and therefore
$|N_G(H):N_{G_\alpha}(H)|=4$.
Finally, if  $p \in \pi(H)$ is such that $p \ge 5$, then
by (a) it contains a Sylow $p$-subgroup of $G$, and (e) applies again.

Now we look at (f):
If $H$ has odd order, then $|G|$ is divisible by $4$ by (d) or (e).
Suppose that $|H|$ is even.
If $G$ has Sylow $2$-subgroup of order
$2$, then $|N_G(H):N_{G_\alpha}(H)| \neq 4$, so by the first statement we know that
 $H$ is a $2$-group, but not a Sylow $2$-subgroup of $G$, which is not
possible because $H \neq 1$.
Therefore $|G|$ is divisible by $4$.

Finally, we prove (h). Of course $N_G(H)$ has order divisible by $p$.
Let $g \in G \setminus N_G(H)$ and suppose that $x\in N_G(H) \cap N_G(H)^g$ is a $p$-element.
By Hypothesis $p \ge 5$, so (f) yields that $x\in H \cap H^g$.
Then $x$ fixes all elements in $\Delta \cup \Delta^g$, so $x=1$ or
$\Delta=\Delta^g$. But if $g$ stabilizes $\Delta$, then it normalizes $H$.
Consequently $N_G(H) \cap N_G(H)^g$ is indeed a $p'$-group (see also
\cite[Proposition]{GLS2}).
\end{proof}

The following lemma will be used frequently for our strategy to connect the 3-structure and the 2-structure of the groups that we analyze.

\begin{lemma}\label{syl3}
Suppose that Hypothesis \ref{4fix} holds and let $P \in \syl_3(G)$.
Let $\Delta$ be the union of all $P$-orbits of $\Omega$ of size at most $3$.
Then one of the following holds:

\begin{enumerate}
\item[(a)] All $P$-orbits are regular and the point stabilizers in $G$ are $3'$-groups.

\item[(b)] $|\Delta| > 4$ and $|P| \le 9$.

\item[(c)] $|\Delta| \leq 4$ and $P$ is of maximal class. There exists some non-regular $P$-orbit on $\Omega \setminus \Delta$ and for every such orbit $\Lambda$ and all  $\lambda \in \Lambda$ it is true that
$|P_\lambda| = 3$ and that $P_\lambda$ fixes exactly three points on
$\Lambda$.

\item[(d)] $\Delta$ is the unique $P$-orbit of length 3 and all orbits of $P$ on
$\Omega \setminus \Delta$ are regular.

\item[(e)] $1 \le |\Delta| \le 4$, there is some $\delta \in \Delta$ such that $P \le G_\delta$,
 and all $P$-orbits on $\Omega \setminus \Delta$ are regular.

\end{enumerate}
In (c), (d) and (e) we see that $P$ possesses an orbit of length at least $9$ and therefore
$|P| \ge 9$.
\end{lemma}

\begin{proof}
 Lemma 10 in \cite{BMW}.
\end{proof}

Here comes a Frattini type lemma that will be useful in several places:

\begin{lemma}\label{helpFrattini}
Let $A$ be a soluble group and $p,q \in \pi(A)$ be two different prime divisors of the order of $A$. Then
there is a prime $r \in \{p,q\}$ and a nontrivial $r$-subgroup $R$ of $A$ that is normalized by
some Sylow $s$-subgroup of $A$, where $\{r,s\} = \{p,q\}$.
\end{lemma}

\begin{proof} By hypothesis $A$ is soluble, so we find $m \in \N$ and normal subgroups $1 = A_1 \le  A_2 \le \cdots \le A_m = A$
of $A$ such that the sections $A_i/A_{i-1}$, where $i\in \{2,...,m\}$, are elementary abelian groups of prime
power order. Let $j \le m$ be minimal such that $r$ divides $|A_j|$ for some $r \in  \{p,q\}$, and let
$R$ be a Sylow $r$-subgroup of $A_j$. Then the Frattini argument gives that $A = A_j N_A(R)$.
Let $s$ be such that $\{r,s\} = \{p,q\}$. Then  $s$ does not divide $|A_j|$, and therefore our assertion holds.
\end{proof}

Finally, we will frequently refer to background knowledge about finite simple groups,
and we often use notation from \cite{Atlas} and \cite{Wilson}. This includes the subgroup structure and isomorphisms, as is explained for example on p. x and xi of \cite{Atlas}, and for details we will give individual references.

%%%%%%%%%%%%%%

\section{From the 3-structure to the 2-structure}

The purpose of this section is to analyze how, given Hypothesis \ref{4fix}, the point stabilizer structure and the 3-structure of $G$
influence each other. It will turn out that the 2-structure of $G$ comes into play here as well, and at the end of the section we will be able to prove a theorem with a natural case distinction along which we can organise the analysis in the following sections. In fact, the natural case distinction that occurs determines the structure of the remainder of the paper and the strategy for the proof of Theorems \ref{3and2simple} and \ref{main4fpsimple}.

First we remark that Lemma \ref{syl3} gives information about the possible orbit sizes for a Sylow $3$-subgroup~$P$ of $G$ on $\Omega$.
In addition to orbits of size $|P|$ and $1$, we see in Part (d) of the lemma that orbit size $3$ is possible, and the only case where another orbit size occurs is Part (c), with orbits of length $|P|/3$.

Most of this section is devoted to a closer look at Cases (c)-(e) of the lemma.

\subsection{Cases (d) and (e) of Lemma \ref{syl3}}

\begin{hyp}\label{hyp3}
In addition to Hypothesis \ref{4fix}, let $P \in \syl_3(G)$, let $\Delta$ denote the union of the $P$-orbits of length at most $3$ and suppose that $P$ satisfies Lemma \ref{syl3}\,(d) or (e).
We define $D$ to be the element-wise stabilizer of $\Delta$ in $G$, $K$ to be the set-wise
stabilizer of $\Delta$ in $G$, and we set $Q:=P \cap D$.
\end{hyp}

\begin{lemma} \label{littlehelp}
Suppose that Hypothesis \ref{hyp3} holds.

(a) If $g \in G \setminus K$, then $D \cap D^g$ is a $3'$-group.

(b) If $|\Delta|=3$, then all $3$-elements of $K \setminus D$ act fixed-point-freely on $\Omega$, whereas the $3$-elements
of $D$ fix exactly three points on $\Omega$.

(c) If $|\Delta|=3$ and  Case (e) of Lemma \ref{syl3} holds, then $P \leq D$.

(d) If $|\Delta| \neq 3$, then Lemma \ref{syl3}~(e) is true. Moreover, if $3$ divides $|K/D|$, then $|\Delta| = 4$ and every $3$-element of $K \setminus D$ has a unique fixed point on $\Omega$.

(e) The $3$-elements of $D$ and of $K \setminus D$ can be distinguished by the number of fixed points on $\Omega$.
In particular $x^G \cap D= \varnothing$ for all $3$-elements $x \in K \setminus D$.
\end{lemma}

\begin{proof}
Throughout, we recall that  \(P\leq K\), because $\Delta$ is a union of $P$-orbits,
and that $P$ acts semi-regularly on $\Omega \setminus \Delta$.

(a) Let $g \in G$ and let $y$ be a nontrivial $3$-element in $D \cap D^g$. Without loss $y \in P$,
and then the fact that $P$ acts semi-regularly on $\Omega \setminus \Delta$ implies that
$y$ cannot fix any points outside of $\Delta$. Hence $\FO(y) \subseteq \Delta$.
But $y$ fixes $\Delta \cup \Delta^g$ point-wise, which yields that $\Delta^g \subseteq \Delta$ and $g \in K$.

(b)
%$|\Delta|=3$ in this situation, so
%$|\Omega| \equiv 3$ modulo $|P|$. In particular $|\Omega|$ is divisible by $3$.
Suppose that $y \in K \setminus D$ is a $3$-element, and again suppose that $y \in P$. Then $y$ does not fix any element outside $\Delta$.
If it fixes any point in $\Delta$, then it fixes all of them (because $|\Delta|=3$), contrary to the fact that $y \notin D$. Therefore $y$ does not fix any point of $\Omega$.

Next suppose that $x \in D$ is a 3-element. Without loss $x \in P$, and then $x$ fixes the three points of $\Delta$ and has no fixed points outside $\Delta$.

(c) By hypothesis $|\Delta|=3$,  and $P$ acts semi-regularly on $\Omega\setminus \Delta$. Therefore $|\Omega| \equiv 3$ modulo $|P|$. In Case (e) of Lemma \ref{syl3} we have that $P$ fixes a point and stabilizes $|\Delta|$, which means that it fixes every element of $\Delta$.

(d) The additional hypothesis contradicts Case (d) of Lemma \ref{syl3}, and hence Case~(e) must hold.
Next suppose that $3$ divides $|K/D|$. Then there exists a $3$-element that stabilizes $\Delta$, but does not fix every element of it, and this is only possible if $|\Delta| \ge 3$. Hence $|\Delta| = 4$ as stated.

(e) follows from all this: If $|\Delta|=3$, then $3$-elements in $D$ have three fixed points on $\Omega$ and $3$-elements of $K \setminus D$ have none.
Otherwise (d) yields that $3$-elements in $K\setminus D$, if any exist, fix a unique point, whereas $D$ fixes four points in this situation.
\end{proof}

\begin{lemma} \label{durch}
Suppose that Hypothesis \ref{hyp3} holds and that $g \in G \setminus K$. Then
$|K \cap K^g|_3 \le 3$.
\end{lemma}

\begin{proof}
Let  $T \in \Syl_3(K \cap K^g)$.
If $T=1$, then there is nothing left to prove.
Hence we suppose that $T \neq 1 $, and  let $T_0 = D \cap T$. Then %$T_0^g \leq D^g $ and on the other hand $T_0^g \leq K $.
$T_0^{g^{-1}} \leq K $.
Now Lemma~\ref{littlehelp}\,(e) implies that %$T_0^g \leq D$.
$T_0^{g^{-1}} \leq D$.
Thus %$T_0^g \leq D \cap D^g$
$T_0 \leq D \cap D^g$
and Lemma~\ref{littlehelp}\,(a) yields that %$T_0 = T_0^g = 1$,
$T_0 = 1$,
which is our assertion.
\end{proof}

\begin{cor} \label{untergruppe}
Suppose that Hypothesis \ref{hyp3} holds and that $R \le P$ is a subgroup of order at least $9$. Then $N_G(R) \le K$.
\end{cor}

\begin{proof}
We recall that $R \le P \le K$ and
let $g \in N_G(R)$. Now $R=R^g \le K \cap K^g$ and $|R| \ge 9$, and then Lemma \ref{durch} forces $g \in K$.
\end{proof}

\begin{rem}\label{remstrong}
Suppose that Hypothesis \ref{hyp3} holds. Then $P \le K$ and therefore, if
Case (e) of Lemma~\ref{syl3} holds and if there are no
$3$-elements in $K \setminus D$, then all $3$-elements of $K$ lie in $D$ and hence
$K \cap K^g$ is a $3'$-group. This actually means that $K$ is strongly $3$-embedded in $G$ in this case.
\end{rem}

For all prime numbers $p$ and all $p$-subgroups $Y$ of $G$,
we denote by $r(Y)$ the \textbf{rank of $Y$}.
%Zweifel: Ich glaube den Rang brauchen wir hier auch noch nicht.
Moreover, we need notation for the fusion arguments that will come up:

\begin{rem}\label{fusionnota}
Whenever we use Alperin-Goldschmidt's Fusion theorem, then we use the version from \cite[Theorem 16.1]{GLS2}. We state it here for convenience and because we want to refer to its notation later.
Given a finite group $G$ and $P \in \syl_3(G)$, we set
$\mathcal{F}$ to be the set of all nontrivial subgroups $F \le P$ such that the following conditions are satisfied:

\begin{itemize}
\item[(a)] $N_P(F) \in \syl_3(N_G(F))$,

\item[(b)] $C_G(F) \leq O_{3',3}(N_G(F))$, and

\item[(c)] $O_{3',3}(N_G(F))=O_{3'}(N_G(F)) \times F$.
\end{itemize}

Then fusion of elements or subgroups of $P$ can be described via $\mathcal{F}$ as follows:
If $A,B \subseteq P$ and $g \in G$ are such that $A^g=B$, then
we find $n \in \N$, $A=A_1,A_2,...,A_n=B \subseteq P$, $c \in C_G(A)$,
$F_1,...,F_{n-1} \in \mathcal{F}$ and $g_1\in N_G(F_1),...,g_{n-1} \in N_G(F_{n-1})$ such that

(i) $g=cg_1 \cdots g_{n-1}$,

(ii) $\langle A_i,A_{i+1}\rangle \le F_i$ for all $i \in \{1,...,n-1\}$ and

(iii) $A_i^{g_i}=A_{i+1}$ for all $i \in \{1,..,n-1\}$.

We notice that (b) and (c) imply that $Z(P) \leq F$ for all $F \in \mathcal{F}$.
\end{rem}

\begin{lemma}\label{fusion-chris}
Suppose that $p$ is prime, that $P \in \syl_p(G)$ and $P \le K \le G$.
Using the notation from Remark~\ref{fusionnota}, suppose for all $F \in \mathcal{F}$ that $N_G(F) \le K$.
Then $K$ control fusion in $P$.
\end{lemma}

\begin{proof}
Let $a,b \in P$ and $g \in G$ be such that $a^g=b$.
Remark \ref{fusionnota}\,(i), with $\{a\}$ in place of $A$ and $\{b\}$ in place of $B$, says that
$g=c \cdot g_1 \cdots g_{n-1}$. Set $k:=g_1 \cdots g_{n-1}$.
Then $a^c=a$ and therefore $a^g=a^{ck}=a^k$.
Since $g_1 \in N_G(F_1) \le K$,..., $g_{n-1} \in N_G(F_{n-1}) \le K$ by hypothesis, we deduce that
$a^k=b$ and $k \in K$.
\end{proof}

\begin{lemma} \label{fusion}
Suppose that Hypothesis \ref{hyp3} holds. Then
$K$ controls fusion in $P$ with respect to $G$.
\end{lemma}

\begin{proof}
We use the notation from Remark~\ref{fusionnota} and we let $F \in \mathcal{F}$.
Then we prove that $N_G(F) \le K$ in order to apply Lemma \ref{fusion-chris}.

Notice that $|P| \geq 9$ and $|\Delta| \leq 4$, whence it follows that $Q=P \cap D \neq 1$. Since $D \unlhd K$ and $P \le K$, we have that
$Q \unlhd P$ and hence $Z(P) \cap D \neq 1$.
Moreover, Remark~\ref{fusionnota} yields that  $Z(P) \cap D  \leq  F \cap D$.
If $F \not\leq D$, then $|F| \geq 9$ and $N_G(F) \leq K$ by Corollary~\ref{untergruppe}.
If $F \leq D$, then $N_G(F) \leq K$ by Lemma~\ref{littlehelp}\,(a).
Now Lemma \ref{fusion-chris} gives the result.
\end{proof}

%In the proof of the next result we refer to \cite{GLS2} and to the definition of extremal $p$-elements. Therefore we recall this notion briefly, with $G$ and $P$ as before: An element $y \in P$ is said  to be {\bf extremal} in $G$ if and only if $C_P(y) \in \Syl_3(C_G(y))$.

\begin{prop}\label{3semireg}
Suppose that Hypothesis \ref{hyp3} holds.
Then one of the
following is true:
\begin{enumerate}
\item[(a)] $K$ is strongly $3$-embedded in $G$.
\item[(b)] $G$ contains a normal subgroup $N$ of index $3$ such that
$N \cap K$ is strongly $3$-embedded in $N$.
\item[(c)] The point stabilizers in $G$ have even order.
\end{enumerate}
\end{prop}

 \begin{proof}
We suppose that point stabilizers in $G$ have odd order.
If, for all $g \in G\setminus K$, it is true that $K \cap K^g$ is a $3'$-group, then $K$ is strongly $3$-embedded in $G$, and (a) holds.
Next we suppose that (a) does not hold, and our objective is to prove (b).
First we show the following statement:

($\ast$) ~~If $g \in G\setminus K$ is such that $|K \cap K^g|_3 > 1$, and if $x \in K \cap K^g$ has order 3, then
$x^G \cap Q=\varnothing$ and $x^G \cap Qx^2 = \varnothing$.

Let $g \in G \setminus  K$ and assume
that $x \in P^g \cap K$ is nontrivial. If $x \in Q$, then $\FO(x) = \Delta$. Hence
$\Delta = \Delta^g$ is the union of all non-regular orbits of $P^g$. It follows that $g \in K$, which is a
contradiction. We conclude that $\langle x \rangle = P^g \cap K$ complements $Q$ in $P$ (and in $P^g$), and this implies ($\ast$).

Now we see that $x^G \cap P \subseteq Qx$.
This shows that \cite[Proposition~15.15]{GLS2} is applicable with $x$ in place of $u$ and that, in its statement, Property (i) cannot be satisfied. Hence there is a normal subgroup $N$ of $G$ of index $3$
such that $G = NX$.

%Entwurf P:
Furthermore the proof of \cite[Proposition~15.15]{GLS2} shows that \(Q\leq N\) and therefore \(Q=N\cap P\).

Assume that there exist \(h\in N\setminus K\) and some nontrivial \(3\)-element \(y\in N\cap K\cap K^h\leq K\cap K^h\).
If \(y\in D\), then Lemma~\ref{littlehelp}\,(e) yields that also \(y^{h^{-1}}\in D\). But this contradicts Lemma~\ref{littlehelp}\,(a) because \(h\in G\setminus K\).
Therefore \(y\in K\setminus D\).
On the other hand \(y^G\cap P\subseteq N^G\cap P=N\cap P=Q\). Since \(P\in \Syl_3(K)\), there exists \(g \in K\) such that \(y^g\in P\).
Hence \(y^g\in Q=D\cap P\), but \(y^g\) has the same number of fixed points as \(y\in K\setminus D\), and this contradicts Lemma~\ref{littlehelp}\,(e).
This final contradiction proves the lemma.
\end{proof}

\subsection{Case (c) of Lemma~\ref{syl3}}

Now we study Case (c) of Lemma~\ref{syl3} more closely.
Given that $3$-groups of maximal class play a role here, we need the following technical lemma.
Recall that a $3$-group of order $3^n$, $n \in \N$, is said to be of \textbf{maximal class} if and only if it has nilpotency class $n-1$.

\begin{lemma}\label{maxclass3}
Let $P$ be a $3$-group of maximal class and let $n \ge 4$ be such that $|P|=3^n$.
Let
 $1 = Z_0 \leq Z_1 = Z(P) \leq \cdots \leq Z_{n-1}= P$ be the upper central series of $P$, which means that
 for all $i \in \{1,...,n-2\}$,
 $Z_i/Z_{i-1}= Z(P/Z_{i-1})$.
 Moreover let $U \le P$ be a subgroup of index $3$, and  set $P_1 := C_P(Z_{n-2}/Z_{n-4})$.
 Then the following hold:
 \begin{itemize}
 \item[(a)] If $i \in \{1,...,n-2\}$, then $|Z_i/Z_{i-1}| = 3$, and $Z_{n-2} = P^\prime$ has index $9$ in $P$.
 \item[(b)] $P_1 = C_P(Z_2)$  and $|P:P_1| = 3$. Furthermore, if $i \in \{1,...,n-2\}$, then $Z_i \le U$.
 \item[(c)]  Let $y\in P\setminus{U}$ be such that $|C_P(y)|=9$. Then $y^U = P^\prime y$.
\item[(d)]   If $U$ is abelian, then $U = P_1$.
\item[(e)]  Suppose that $U$ is elementary abelian. Then $n = 4$.
\item[(f)] $\aut(P)$ is a $\{2,3\}$-group.
\item[(g)] Suppose that $h \in \aut(U)$ has prime order $r \geq 5$.  Then $U = P_1$, $n =4$, $U$ is elementary abelian and $r= 13$.

 \end{itemize}
 \end{lemma}

 \begin{proof}  Statement (a) follows from  \cite[III Hilfssatz 14.2]{Huppert-I}, while the first part of (b) follows from   \cite[III Hauptsatz 14.7 and Satz 14.14]{Huppert-I}.
 %Zweifel: 14.7. hat n>4 als Voraussetzung, aber ich denke, man braucht stattdessen nur die Definition von Ausnahmegruppe
  If $U = P_1$, then  $Z_{n-2} \leq U$, and if
  $U \neq P_1$, then $P_1 \cap U$ is a normal subgroup of $P$ of index $9$,
  and therefore $P_1 \cap U = Z_{n-2}$ by  \cite[III Hilfssatz 14.2]{Huppert-I}, which gives  the second part of (b).

  In (c) we first observe that $P=\langle y \rangle U$ and therefore $y^P = y^U$.
   The length of this conjugacy class is $|P:C_P(y)| =|P|/9 = |P^\prime|$, which means that $|[U,y]|=|y^U|=|P^\prime|$. Moreover $[U, y] \le P^\prime$, and then equality holds.
   Now $P'=[U,y]$ and therefore $P'y=yP'=y^U$, which is (c).

The definition of $P_1$ and (b) yield (d).

Now we turn to (e).  Since $U$ is elementary abelian, we have that $U = P_1$ by (d). Let $y \in P\setminus{U}$ and $y_1 \in U\setminus{Z_{n-2}}$, and for all $i \in \{2,\ldots,n\}$ set
$y_i := [y,y_{i-1}]$. Then  $P = U \langle y \rangle$, and
for all \(i\in \{2,\ldots, n\}\), we see that
$Z_{n-i} = \langle y_i \rangle Z_{n-i-1}$ by \cite[Lemmas~ 3.2.4, 3.2.7]{LGM} and \cite[III Satz 14.17]{Huppert-I}.
Thus we may suppose that $\{y_1, \ldots , y_{n-1}\}$ is a generating set for $U$.
 Therefore, $y$ induces an automorphism $\tilde{y}$ on $U$ of order  $o(\tilde{y}) = 3^j \geq n-1 > 3^{j-1}$, where $j \in \{1,...,n\}$.
 Since $|P:U| = 3$ and $U$ is abelian, the latter automorphism has order $3$ and $n = 4$, which is (e).

We notice that (f) follows from coprime action and the fact that every automorphism of $P$ stabilizes
its upper central series, with factors of order 3 or 9. Here we do not even need the hypothesis that $n \ge 4$.

For the proof of (g), we
%Zweifel: Wir brauchen bereits fpr 3.3.6 in LGM n>4, müssen das also schon hier voraussetzen.
assume first that $U \neq P_1$. Then by \cite[Corollary 3.3.6]{LGM} we know that $U$ has maximal class, and we have just seen that
$\aut(U)$ is a $\{2,3\}$-group then, contrary to our choice of $h$. Hence $U=P_1$.
Now $\aut(U)$ stabilizes the characteristic subgroup $\Omega_1(U)$.
If $n \geq 5$, then $\Omega_1(U)$ is elementary abelian of order $9$ by \cite[III Satz 14.16]{Huppert-I}, which
yields (f) in this case.
If $n \le 4$, then we must have that $n=4$ and $|U| = 3^3$. Since $\GL_2(3)$ is a $\{2,3\}$-group, it follows that $U$ is elementary abelian and $r = 13$.
\end{proof}

For the remainder of this subsection we work with the following hypothesis.

\begin{hyp}\label{hyp4}
	In addition to Hypothesis \ref{4fix}, let $P \in \syl_3(G)$, let $\Delta$ denote the union of the $P$-orbits of length at most $3$ and suppose that $P$ satisfies Lemma \ref{syl3}~(c), with all the notation given there.
	We define $D$ to be the element-wise stabilizer of $\Delta$ in $G$, $K$ to be the set-wise
	stabilizer of $\Delta$ in $G$, and we set $Q:=P\cap D$. Let $\Lambda$ be a non-regular $P$-orbit on $\Omega \setminus \Delta$ and $\lambda \in \Lambda$.
	Furthermore, we denote by $x$ an element that generates $P_\lambda$, and we keep the notation from Lemma \ref{maxclass3} for the upper central series of $P$, including $P_1:=C_P(Z_{n-2}/Z_{n-4})$.
\end{hyp}

\begin{lemma}\label{tiny}
Suppose that Hypothesis \ref{hyp4} holds. Then
$P$ fixes at most one point in $\Delta$ and
$|\Delta| \in \{0, 1,3,4\}$. Moreover $|P| \geq 27$, $|Z(P)| = 3$, $C_P(P_\lambda) = Z(P) \times P_\lambda$,
$Z(P) \le Q$ and $P_\lambda \nleq P_1$.
Finally, $Z(P)$ has no fixed points outside $\Delta$ and, if $|\Delta| \in \{3,4\}$, then $|P:Q|=3$.
\end{lemma}

\begin{proof}
We know that $|\Delta| \leq 4$ by hypothesis.
If $P$ fixes two or more points on $\Delta$, then Lemma \ref{syl3}~(c) gives that $P_\lambda$ fixes five or more elements in $\Omega$, which contradicts Hypothesis~\ref{4fix}.
In particular this implies that $|\Delta| \neq 2$, which
proves the first two statements. If $|\Delta| \in \{3,4\}$, then $P \neq Q$ and hence $|P:Q|=3$, which is the last statement.
Since $|\Lambda| = |P:P_\lambda|=|P|/3 >3$, it also follows that $|P| \geq 27$.
Now we recall that $P$ is of maximal class and therefore $|Z(P)| = 3$.
Moreover $Q \unlhd P$ and therefore $Z(P) \cap Q \neq 1$, which together with the fact that $|Z(P) | = 3$ forces $Z(P) \leq Q$.

Assume that $Z(P)$ fixes a point outside of $\Delta$. Then
it fixes three points outside of $\Delta$ in total, and $P$ stabilizes the set of these fixed points,
which means that they give a $P$-orbit of size $3$ outside of $\Delta$. This is impossible.

Next we note that \(C_P(P_\lambda)\) stabilizes \(\FO(P_\lambda)\) (of size 3), which means that it induces a group of order at most \(3\) on \(\FO(P_\lambda)\). This implies that \(|C_P(P_\lambda)|\leq 9\), hence $|C_P(P_\lambda)|=9$ and
%Since \(Z(P)\neq P_\lambda\), because otherwise the \(P\)-orbit of \(\lambda\) would be \(\FO(Z(P))=\FO(P_{\lambda})\),
%it follows $C_P(P_\lambda) = Z(P) \times P_\lambda$.
$C_P(P_\lambda) = Z(P) \times P_\lambda$.
If we assume that  $P_\lambda \le P_1$, then $[P_\lambda,Z_2]=1$ which, by the previous statement, forces
$Z_2=C_P(P_\lambda)$ and then $P_1\le C_P(P_\lambda)$. This is false, and therefore $P_\lambda \nleq P_1$ as stated.
\end{proof}

\begin{rem}\label{ConjugacyPlambda}
  Since $ [P_\lambda, Z_2] = Z_1 = Z(P)$, all the subgroups of order $3$ in $F:= Z(P) \times P_\lambda$ that are different from $Z(P)$
are conjugate in $N_P(F)$.
If $|\Delta| = 0$, then $|\FO(P_\lambda)| = 3$ and therefore  $Z(P)$ and $P_\lambda$ are not conjugate in $G$.
If $|\Delta| = 1$, then $|\FO(P_\lambda)| = 4$ and again  $Z(P)$ and $P_\lambda$ are not conjugate in $G$.
\end{rem}

We can now draw some conclusions.

\begin{lemma}\label{OddNormalizersOrder9}
Suppose that Hypothesis \ref{hyp4} holds.
	Let $F = Z(P) \times P_\lambda$, and suppose that the point stabilizers have odd order.
\begin{itemize}
 \item[(a)] If $|\Delta| \in \{0,1,3\}$, then $N_G(F)= N_P(F)C_G(F) $ and $N_G(F)$ has two orbits of length
	$1$ and $3$, respectively, on the set of subgroups of order $3$ of $F$ (see Remark~\ref{ConjugacyPlambda}).
 \item[(b)] If $|\Delta| \in \{1,4\}$, then $N_G(F)$ has odd order.
\end{itemize}
\end{lemma}
\begin{proof}
If 	$N_P(F)C_G(F) $ is a proper subgroup of $N_G(F)$, then $N_G(F)$ acts transitively on the set of subgroups of order $3$ of $F$ and $|\Delta| = 3$ by Remark~\ref{ConjugacyPlambda}.
Then $4$ divides $|N_G(F)/C_G(F)|$, which implies that $N_G(F)$ induces the central involution of $\aut(F) \cong \GL_2(3)$ on $F$. This involution inverts every element
in $F$, in particular it stabilizes the fixed point set of $P_\lambda$ and therefore it fixes a point. This contradicts our hypothesis that point stabilizers have odd order. we conclude that (a) holds.

If $|\Delta|  \in \{1,4\}$, then $|\FO(F)| =1$ and we conclude (b).
\end{proof}

\begin{lemma}\label{OddNormalizers}
Suppose that
Hypothesis~\ref{hyp4} holds, that the point stabilizers have odd order and that
$U \in \{P_\lambda,  Z(P)\}$. Then the following hold:
\begin{itemize}
\item[(a)] If $|\FO(U)| = 4$, then $N_G(U)/C_G(U)$ has odd order.
\item[(b)] If $|\FO(U)|  \in \{1,3\}$, then $N_G(U)$ has odd order.
\end{itemize}
\end{lemma}
\begin{proof}
We know that $N_G(U)$ stabilizes the set $\FO(U)$, and then assertion (b) follows.
Next suppose that $|\FO(U)| = 4$. We assume for a contradiction that  $|N_G(U)/C_G(U)|$  is even and we let
 $g \in N_G(U)\setminus C_G(U)$ be a nontrivial $2$-element. Then $g$ inverts the elements in $U$.
Let $T \leq P$ be such that $\{T,U\}= \{Z(P), P_\lambda \}$. Then $T$ and $g$ act nontrivially on $\FO(U)$,
in fact $\langle T,g\rangle$ induces $\Alt(4)$ on $\FO(U)$ by Lemma~\ref{Inv-lem:TB}. Thus, if we set
$V := \langle g^T\rangle$, then $VT \cong \Alt(4)$.
But we also know that $[T,U] = 1$, and then $g \in V = [V,T] \leq C_G(U)$, in contradiction to the choice of $g$.
\end{proof}

\begin{lemma} \label{NormalIndex3Subgroup}
Suppose that
Hypothesis~\ref{hyp4} holds, that $|P| \ge 3^4$, and that the point stabilizers in $G$ have odd order. Then $P \cong C_3 \wr C_3$ or $G$ has a normal subgroup of index $3$.
\end{lemma}

\begin{proof} Suppose that $G$ has no normal subgroup of index $3$. Notice that if $P$ has a quotient isomorphic to $C_3 \wr C_3$, then $ P \cong C_3 \wr C_3$ by \cite[III,~Theorem~14.20]{Huppert-I}. Hence we may suppose that $P$ does not have any quotient
isomorphic to $C_3 \wr C_3$. Now Yoshida's Transfer Theorem \cite[Theorem~15.19]{GLS2} implies that $[P,N_G(P)] = P$. Since $\aut(P)$ is a $\{2,3\}$-group, there exists a $2$-element $g \in N_G(P)$ such that
$P/P^\prime$ is inverted by $g$. Let $\lambda$ be as in Hypothesis~\ref{hyp4}.
Then $g$ normalizes $P^\prime P_\lambda$. We recall that all subgroups of order $3$ in
$P^\prime P_\lambda \setminus{P^\prime}$ are $P$-conjugate, by Lemma~\ref{maxclass3}
(c), and then a Frattini argument gives that $N_{P\langle g\rangle}(P_\lambda)P = P \langle g \rangle$. Since $g$ inverts $P/P^\prime$, it follows that some $2$-element in $N_{P\langle g\rangle}(P_\lambda)$ inverts $P_\lambda$. Moreover, $N_{P\langle g\rangle}(P_\lambda)$ stabilizes the set $\FO(P_\lambda)$, which has size $3$, and consequently $G_\lambda$ has even order. This is a contradiction.
\end{proof}

\begin{lemma}\label{3:xNotInN}
Suppose that
Hypothesis~\ref{hyp4} holds, that $G$ has a normal subgroup $N$ of index $3$ and that $\Delta = \varnothing$. Moreover suppose that the point stabilizers have odd
order. Then
$N$ does not contain a subgroup $P_\mu$, where $\mu \in \Omega$.
In particular, $N \cap P$ acts semi-regularly on  $\Omega$.
\end{lemma}

\begin{proof}
By hypothesis $\Delta = \varnothing$, which yields that $P_\mu$ is a Sylow $3$-subgroup of $G_\mu$ and that all these subgroups are conjugate in $G$. Let $H \leq G_\mu$  be a nontrivial stabilizer of four points in $\Omega$. By hypothesis and by Lemma~\ref{Inv-lem:TB} there is a prime $r \in \pi(H)$ such that $r \geq 5$, and a Sylow $r$-subgroup $R$ of $G$ that is a subgroup of $H$. By Frattini we know that $G = NN_G(R)$. Therefore there is a $3$-element $y$ in $N_G(R)\setminus{N_N(R)}$, and $y$ stabilizes the set $\FO(R)$, which has size 4. This means first that $y$ fixes a point in $\Omega$ and then
that one and therefore  all the conjugates of $P_\lambda$ intersect $N$ trivially.
\end{proof}

\begin{lemma}\label{CaseDelta=3}
Suppose that
Hypothesis~\ref{hyp4} holds, that
$|\Delta| = 3$ and that the point stabilizers have odd order. Then $P \cong C_3 \wr C_3$.
In addition, if $\delta \in \Delta$, then  $G_\delta$ is a Frobenius group where the complements
are cyclic of order $13$, and the Sylow 3-subgroups of $G_\delta$ are elementary abelian of order 27.
\end{lemma}

\begin{proof}
    Let $\alpha \in \Delta$ and first note that $P \cap G_\alpha = Q = P \cap D $.
    Let $H$ be a subgroup of $G_\alpha$ such that $|\FO(H)| = 4$, and let $r \in \pi(H)$. Since elements of order $3$ fix either zero or three elements in $\Omega$, we see that $r \geq 5$.
We recall that $G_\alpha$ has odd order and is, therefore, soluble. Now one of the following holds: There is a nontrivial $r$-subgroup $R$ of $H \le G_\alpha$ that is normalized by some Sylow $3$-subgroup of $G_\alpha$, or there is a nontrivial $3$-subgroup  of $G_\alpha$ that is normalized by some Sylow $r$-subgroup of $H$. Here we use Lemma~\ref{helpFrattini}.
In the first case, we may suppose that $Q$ normalizes $R$, and therefore it stabilizes $\FO(R) = \FO(H)$, a set of size 4. Now, since $|Q| \geq 3^2$, there is a nontrivial element in $Q$ which acts trivially on $\FO(R)$, and this is impossible. Therefore we find a nontrivial $3$-subgroup $P_0$ of $G_\alpha$ that is normalized by some Sylow $r$-subgroup $R$ of $H$.
We may suppose that $P_0 \le Q$.  Then $R$ stabilizes the set $\FO(P_0) = \Delta$, which implies that $\Delta \subseteq \FO(R)$.
If $P_0 \neq Q$, then there is a nontrivial $3$-subgroup $T$ of $G_\alpha$ that normalizes $P_0 \cdot R$, with a Frattini argument, because $G_\alpha$ is soluble.
Therefore $T$ stabilizes $\FO(P_0)$ and  $\FO(R)$. Then the fact that  $\alpha \in \Delta$ implies that $T$ acts trivially on $\FO(H)$, which is false. Therefore $P_0 = Q$.
This shows, in particular, that $QH$ is a Frobenius group.
If $|P| = 3^3$, then $Q$ has order $9$ and therefore it does not have an automorphism of order $r$. This forces $|P| \geq 3^4$. Now our assertion follows from Lemma~\ref{maxclass3} (g), because $Q$ is elementary abelian and then $P \cong C_3 \wr C_3$.
\end{proof}

\begin{prop} \label{maxKlasse}
Suppose that Hypothesis~\ref{hyp4} holds.
Then one of the following is true:

\begin{enumerate}

\item[(a)] The point stabilizers have even order.% All $P$-orbits have length at least $\frac{|P|}{3}$.

\item[(b)]  $|\Delta| = 0 $,  and $G$ contains a normal subgroup $N$ of index $3$ or $9$ such that $|N_\lambda|$ is coprime to $6$.
%$K$ is strongly $3$-embedded in $G$.

\item[(c)] $|\Delta| = 1$ and $G$ contains a subnormal subgroup $N$ of index $3$ or $9$ such that $N \cap K$ is strongly $3$-embedded in $N$.

\item[(d)] $|\Delta| = 4$  and $G$ contains a normal subgroup $N$ of index $3$ or $9$ such that $N \cap K$ is strongly $3$-embedded in $N$.

\item[(e)]  $|\Delta| = 0$, $P$ is extra-special of exponent $3$ and order $27$, and $|G_\lambda|_3=3$.

\item[(f)] $|\Delta| = 3$ and $P \cong C_3 \wr C_3$. If $\delta \in \Delta$, then  $G_\delta = D$ is a Frobenius group with Frobenius complements
of order $13$, and the Sylow 3-subgroups of $G_\delta$ are elementary abelian of order 27.
\end{enumerate}
\end{prop}

\begin{proof}
We suppose that the point stabilizers have odd order.
If $|\Delta| = 3$, then (f) holds by Lemma~\ref{CaseDelta=3}.
Therefore, and by Lemma \ref{tiny},
we suppose from now on that $|\Delta| \in \{0,1,4\}$.
We first consider the cases where
$P \cong 3^{1+2}$ or $P \cong C_3\wr C_3$.
If $|\Delta| = 0$, $P \cong 3^{1+2}$ and $P$ has exponent $3$, then
(e) holds. Thus we exclude this case now.

Our strategy is, as in the previous subsection, to understand the fusion of $P_\lambda$ in $P$, which is why we analyze a
Alperin-Goldschmidt conjugating family, as explained in Remark~\ref{fusionnota}.

%in order to produce a normal subgroup $N$ of $G$ of index $3$ or $9$ such that (b), (c) or (d) is satisfied. More precisely we will consider the normalizers of the elements $F$ in the fusion family ${\cal F}$, see Remark~\ref{fusionnota}, that contain $P_\lambda = \langle x \rangle$ in order to understand the fusion of $P_\lambda$ in $P$. All the remaining cases will then be handled with the help of Lemma~\ref{NormalIndex3Subgroup}.
\medskip

\textbf{Case 1: (i) $P \cong 3^{1+2}$, $P$ has exponent $9$ and $|\Delta| =0$, or \\
(ii) $P \cong C_3 \wr C_3$.}

Let $F \in {\cal F}$ be such that $P_\lambda = \langle x \rangle \le F$ (see Remark~\ref{fusionnota}).

First we suppose that there is $y \in N_G(F)$ such that $P^\prime x^y \neq P^\prime x$.
%Then in particular, the action of $y$ on $F$ is not induced by an element in $P$.
\medskip

If $|F|= 3$, then we recall that $\langle x \rangle \le F$ and that $x^y \neq x$, which means that
$y$ inverts $x$. In particular $y$ has even order, which contradicts
Lemma~\ref{OddNormalizers}.
\medskip

Suppose that $|F|= 9$. If $F$ is cyclic, then $y$ normalizes its unique subgroup of order 3, which is $\langle x \rangle$, and we obtain a contradiction as in the case where $|F| = 3$. Next suppose that $F$ is not cyclic. Then $F = Z(P) \times P_\lambda$ and $\aut(F) \cong \GL_2(3)$ is a
$\{2,3\}$-group.
Together with Lemma~\ref{OddNormalizersOrder9}, our hypotheses yield that $y$ is a $3$-element. Since $Z_2(P)$ induces a Sylow $3$-subgroup of $\aut(F)$ on $F$, it follows that $y$ acts like an element of $Z_2(P)$ on $F$, which contradicts the choice of $y$.

\medskip

Suppose that $|F| = 3^3$ and that $F \lneq P$. Then $|P| >3^3$ and hence Case 1\,(ii) holds, i.e. $P \cong C_3 \wr C_3$. Now $P^\prime \cong C_3^2$
and $F = P^\prime P_\lambda$ by Lemma~\ref{maxclass3}.
We can see that $F$ is an extra-special group of exponent $9$, and then it follows that $y$ is a $3$-element by \cite[Theorem~1 (b)]{Winter}.
Thus $F \langle y \rangle  = P$, contrary to our choice of $y$.
\medskip

Next suppose that $|F|=3^3$ and that $F = P$. Then we may suppose that $y$ is a $2$-element by Lemma~\ref{maxclass3} (f) and by \cite{Winter}. Since  $y$  normalizes $F$, hence $P$ and then $Z(P)$, it stabilizes $\Delta$. Then it follows that $|\Delta| \notin  \{1,3\}$. We conclude that $|\Delta| = 0$ or $|\Delta| = 4$ in this case.
The same is true if $|F|=3^4$ and $F = P$: Again $y$ acts like an element from $P$ (which is impossible) or it can be chosen as a $2$-element (and then $|\Delta| = 0$ or $|\Delta| = 4$).

 \medskip
We summarize the situation for the following arguments:

In Case 1 (i) or (ii), we have a $2$-element $y \in N_G(F)$ such that $P^\prime x^y \neq P^\prime x$, and $F=P$.
In particular $y \in N_G(P)$.

First we consider the case where $|\Delta| = 4$. We recall that $\langle x \rangle=P_\lambda$, and now $\langle y,x \rangle$ induces a group $B \cong \Alt(4)$ on $\Delta$.
Let $V = \langle y, y^x\rangle$.
Since $y \in N_G(P)$, this yields that $[x,y] \in P \cap DV = P \cap D$, and $V$ induces $O_2(B)$ on $\Delta$. But this is impossible because $B \cong \Alt(4)$ and $P \cap D$ is the point-wise stabilizer of $\Delta$ in $P$.

Now $|\Delta| = 0$, and still we have two possible sub-cases for $P$.

First suppose that $P \cong 3^{1+2}$ and that $P$ has exponent $9$ (Case 1\,(i)).
Then $P = TP_\lambda$, where $T$ is a cyclic group of order $9$. We recall that $P=F$ and that $y$ is a $2$-element that induces a non-trivial automorphism on $P$.
% and that, therefore, $|N_G(P)/C_G(P)|$ is even. Then we apply
In \cite[Theorem~1]{Winter} we see the structure of $\aut(P)$ and in particular that it has a unique class of involutions, and therefore we may suppose that $y$
 inverts $T$ and centralizes $x$. But this  contradicts our choice of $y$.
We conclude that, in Case 1\,(i), all $y \in N_G(P)$ have the property that
$P^\prime x^y = P^\prime x$. This contradicts Remark \ref{ConjugacyPlambda}.

Now suppose, still in the situation summarized above, that $P \cong C_3 \wr C_3$.
Then $P$ has a characteristic subgroup of structure $C_3^3$, which must be $P_1$ (see Lemma \ref{maxclass3}(d)), and $P_1$ is normalized by $y$.
Moreover $y$ acts nontrivially on the section $P/P^\prime \cong C_3^2$.
 %As the square of every element of order $4$ in
 %$GL_2(3)$ is its central involution, $y$ is of order $2$ and $|N_G(P):P| = 2$.
 If $P^\prime x^y = P^\prime x^2$, then
 by Lemma~\ref{maxclass3} (c) we may suppose that $y$
 inverts $x$. This contradicts  Lemma~\ref{OddNormalizers}.
Recalling the choice of $y$, we find $v \in P \setminus{P^\prime P_\lambda}$ such that $P^\prime x^y = P^\prime v$.
Set $R:= P^\prime \langle vx^2 \rangle $. Then $|P:R|=3$, and moreover
 $R \cap x^G = \varnothing$ and $x^G \cap P \subseteq Rx$  by Alperin-Goldschmidt \cite[Theorem~16.1]{GLS2}. We will come back to this.
\medskip

In order to finish Case 1, we now
suppose that, in the case where $P \cong C_3 \wr C_3$, such an element $y$ as taken at the beginning of Case 1 does not exists.
If $|\Delta| = 4 $, then we set $R: = D \cap P$, and if
 $|\Delta| \in \{0,1\}$, then $R:= P_1$, where $P_1$ is the characteristic subgroup of $P$ of structure $C_3^3$ (see Lemma \ref{maxclass3} and previous paragraph).

 Again it follows that $R \cap x^G = \varnothing$ and that $x^G \cap P \subseteq Rx$ by Alperin-Goldschmidt \cite[Theorem~16.1]{GLS2}.
 \medskip

In all remaining situations of Case 1, we can now apply \cite[Proposition~15.15]{GLS2}. It yields that
there is a normal subgroup $N$ of index $3$ in $G$ such that
 $G = N P_\lambda$ and $N \cap P = R$.
\medskip

Next we show that there does not exist any   $\mu \in \Omega\setminus{\Delta}$ such that
 $P_\mu \le R$.
Assume the contrary. Then $\Delta \neq \varnothing$ by Lemma~\ref{3:xNotInN}. Furthermore, by Lemma~\ref{tiny}, we know that $P_\mu$ is not a subgroup of $P_1$. Also, $P_\mu$ is not a subgroup of $D \cap P$ by our global fixity 4 hypothesis.
It follows that there is some $2$-element $y$ in $N_G(P)$ and some $v \in P \setminus{P^\prime \langle x \rangle}$ such that $P^\prime x^y = P^\prime v$, and
 $R = P^\prime \langle vx^2 \rangle$.
If $|\Delta| =4$, then $R\le D$ does not contain $P_\mu$.
If $|\Delta| =1$, then it is not possible that $R = P^\prime \langle vx^2 \rangle$ has structure $C_3^3$.
Since $|\Delta| \in \{0,1,4\}$ (as seen at the beginning of Case 1), we have a contradiction and therefore, for all $\mu \in \Omega\setminus{\Delta}$, we have that
 $P_\mu \nleq R$.

 The previous paragraph shows that, if $|\Delta| \in \{0,1\}$, then every $3$-element of $N^\#$
  fixes no point or exactly one point, respectively, and therefore (b) or (c) holds, respectively. If $|\Delta| = 4$, then our construction of $N$ implies that $R = N \cap P = D \cap P$.
 Therefore every $3$-element in $N$ fixes $\Delta$ element-wise, and (d) holds.
 \bigskip

\textbf{Case 2: $|P| \geq 3^4$ and $P \not\cong C_3\wr C_3$.}
%If $|\Delta| = 3$, then (f) holds by Lemma~\ref{CaseDelta=3}.

Then Lemma~\ref{NormalIndex3Subgroup} tells us that there is a normal subgroup $N$ of index $3$ in $G$.
Assume that there does not exist $\mu \in \Omega\setminus{\Delta}$
such that $P_\mu \leq N \cap P$. If $|\Delta| \in \{0,1\}$,
 then (b) or (c) holds.
Now let $|\Delta| = 4$. Then $N \cap P$ is as in  Lemma~\ref{syl3} (e), and we conclude (d) from Proposition~\ref{3semireg}.
\medskip

Finally assume that there is some $\mu \in \Omega\setminus{\Delta}$
such that $P_\mu \leq N \cap P$. Then $\Delta \neq \varnothing$ by Lemma~\ref{3:xNotInN}, and $N  \cap P \neq P_1$ by Lemma~\ref{tiny}, and therefore $N \cap P$ has maximal class. (This would also follow from Lemma~\ref{syl3}). If $N \cap P $ has order at least $ ^4$, but is not isomorphic to $C_3 \wr C_3$,
then we apply Lemma~\ref{NormalIndex3Subgroup} to $N$ and we obtain a normal subgroup $N_1$ of index $3$ in $N$. Then $N_1 \cap P = P^\prime$, which yields (b), (c) or (d), because $P^\prime$ does not contain any $P_\mu$ (where $\mu \in \Omega\setminus{\Delta}$) by Lemma~\ref{tiny}.

If $|N \cap P| = 3^3$ or $N \cap P \cong C_3 \wr C_3$, then we
can show exactly as in Case 1 that there is a normal subgroup of index $3$ in $N$ as stated  in (c) or (d). In conclusion, (c) or (d) holds.
\end{proof}

Finally, we can connect the 2-structure and the 3-structure of $G$ and pave the way to the proof of our main results.

\begin{thm}\label{3and2}
Suppose that $G$ is a finite group acting transitively, faithfully and with fixity $4$ on a set~$\Omega$.
Let $P$ be a Sylow $3$-subgroup of $G$, let $\Delta$ be the union of $P$-orbits of length at most $3$ on $\Omega$ and let $K$ denote the stabilizer of the set $\Delta$ in $G$.
Moreover let $\alpha \in \Omega$ and let $f$ denote the maximum number of fixed points of involutions in $G$.
Let %$F:=\langle t^g \mid t,g \in G, o(t)=2,\FO(t)=f\rangle$.
$F:=\langle t \mid t \in G, o(t)=2,\FO(t)=f\rangle$.

 Then one of the following holds:
\begin{enumerate}
\item[(a)]    There is a normal subgroup $N$ of $G$ of index dividing $9$ such that
\begin{enumerate}
\item[(i)]  $(|N_\alpha|, 6) = 1$ or
\item[(ii)]  $K \cap N$ is strongly $3$-embedded in $N$ and $1 \leq |\Delta| \leq 4$.
\end{enumerate}
\item[(b)] $|P| \le 3^4$, more precisely
\begin{enumerate}
\item[(i)]  $|P| = 3$ and $|N_G(P)|_2 \le 4$, or
\item[(ii)] $P$ is  elementary abelian of order $9$, or
\item[(iii)] $P$ is extra-special of exponent $3$ and order $27$, or
\item[(iv)] $P \cong C_3 \wr C_3$.
\end{enumerate}

\item[(c)] $1\le f \le 4$ and $F$ has a strongly embedded subgroup. More precisely, if $S \in \syl_2(G)$, if $\Gamma$ is the union of $S$-orbits of length at most $f$
and if $F_0$ denotes the stabilizer in $F$ of the set $\Gamma$, then $F_0$ is strongly embedded in $F$.
%and contains $N_F(S)$.

\item[(d)] $2\le f \le 3$ and $G$ has dihedral or semi-dihedral Sylow $2$-subgroups.

\item[(e)] $f=4$ and $G$ has sectional $2$-rank at most $4$.
\end{enumerate}
\end{thm}

\begin{proof}
Let $\alpha \in \Omega$ and first suppose that $G_\alpha$ has odd order.
We go through the cases of Lemma \ref{syl3}.

\smallskip
\textbf{Lemma \ref{syl3}\,(a):} Then $G_\alpha$ has order coprime to $6$ and therefore
(a)(i) of the theorem is satisfied in the case where $N=G$.

\smallskip
\textbf{Lemma \ref{syl3}\,(b):}
Then  $|\Delta| > 4$ and $|P| \leq 9$.

If $P=1$, then $G$ is a $3'$-group and therefore the point stabilizers have order coprime to $6$. This is
(a)(i) in our theorem, again, in the case where $N=G$.

If $|P| = 3$, then $P$ is contained in a point stabilizer or the point stabilizers have order coprime to $6$, which is again (a)(i).
Hence we suppose without loss that $P$ fixes $\alpha$.
%Since $|\Omega|=|G:G_\alpha|$ is coprime to $3$ now and all nontrivial $P$-orbits have length $3$, we see that $|\FO(P)| \in \{1,2,4\}$. If $T \le N_G(P)$ is a $2$-subgroup, then it stabilizes $\FO(P)$, a set of size at most $4$.  The elements of $T^\#$ do not fix any points because the point stabilizers have odd order,  and then
The fact that $G_\alpha$ has order yields that $|N_{G_\alpha}(P)|$ is odd, and
Lemma \ref{Inv-lem:TB} gives that $|N_G(P):N_{G_\alpha}(P)| \le 4$. Together this forces $|N_G(P)|_2 \le 4$, as stated in (b)(i) in the theorem.

If $P$ has order $9$, then there are two cases:
First assume that $P$ is cyclic. Then $\Omega_1(P)$ acts trivially on each $P$-orbit of length at most $3$, which means that it acts trivially on $\Delta$. But in the present case
$|\Delta| > 4$, contrary to Hypothesis~\ref{4fix}.
Therefore $P$ is elementary abelian, which is (b)(ii) in our theorem.

\smallskip
\textbf{Lemma \ref{syl3}\,(c):} This situation is treated using Hypothesis \ref{hyp4}, which gives the possibilities in~Proposition \ref{maxKlasse}. \ref{maxKlasse}\,(b) gives Case (a)(i) of our theorem, \ref{maxKlasse}\,(c) and (d) lead to Case (a)(ii) of our theorem and \ref{maxKlasse}\,(e) is contained in Case (b)(iii). Finally, \ref{maxKlasse}\,(f) gives Case (b)(iv) of our theorem.

\smallskip
\textbf{Lemma \ref{syl3}\,(d) and (e):}
This is captured by Hypothesis
\ref{hyp3} and analyzed completely in Proposition \ref{3semireg}. As $G_\alpha$ has odd order, all the possibilities there are included
in Case (a)(ii) of the theorem.

\smallskip
Next we suppose that $G_\alpha$ has even order.

Then $1 \le f \le 4$ and we refer to the main theorem in \cite{Ronse}, which leads directly to the cases (c), (d) or (e) of our theorem. The additional information in (c) comes from Ronse's Proposition 3.1 (see \cite{Ronse}).
\end{proof}

We remark that this result immediately implies Theorem \ref{3total}, as stated in the introduction, because if the point stabilizers have odd order, then (a) or (b) of Theorem \ref{3and2} hold and all the possibilities are captured in Theorem \ref{3total}.\\

\textbf{Proof of Theorem \ref{3and2simple}:}

Here we go through the cases of Theorem \ref{3and2} above with the additional hypothesis that $G$ is non-abelian simple, which means that the normal subgroups $N$ and $F$ mentioned in the statement are $G$ itself.

Case (a)\,(i) means that point stabilizers have order coprime to 6, which is Theorem \ref{3and2simple}\,(5).
Case (a)\,(ii) with a strongly $3$-embedded subgroup is captured by Theorem \ref{3and2simple}\,(4)\,(a).
Case (b) gives exactly the possibilities mentioned in Theorem \ref{3and2simple}\,(4)\,(b)--(e).
Case (c) is Theorem \ref{3and2simple}\,(1), Case (d) is Theorem \ref{3and2simple}\,(2) and Case (e) is Theorem
\ref{3and2simple}\,(3). \qed

\medskip
At the beginning of the next section we will discuss the consequences of this result and how we treat the individual cases.

\section{Strategy and tools}

If we consider Theorem \ref{3and2} in the special case where $G$ is a non-abelian simple group, then
we see a natural case distinction: %TODO: Dafür haben wir mittlerweile einen eigenen Satz

Case (a) leads to point stabilizers of order coprime to 6,
Case (b) leads to small and explicitely described types of Sylow $3$-subgroups
and
Case (c) gives a strongly embedded subgroup.
Here we can apply Bender's main results from \cite{Ben}, which gives specific possibilities
for the group $G$.
In Cases (d) and (e) we also have classification results available, namely
the theorems of Gorenstein and Walter (\cite{GW}), of Alperin, Brauer and Gorenstein (\cite{ABG2}) and of Gorenstein ad Harada (\cite{GoHa1974}), respectively.

We will consider groups of Lie type and small Lie type first, because they appear several times in our analysis.
Then there will be separate sections for the different cases for the point stabilizer order.
In all cases, there are some small groups to consider, and also for some of the sporadic groups we did not really see a convenient generic argument, which is why we decided to analyze some relatively small groups with
%\texttt{GAP}.
GAP~\cite{GAP}.
The code can be found in the Appendix along with examples for output.

The general analysis in later sections also requires some general techniques and basic results for groups with low fixity.

So the remainder of this section is dedicated to results about small groups (with %\texttt{GAP})
GAP)
and to technical background.

\begin{rem}\label{rem4.1}
	The information for the following small examples have been calculated with GAP, %TODO GAP-Code einfügen
	using~\cite{TomLib} and the code in the appendix.

\begin{longtable}{ll}	%HEAD
	$Group$ & possible point stabilizer structure\\[0.5ex]\hline\hline
	\endhead \\[-2ex]
	%MAIN
	\(\Alt_6 \cong \PSL_2(9)\)            & \(C_2\),
	\(\Sym_3\), \(C_3 \times C_3\)\\
	& \(D_{10}\) \text{or}
	\((C_3 \times C_3) : C_2\) \\ \hline
	\(\Alt_7\)                            & \(C_5\) \text{or}
	\(\Alt_6\) \\  \hline
	
	\(\PSL_2(7) \cong \PSL_3(2)\)         & \(C_2\) \text{or}
	\(\Sym_3\) \\ \hline
	
	\(\PSL_2(8)\)                         & \(C_2\), \(\Sym_3\),
	\(D_{14}\) \text{or}
	\(D_{18}\) \\ \hline
	
	\(\PSL_2(11)\)                        & \(C_3\) \text{or}	
	\(\Alt_4\)                                             \\ \hline

	\(\PSL_2(13)\)                        & \(C_3\), \(\Alt_4\) \text{or}
	\(C_{13} : C_3\)                                             \\\hline

	\(\PSL_2(17)\)                        & \(C_4\) \text{or} \(C_{17} : C_4\)                                             \\\hline
	
	\(\PSL_2(19)\)                        & \(C_5\)                                                                        \\\hline
	
	\(\PSL_2(23)\)                        & \(C_6\)                                                                        \\\hline
	
	\(\PSL_2(25)\)                        & \(C_6\) \text{or} \((C_5\times C_5) : C_4\)                        \\\hline
	
	\(\PSL_2(27)\)                        & \(C_7\)                                                                        \\\hline

	\(\PSL_2(29)\)                        & \(C_7\) \text{or} \(C_{29} : C_7\)                                             \\\hline

	\(\PSL_2(31)\)                        & \(C_8\)                                                                        \\\hline
	
	\(\PSL_2(37)\)                        & \(C_9\) \text{or} \(C_{37} : C_9\)                                             \\\hline
	
	\(\PSL_2(41)\)                        & \(C_{10}\) \text{or} \(C_{41} : C_{10}\)                                             \\\hline

	\(\PSU_3(3)\)                         & \(((C_3 \times C_3) :C_3) :C_8\)
	\\\hline

	\(\PSU_4(2)\cong \PSp_4(3)\)          & \(C_5\)                                                                           \\\hline
	\(\PSU_4(3)\)                         & \(C_5\)                                                                           \\\hline
	\(\PSp_4(4)\)                         & \(C_{17}\)                                                                           \\\hline
	\(\PSp_4(5)\)                         & \(C_{13}\)                                                                           \\\hline
	\(\Sz(8)\)                            & \(C_5\) \text{or}
	\(C_{13}\)                                                                                                                \\\hline
	
	\(M_{11}\)                            & \(C_5\), \(C_{11} :C_5\) \text{or}
	\(\PSL_2(11)\)                                          \\\hline
	\(M_{12}\)                            & \(M_{11}\)                                                                       \\
    \hline
	\(M_{22}\)                            & \(C_5\) \text{or}
	\(C_{11}:C_5\)                                                                                                                \\\hline

	\(J_1\)                               & \(C_{15}\)
     \\[0.5ex]\hline\hline\\[-1ex]
	\caption[Small Groups Table for fixity 4]{For each group, we give the possible point stabilizer structures for fixity 4 actions.}
	\label{SmallGroupsTable}

\end{longtable}

Using GAP and the table of marks again (TomLib package, see Appendix with GAP code), we can check that the following groups do not exhibit any fixity 4 actions:

\(\PSL_2(4)\cong \PSL_2(5)\),

	\(\Alt_8\), \(\Alt_9\),   \(\Alt_{10}\),      \(\Alt_{11}\),
	
	\(\PSL_2(16)\),     \(\PSL_2(32)\) , \(\PSL_2(64)\),
	
	\(\PSL_3(3)\),                   \(\PSL_3(4)\),      \(\PSL_3(5)\),%?
	\(\PSL_3(7)\),          \(\PSL_3(11)\),
	
	\(\PSL_4(2)\), \(\PSL_4(3)\),     \(\PSL_5(2)\),
	
	\(\PSU_3(4)\),              \(\PSU_3(5)\),       \(\PSU_3(7)\),%?
	\(\PSU_3(9)\),
	
	\({}^2\textrm{F}_4(2)'\),
	\(\Mat_{23}\),                       \(\Mat_{24}\),     \(\Jan_2\),    \(\Jan_3\),    \(\McL\),       \(\HS\),     \(\He\).
\end{rem}

%%P: Ich glaube wir benutzten das so selten, dass wir das gar nicht benutzen wollen
%\begin{definition}
%Let $G$ be a finite group and let $p,q \in \pi(G)$ be distinct.\\
%We write $p \vdash q$ if and only if one of the following holds:
%
%$p,q \ge 5$ and there exists $1  \neq X \le G$ such that $p \in \pi(X)$ and $q \in \pi(N_G(X))$.
%
%$p \ge 5, q = 2$ and there exists $1  \neq X \le G$ such that $p \in \pi(X)$ and
%$8$ divides $|N_G(X)|$.
%
%$p \ge 5, q = 3$ and there exists $1  \neq X \le G$ such that
%$p \in \pi(X)$ and $9$ divides $|N_G(X)|$.
%
%$p \in \{2,3\}, q \ge 5$ and for all elements $x \in G$ of order $p$, $q \in \pi(C_G(x))$.
%
%$p=2,q=3$ and for all involutions $t \in G$, $9$ divides $|C_G(x)|$.
%
%$p=3,q=2$ and for all elements $y \in G$ of order $3$, $8$ divides $|C_G(x)|$.\\
%Finally, we let $\rightarrow$ denote the transitive extension of $\vdash$.
%\end{definition}
%
%\begin{lemma}
%Suppose that Hypothesis \ref{4fix} holds, let $\alpha \in \Omega$ and let $p,q \in \pi(G)$.
%
%If $p \in \pi(G_\alpha)$ and $p \rightarrow q$, then $q \in \pi(G_\alpha)$.
%\end{lemma}
%
%\begin{proof}
%This follows from Lemma \ref{Inv-lem:TB}.  If $p,q \ge 5$, then (a) and (c) give the result immediately. If one or both of $p,q$ are $2$ or $3$, then
%our more specialised definition for $\vdash$ and Part (c) of the lemma need to be used.
%\end{proof}

The following three lemmas describe ways to count fixed points and will be heavily used in our analysis of specific groups.

\begin{lemma}\label{numberfpt}
Suppose that $G$ is a group, that $1\neq U\lneq G$ and that $x\in G^\#$.
Then, in the action of $G$ on the coset space $G/U$ by right multiplication, the number of fixed points of
$x$ is exactly
\[\frac{|\{\langle x \rangle^g \leq U\mid g\in G\}|\cdot |N_G(\langle x \rangle)|}{|U|}\,.\]
\end{lemma}

\begin{proof}
Let \(y\in G\). Then \(Uy\) is a fixed point of \(x\) if and only if \(x^{y^{-1}}\in U\). So the number of fixed points of \(x\) is exactly
\[|\{ Uy \mid x^{y^{-1}} \in U\}| =\frac{|\{y\in G\mid x^{y^{-1}}\in U\}|}{|U|}=\frac{\{g\in G\mid x^g\in U\}|}{|U|}  = \frac{|\{\langle x \rangle^g \leq U\mid g\in G\}|\cdot |N_G(\langle x \rangle)|}{|U|}\]
For the last equality, we note that if $\langle x \rangle^G \cap U
= \{\langle x^{g_1}\rangle, \ldots , \langle x^{g_m}\rangle\}$, then $x^g \in U$
if and only if $g \in N_G(\langle x\rangle) g_i$ for some $1 \leq i \leq m$.
\end{proof}

\begin{lemma}\label{Inv-lem:AnzFP}~\\
	Suppose that $G$ is a group, that $1\neq U\lneq G$ and that $x\in G^\#$. \par
	Moreover, let $M_x:=\Menge{\langle y\rangle^U}{y\in x^G\cap U}$ and $k\in\N_0$ be such that $\abs {M_x}=k$, and let
$x_1,\ldots, x_k\in U$ be such that $x=x_1$ is conjugate to $x_2,\ldots, x_k$ in $G$ and that $M_x=\{\langle x_1\rangle^U,\ldots, \langle x_k\rangle^U\}$.
Then, in the action of $G$ on the coset space $G/U$ by right multiplication, the number of fixed points of
$x$ is exactly
\[\sum_{i=1}^k\abs{N_G(\langle x_i\rangle)\colon N_U(\langle x_i\rangle)}.\]
\end{lemma}

\begin{proof}~\\
If $x^G\cap\,U=\varnothing$, then there is nothing to prove.
From now on we suppose that $x^G\cap U\,\neq \varnothing$.

Then without loss $x\in U$.
Further let $y_i \in G$ be such that $x_i^{y_i} = x$.
Let $N:= N_G(\langle x \rangle)$ and $N_i:= N_G(\langle x_i \rangle)$.
Then we see that $N_i^{y_i} = N$.

We claim that $\Lambda := \bigcup\limits_{i = 1}^k UN_i y_i$ is the union of the fixed points of $x$ on
the coset space $G/U$.
Let $y \in G$ such that $Uy$ is fixed by $x$. Then $Uyx = Uy$ implies
$yxy^{-1} \in U$ and therefore  $\langle x \rangle $ is mapped  onto $\langle x_i \rangle$ by $y^{-1}u$ for some $i \in \{1, \ldots , k\}$ and some $u \in U$. This implies  $n:= y_i y^{-1}u \in N_i$. Thus $y =  u n^{-1} y_i \in UN_iy_i$,
which yields that $\Lambda$ contains all the cosets that are fixed by $x$.
On the other hand it is a straightforward calculation that every coset in
$\Lambda$ is fixed by $x$.

Observe that $N$ acts on $UN_iy_i$ and that this action is transitive
as $Uy_i N = U(y_iN y_i^{-1}) y_i = UN_i y_i$. Furthermore, $Uy_i N$ and
$Uy_l N$ are disjoint, if $i \neq l$ by the choice of the $y_i$.

Thus it just remains to show that
$UN_iy_i$ is the union of $|N_G(\langle x_i \rangle):
N_U(\langle x_i \rangle)|$ pair-wise distinct $U$-cosets.
Let $T:=\{g \in N \mid Uy_ig=Uy_i\}$.

The transitive action of $N$ on $UN_iy_i$ by right multiplication
yields that this number of cosets equals $|N:T|$. We note that $T = %(U \cap N_i)^{y_i} =
U^{y_i} \cap N$.
Thus \[|N:T| = |N| / |T| = |N|/|U \cap N_i| = |N_i|/ |N_U(\langle x_i \rangle)|
= |N_i:N_U(\langle x_i \rangle)| ,\] which is the assertion.

\end{proof}

We apply this lemma to the special situation that $U$ is a Frobenius group.

\begin{lemma} \label{FrobCyclic}
	Let \(G\) be a finite group and suppose that \(U \leq G\) is a Frobenius group with Frobenius kernel
	\(K\) and cyclic Frobenius complement \(J\). Let \(x\in U \setminus K\).
	Then \(x\) has exactly \(\frac{|N_{G}(\langle x \rangle)|}{|J|}\) fixed points on the coset space \(G/U\).
\end{lemma}

\begin{proof}
%Der direkte Beweis der Aussage:
%It is \(U y\in G/U\) a fixed point of $x$ if and only if $x \in y^{-1} U y$ if and
%only if $yxy^{-1} \in U$. As $U$ is a Frobenius group we have $(|K|,|J|) = 1$ and get due to Hall's theorem
%that $yxy^{-1}$ is conjugate in $U$ to some element in $J$. Now the fact that $J$ is cyclic
%implies $y^{-1}u \in N:= N_G(\langle x \rangle)$ for some $u \in U$. This shows that $UN$ is
%the union of the cosets that are fixed by $x$. The assertion follows as $N$ acts transitively
%on this set and as the stabilizer of $U$ in $N$ is $U \cap N = N_U(\langle x \rangle = J$.
Let $M_x:=\{\langle y \rangle \mid y\in x^G \cap U\}$, as in Lemma~\ref{Inv-lem:AnzFP}.
We recall that all Frobenius complements are conjugate by 4.3.7 in~\cite{KurzStell} and cyclic by hypothesis.
Therefore, if $y\in x^G \cap U$, then
$\langle y\rangle$ is $U$-conjugate to $\langle x \rangle$ and
it follows that $|M_x| = 1$.
Moreover $N_U(\langle x \rangle) = J$ because $J$ is abelian and because of the Frobenius group property. Now
Lemma~\ref{Inv-lem:AnzFP} implies that $x$ fixes precisely
$|N_G(\langle x \rangle )|/|J|$ cosets in $G/U$.
\end{proof}

For the analysis of Lie type groups later in this article we need generic information about the relevant groups, for example
information about Sylow subgroups in defining characteristic, the structure of involution centralizers and the structure of Levi subgroups of maximal parabolic subgroups.
This information is mostly available in the literature, but not always explicitly, and therefore we collect it here.

\begin{rem}\label{Inv-rem:Ser}~\\
Let $p\in\N$ be an odd prime, $f\in \N$ and $q:=p^f\ge 3$. Let $G$ be one of the groups $\PSL_4(q)$, $\PSL_5(q)$, $\PSU_4(q)$, $\PSU_5(q)$, $\PSp_4(q)$, $G_2(q)$, $^2G_2(q)$ or $^3D_4(q)$ and let $S\in\Syl_p(G)$. \smallskip
	
For the definition of (maximal) parabolic subgroups, the unipotent radical and Levi subgroups we refer to Definition 2.6.4,
Theorem 2.6.5 and Definition 2.6.6 in \cite{GLS3}.

If $t\in G$ is an involution, then Table 4.5.1 in \cite{GLS3} gives information about the structure of $C_G(t)$.
We keep this brief because the only relevant pieces of the table are those where the \textit{coset mod $Inn(K)$} column has the entry ``1'', which means that $t$ is an inner involution. Then we need $L^*=O^{p'}(C_G(t))$ (from columns 6 and 7 of the table) and in the remainder we will use the notation for the conjugacy classes of involutions as well as for $L^*$ from Table 4.5.1 in~\cite{GLS3}. For the groups $G_2(q)$ and $^3D_4(q)$ we refer to \cite{Wilson} for structure information about $C_G(t)$.
The nilpotency class of the Sylow subgroups in defining characteristics is given by
Theorems 3.2.2 and 3.3.1 and Remark 1.8.8 in \cite{GLS3}, except for $G_2(q)$.
We need the height of the highest root of the root system for $G$, again following the notation in \cite{GLS3}.
For the twisted groups we also refer to
13.3 and 3.6 of \cite{Car}).
In the following paragraphs we let $S \in \Syl_p(G)$ and we collect the relevant references, sometimes with additional details.
We let $n\in\{4,5\}$ and we denote the isomorphism type of the Levi subgroups of
$\GL_n(q)$, $\SL_n(q)$ and $\PSL_n(q)$ by $L_{\GL}$, $L_{\SL}$ and $L_{\PSL}$, respectively. We choose similar notation for  $\SU_n(q)$ and $\PSU_n(q)$.
\smallskip

At the end there will be a table collecting the information.

	\smallskip
	\underline{$G= \PSL_4(q)$ and $q\not\equiv 1$ modulo 8}:\par
	We use Sections 3.3.1 (p.44) and 3.3.3 (p.47) in \cite{Wilson} to see that $\abs S=q^6$ and $L_{\PSL} \cong \GL_3(q)/C_{(4,q-1)}$ or $(\GL_2(q)\times\SL_2(q))/C_{(4,q-1)}$.\par
	
	%We use Sections 3.3.1 in \cite[S. 44, Z. 28]{Wilson} and gives that \[\abs G=\frac{1}{(4,q-1)}\cdot q^6\cdot(q^2-1)\cdot (q^3-1)\cdot (q^4-1),\] so $\abs S=q^6$. Then Section 3.3.3 in \cite[S. 47, Z. 3ff]{Wilson} implies that $L_{\PSL} \cong \GL_3(q)/C_{(4,q-1)}$ or $(\GL_2(q)\times\SL_2(q))/C_{(4,q-1)}$.\par
	
	%The Levi subgroups $L_{\GL}$ are isomorphic to $\GL_1(q)\times \GL_3(q)\cong C_{q-1}\times \GL_3(q)$ or to $\GL_2(q)\times \GL_2(q)$. Then $L_{\SL} \cong \GL_3(q)$ or $\GL_2(q)\times\SL_2(q)$.
	
Next let $t\in G$ be an involution.
%As $q-1$ is not divisible by 8, we see that $t \notin t_1^G$, again with notation from
Then Table 4.5.1 in \cite{GLS3}
% In addition, $m=3\equiv -1$ modulo 4, so the same table
gives two possibilities:

$t$ is conjugate to $t_{\frac{m+1}2}=t_2$ in $G$ and, therefore, $L^*\cong\SL_2(q)\circ\SL_2(q)$ or

$q\equiv 3$ modulo 4, $t$ is conjugate to $t_{\frac{m+1}2}'=t_2'$ and $L^*\cong \PSL_2(q^2)$.
(The conjugacy classes denoted by $\gamma_1$, $\gamma_1'$, $\gamma_2$ and $\gamma_2'$ are irrelevant for us because they belong to outer involutions.)

Moreover $S$ has class $3$.
%$\Pi=\{a_1-a_2,a_2-a_3,a_3-a_4\}$ is a system of fundamental roots for $A_3$ (for $\PSL_4(q)$) and \[\alpha=a_1-a_4=\boldsymbol{1}\cdot (a_1-a_2)+\boldsymbol 1\cdot(a_2-a_3)+\boldsymbol 1\cdot(a_3-a_4)\] is the highest root in the root system. Thus $h(\alpha)=3$ is the class of $S$.
\smallskip

	\underline{$G= \PSL_5(q)$ and $q\equiv -1$ modulo 4}:\par
	%Section 3.3.1 in \cite[S. 44, Z. 28]{Wilson} gives that \[\abs G=\frac{1}{(5,q-1)}\cdot q^{10}\cdot(q^2-1)\cdot (q^3-1)\cdot (q^4-1)\cdot(q^5-1)\] and thus $\abs S=q^{10}$. For the Levi subgroups we use Section 3.3.3 in \cite[S. 47, Z. 3ff]{Wilson}:
	Sections 3.3.1 and 3.3.3 in \cite{Wilson} give that $\abs S=q^{10}$ and that
	%$L_{\GL}$ is isomorphic to $\GL_1(q)\times \GL_4(q)\cong C_{q-1}\times \GL_4(q)$ or to $\GL_2(q)\times \GL_3(q)$, $L_{\SL}$ is isomorphic to $\GL_4(q)$ or $\GL_3(q)\times\SL_2(q)$ and then
	$L_{\PSL}\cong\GL_4(q)/C_{(5,q-1)}\quad\text{or}\quad L_{\PSL}\cong(\GL_3(q)\times\SL_2(q))/C_{(5,q-1)}.$

If $t\in G$ is an involution, then
% Since $q-1$ is even and $5$ is odd,
Table 4.5.1 in \cite{GLS3} shows that $t \in t_1^G$ and $L^*\cong\SL_4(q)$ or $t \in t_2^G$ and $L^*\cong \SL_2(q)\times\SL_3(q)$.
%As $q\equiv -1$ modulo 4 holds, these are the only two possibilities. \par
The height of the highest root, and hence the class of $S$, is 4.
%$\Pi=\{a_1-a_2,a_2-a_3,a_3-a_4, a_4-a_5\}$ is a fundamental system of roots for $A_4$ (for $\PSL_5(q)$) and \[\alpha=a_1-a_5=\boldsymbol{1}\cdot (a_1-a_2)+\boldsymbol 1\cdot(a_2-a_3)+\boldsymbol 1\cdot(a_3-a_4)+\boldsymbol 1\cdot(a_4-a_5)\] is the highest root. This implies that $h(\alpha)=4$ is the class of $S$.
\smallskip
	%
	%
	%

	%
	%1 \mod 4$ --> 1$ modulo $4$
	%
	\underline{$G= \PSU_4(q)$ and  $q\not\equiv -1$ modulo 8}:\par
	Again Section 3.6 in \cite{Wilson} applies:
	
	$\abs S=q^{6}$
	%For the Levi subgroups we look at Section 3.6.2 in \cite[S. 67, Z. 36ff]{Wilson}. Then $L_{\SU}$ is isomorphic to \[\GL_2(q^2)\qquad\text{or}\qquad\GL_1(q^2)\times \SU_2(q)\cong C_{q^2-1}\times \SL_2(q)\]
	and $L_{\PSU}\cong\GL_2(q^2)/C_{(4,q+1)}$ or $L_{\PSU}\cong(C_{(q^2-1)}\times \SL_2(q))/C_{(4,q+1)}$.\par

Let $t\in G$ be an involution.
Then there are two possibilities
%By hypothesis $8$ does not divide $q+1$, whence $t \notin t_1^G$, with notation from
 Table 4.5.1 in \cite{GLS3} gives two possibilities for involutions in $G$, namely
%Moreover, $m=3\equiv -1$ modulo $4$, so by the same table
the class of $t_{\frac{m+1}2}=t_2$, where $L^*\cong\SL_2(q)\circ\SL_2(q)$, or $q\equiv 1$ modulo $4$ and there is also the class of $t_{\frac{m+1}2}'=t_2'$, where
$L^*\cong \PSL_2(q^2)$.
%There are no more possibilities for the type of $t$.

%We also have that $\Pi=\{a_1-a_2,a_2-a_3,a_3-a_4\}$ is a fundamental system for the root system for $A_3$ (for $\PSL_4(q)$) and \[\alpha=a_1-a_4=\boldsymbol{1}\cdot (a_1-a_2)+\boldsymbol 1\cdot(a_2-a_3)+\boldsymbol 1\cdot(a_3-a_4)\] is the highest root.
For a root system for $^2A_3$ (for the group $\PSU_4(q)$) we apply a graph automorphism
that interchanges the roots $p_1=a_1-a_2$ and $p_3=a_3-a_4$ and fixes the root $p_2=a_2-a_3$
%This yields $\hat\Pi=\left\{\frac 12(p_1+p_3),p_2\right\}$
(see also Section 13.3 in \cite{Car}), and then we see that
% and for the highest root we have \[\hat \alpha=a_1-a_4=\boldsymbol 2\cdot \bigr(\tfrac 12(p_1+p_3)\bigr)+\boldsymbol 1 \cdot p_2.\]
 the nilpotency class of $S$ is  $3$.

	\smallskip
	\underline{$G= \PSU_5(q)$ and $q\equiv 1$ modulo $4$}:\par
Section 3.6 in \cite{Wilson} gives that
%\[\abs G=\frac{q^{10}}{(5,q+1)}\cdot (q^2-1)\cdot (q^3+1)\cdot (q^4-1)\cdot (q^5+1)\], so
$\abs S=q^{10}$ and
%Section 3.6.2 in \cite[S. 67, Z. 36ff]{Wilson} gives information about the Levi subgroups $L_{\SU} \cong \GL_1(q^2)\times \SU_3(q)\cong C_{q^2-1}\times \SU_3(q)$ or $\GL_2(q^2)$ and
$L_{\PSU}\cong (C_{(q^2-1)}\times \SU_3(q))/C_{(5,q+1)}$ or $\GL_2(q^2)/C_{(5,q+1)}$.\par
	
%Let $t\in G$ be an involution. As $q+1$ is even and 5 is odd,
According to Table 4.5.1 in \cite{GLS3}, there are two classes of involutions: $t_1^G$, where $L^*\cong\SU_4(q)$, and $t_2^G$, where $L^*\cong \SU_2(q)\times\SU_3(q)$.

%With the fundamental root system $\Pi=\{a_1-a_2,a_2-a_3,a_3-a_4,a_4-a_5\}$ of $A_4$ (for the group $\PSL_5(q)$) we calculate the highest root \[\alpha=a_1-a_5=\boldsymbol{1}\cdot (a_1-a_2)+\boldsymbol 1\cdot(a_2-a_3)+\boldsymbol 1\cdot(a_3-a_4)+\boldsymbol 1\cdot(a_4-a_5)\].
For a root system for $^2A_4$ (for $\PSU_5(q)$) we use a graph automorphism that interchanges $p_1=a_1-a_2$ with $p_4=a_4-a_5$ and $p_2=a_2-a_3$ with $p_3=a_3-a_4$, respectively, see also
%With $\hat\Pi=\left\{\frac 12(p_1+p_4),\frac 12(p_2+p_3)\right\}$ (see
Section 13.3 in \cite{Car}.
%we obtain the highest root $\hat \alpha=a_1-a_5=\boldsymbol 2\cdot \bigr(\frac 12(p_1+p_4)\bigr)+\boldsymbol 2\cdot \bigr(\frac 12(p_2+p_3)\bigr)$.
Then it follows that $S$ has class $4$.

	\smallskip
	\noindent\underline{$G= \PSp_4(q)$}:\par
	Section 3.5 in \cite{Wilson} yields that
	%$\abs G=\frac{1}2\cdot q^{4}\cdot(q^2-1)\cdot (q^4-1)$ and hence
	$\abs S=q^{4}$, and with Theorem 3.8 in \cite{Wilson}, the Levi subgroups are isomorphic to $\GL_1(q)\circ \SL_2(q)$ or $\GL_2(q)/C_2$.\par
 Table 4.5.1 in \cite{GLS3} gives that an involution  $t\in G$ is conjugate to $t_{\frac 22}=t_1$, and  $L^*\cong \SL_2(q)\circ\SL_2(q)$ or, otherwise,  $t \in t_2^G$ (if $q\equiv 1$ modulo $4$) or $t \in t_2'^G$ (if $q\equiv 3$ modulo $4$) and in both cases $L^*\cong\PSL_2(q)$.\par
 The class of $S$ is $3$.

%For the fundamental roots we have $\Pi=\{a_1-a_2,2a_2\}$ for the root system $C_2$ (for $\PSp_4(q)$) and the highest root $\alpha=2a_1=\boldsymbol{2}\cdot (a_1-a_2)+\boldsymbol 1\cdot(2a_2)$, giving that $h(\alpha)=2+1=3$. So $S$ has class $3$.
\smallskip
	\underline{$G= \,^3D_4(q)$}:\par
This time Section 4.6.2 and Equation (4.67) in \cite{Wilson} imply that
%\[\abs G=q^{12}(q^8+q^4+1)(q^6-1)(q^2-1)\] and hence
$\abs S=q^{12}$, and by Theorem 4.3 in \cite{Wilson} the Levi subgroups are isomorphic to $\SL_2(q^3).C_{q-1}$ or to $\SL_2(q).C_{q^3-1}$.\par
	%Let $t\in G$ be an involution. Then
	Table 4.5.1 in \cite{GLS3}, together with Section 4.6.5 in \cite{Wilson}, gives that $G$ has only one class of involutions, namely $t_2^G$, with involution centralizer isomorphic to $C_2\dot~(\PSL_2(q^3)\times\PSL_2(q)).C_2$.\par
	%With the fundamental system $\Pi=\{a_1-a_2,a_2-a_3,a_3-a_4,a_3+a_4\}$ for $D_4$ (for $\POm_8^+(q)$) and the highest root \[\alpha=a_1+a_2=\boldsymbol{1}\cdot (a_1-a_2)+\boldsymbol 2\cdot(a_2-a_3)+\boldsymbol{1}\cdot (a_3-a_4)+\boldsymbol{1}\cdot (a_3+a_4)\] we obtain the information for
Starting with a fundamental system for	 $D_4$ we move to $^3D_4$ with a graph automorphism of order~3, acting as a 3-cycle on the set $\{p_1,p_3,p_4\}$ of roots $p_1=a_1-a_2$, $p_3=a_3-a_4$ and $p_4=a_3+a_4$ and fixing the root $p_2=a_2-a_3$ %This gives the fundamental system $\hat\Pi=\left\{\frac 13(p_1+p_3+p_4),p_2\right\}$
(see Section 13.3 in \cite{Car}).
It follows that %and the highest root \[\hat \alpha=a_1+a_2=\boldsymbol 2\cdot p_2+\boldsymbol 3\cdot \bigr(\tfrac 13(p_1+p_3+p_4)\bigr).\] Now $h(\hat\alpha)=2+3=5$ gives the nilpotency class of
$S$ has class $5$.\\
We still need to look at the groups $^2G_2(q)$ and $G_2(q)$, and here we use a different strategy for determining
the nilpotency class of $S$.\smallskip

	\underline{$G= G_2(q)$}:\par
	Section 4.3.3 and Equation (4.25) in \cite{Wilson} yield that
	%\[\abs G=q^{6}(q^6-1)(q^2-1)\] and therefore
	$\abs S=q^{6}$. The Levi subgroups of $G$ are isomorphic to $\GL_2(q)$, by Section 4.3.5 and Table 4.1 in \cite[S.~125, 127]{Wilson}.\par
	By Table 4.5.1 in \cite{GLS3} there is only one class of involutions, with centralizer of structure $C_2\dot~(\PSL_2(q)\times\PSL_2(q)).C_2$, see Section 4.3.6 in \cite{Wilson}.\par
	If $p\ge 5$, then Theorems 3.2.2 and 3.3.1 in \cite{GLS3} give a way to calculate the nilpotency class of $S$ via the height of a highest root, which gives that
	%so we consider a fundamental system $\Pi=\{\alpha_1,\alpha_2\}$ of roots for $G_2$. The highest root is $\alpha=\boldsymbol 2\alpha_1+\boldsymbol 3\alpha_2$ (by Remark 1.8.8 in \cite{GLS3}) and $h(\alpha)=5$, so
	$S$ has class $5$.\par
If $p=3$, then we adopt
the notation on p. 65 in \cite{GLS3}, so we denote by $P_{\{\alpha_2\}}$ a maximal parabolic subgroup and by $U_{\{\alpha_2\}}$ its unipotent radical.  By Definition 2.6.4 and Remark 1.8.8 in \cite{GLS3} we see that
	\begin{align*}
		U_{\{\alpha_2\}}=U_{\{\alpha_2\}}^1&=\langle X_{\alpha_1}, X_{\alpha_1+\alpha_2}, X_{\alpha_1+2\alpha_2}, X_{\alpha_1+3\alpha_2}, X_{2\alpha_1+3\alpha_2}\rangle,\\
		U_{\{\alpha_2\}}^2&=\langle X_{2\alpha_1+3\alpha_2}\rangle=U_\varnothing^5\ug Z(U_\varnothing),\\
		U_{\{\alpha_2\}}^3&=1.
	\end{align*}

Theorem 3.2.2 in \cite{GLS3} gives that $1=U_{\{\alpha_2\}}^3 \neq U_{\{\alpha_2\}}^2\neq U_{\{\alpha_2\}}^1$ is a central series of $U_{\{\alpha_2\}}$ with elementary abelian sections. Moreover $[U_{\{\alpha_2\}},U_{\{\alpha_2\}}]\ug U_{\{\alpha_2\}}^2=\langle X_{2\alpha_1+3\alpha_2}\rangle$ and therefore, for all roots \[\beta\in \{\alpha_1,\alpha_1+\alpha_2,\alpha_1+2\alpha_2,\alpha_1+3\alpha_2,2\alpha_1+3\alpha_2\}=:\Sigma_1,\] it is true that $[X_{\alpha_1},X_\beta]\ug X_{2\alpha_1+3\alpha_2}=U_{\{\alpha_2\}}^2$.\par
We refer to p.99 in \cite{GLS3} and deduce that $P_{\{\alpha_2\}}\cong P_{\{\alpha_1\}}$ and hence $U_{\{\alpha_2\}}\cong U_{\{\alpha_1\}}$. Here it is relevant that all maximal parabolic subgroups are isomorphic even though there are two conjugacy classes of them, see Sections 4.3.5 to 4.3.7 in \cite{Wilson}. There is a graph automorphism $\gamma$  of order 2 that interchanges
 $\alpha_1$ and $\alpha_2$. We also note that $\gamma$ maps positive roots to positive roots and root subgroups to root subgroups and that $U_\varnothing \in \Syl_p(G)$. By Definition 2.6.4 in \cite{GLS3} we have that
	\begin{align*}
		U_\varnothing=U_{\varnothing}^1&=\langle X_{\alpha_1}, X_{\alpha_2}, X_{\alpha_1+\alpha_2}, X_{\alpha_1+2\alpha_2}, X_{\alpha_1+3\alpha_2}, X_{2\alpha_1+3\alpha_2}\rangle,
		U_\varnothing^2=\langle X_{\alpha_1+\alpha_2}, X_{\alpha_1+2\alpha_2}, X_{\alpha_1+3\alpha_2}, X_{2\alpha_1+3\alpha_2}\rangle.
	\end{align*}
Moreover the section $U_\varnothing/U_\varnothing^2$ is elementary abelian by Theorem 3.2.2 in \cite{GLS3} and $[U_\varnothing,U_\varnothing]\ug U_\varnothing ^2$. We also know that $[X_{\alpha_1},X_\beta]\ug X_{2\alpha_1+3\alpha_2}$ for all $\beta\in\Sigma_1$. Before we construct $[U_\varnothing,U_\varnothing, U_\varnothing]\ug [U_\varnothing, U_\varnothing^2]$, we need $[X_{\alpha_2},X_\beta]$ for all $\beta \in\Sigma_1\setminus \{\alpha_1\}$. Using the graph automorphism $\gamma$ of $G_2$ we calculate
	\begin{align*}
		[X_{\alpha_2},U_\varnothing^2]&=([X_{\alpha_2},U_\varnothing^2]^\gamma)^{\gamma}=\left[{X_{\alpha_2}}^\gamma,(U_\varnothing^2)^\gamma\right]^\gamma
		=[X_{\alpha_1},U_\varnothing^2]^\gamma\ug (U_{\{\alpha_2\}}^2)^\gamma\ug Z(U_\varnothing)^\gamma=Z(U_\varnothing).
	\end{align*}

Next we see that
	$[U_\varnothing,U_\varnothing, U_\varnothing]\ug [U_\varnothing, U_\varnothing^2]\ug Z(U_\varnothing)$ and $[U_\varnothing,U_\varnothing, U_\varnothing, U_\varnothing]\ug [U_\varnothing, Z(U_\varnothing)]=1$.

Then the nilpotency class of $U_\varnothing$ is at most 3, because we found a descending central series of length four, and it is
 at least 3, because $U_\varnothing$ is isomorphic to a subgroup of $P_{\{\alpha_2\}}$,  $U_{\{\alpha_2\}}$ has class 2, the $p$-part of $\abs{P_{\{\alpha_2\}}/U_{\{\alpha_2\}}}$ is exactly $q$ and the Sylow $p$-subgroups of the Levi subgroup $L_{\{\alpha_2\}}\ug P_{\{\alpha_2\}}$ act nontrivially on $U_{\{\alpha_2\}}$.

\smallskip
This completes our analysis for now and we show the relevant information in a table.%\newpage

\vspace{0.5cm}
	{\small
	\begin{longtable}{llllllp{3.5cm}}%
	%HEAD
	$G$ & $t^G$ & section $L^*$ in $C_G(t)$ & $\abs S$ & class of $S$ & Levi subgroup $L$ & remarks\\[0.5ex]\hline\hline
	\endhead \\[-2ex]
	%MAIN
	$\PSL_4(q)$	&	$t_2$		& $\SL_2(q)\circ \SL_2(q)$	& $q^6$ 	& 3 & $\GL_3(q)/C_d$&$d=(4,q-1)$\\
				&	$t_2'$	&   $\PSL_2(q^2)$	&	& & $(\GL_2(q)\times \SL_2(q))/C_d$& \\[0.2ex]\hline\\[-2ex]
	$\PSL_5(q)$	&	$t_1$		& $\SL_4(q)$	& $q^{10}$ 	& 4& $\GL_4(q)/C_d$&$d=(5,q-1)$\\
				&	$t_2$		& $\SL_2(q)\times \SL_3(q)$	&	& & $(\GL_3(q)\times \SL_2(q))/C_d$&\\[0.2ex]\hline\\[-2ex]
	$\PSU_4(q)$	&	$t_2$		& $\SL_2(q)\circ \SL_2(q)$	& $q^6$ 	& 3 & $\GL_2(q^2)/C_d$&$d=(4,q+1)$\\
				&	$t_2'$	& $\PSL_2(q^2)$	&	& & $(C_{(q^2-1)}\times \SU_2(q))/C_d$ &\\[0.2ex]\hline\\[-2ex]
	$\PSU_5(q)$	&	$t_1$		& $\SU_4(q)$	& $q^{10}$ 	& 4 & $\GL_2(q^2)/C_d$&$d=(5,q+1)$\\
				&	$t_2$		& $\SU_2(q)\times \SU_3(q)$	&	& & $(C_{(q^2-1)}\times \SU_3(q))/C_d$&\\[0.2ex]\hline\\[-2ex]
	$\PSp_4(q)$	&	$t_1$		& $\SL_2(q)\circ \SL_2(q)$	& $q^4$ 	& 3& $\GL_1(q)\circ \SL_2(q)$&\\
				&	$t_2,t_2'$	& $\PSL_2(q)	$	&	& &  $\GL_2(q)/C_2$&\\[0.2ex]
     \hline\\[-2ex]

	$^3D_4(q)$	&	$t_2$		& $C_2\dot~(\PSL_2(q^3)\circ \PSL_2(q)).C_2$	& $q^{12}$ 	& 5&$\SL_2(q^3). C_{q-1}$&\\
	&&&&& $\SL_2(q)/C_{q^3-1}$&\\[0.2ex]

    \hline\\[-2ex]

	$G_2(q)$	&	$t_1$		& $C_2\dot~(\PSL_2(q)\circ \PSL_2(q)).C_2$	& $q^6$ 	& 3 ($p=3$)& $\GL_2(q)$&\tiny{Two classes of mps (Section 4.3.5. and Table 4.1 in \cite[S. 123ff, 127]{Wilson})}\\
			&				&											&		& 5 ($p>3$)\\[0.5ex]

  \hline\hline\\[-1ex]
	\caption[Sections of involution centralizers, order and nilpotency class of Sylow $p$-subgroups of simple groups of Lie type of sectional 2-rank at most 4 in odd characteristic $p$.]{Conj. class of $t$ in $G$, sections $L^*$ of $C_G(t)$ from Table 4.5.1 in \cite{GLS3}, order of $S\in\Syl_p(G)$ from \cite{Wilson}, nilpotency class of $S$ and Levi subgroup structure for groups $G$ of Lie type of sectional 2-rank at most 4 in odd characteristic $p$. In the remark we abbreviate ``maximal parabolic subgroups'' with \textbf{mps}.}
	\label{Inv-tab:Ser}
	\end{longtable}}
	
\end{rem}

%%%%%%%%%%%%%%%

\section{Small Lie Type Groups}

Here we look at some finite simple groups of Lie Type of small Lie rank. This is because they occur in various places in Theorem \ref{3and2} and therefore they need careful attention. More precisely, we consider the series $\PSL_2(q)$, $\Sz(q)$, $\PSL_3(q)$ and $\PSU_3(q)$, in this ordering, and we classify all possibilities for them to satisfy Hypothesis \ref{4fix}.

\begin{lemma}\label{L2q-main}
    Suppose that Hypothesis \ref{4fix} holds, that $p$ is a prime and that $n\in \N$ and $q:=p^n$ are such that
$G=\PSL_2(q)$ is simple. Let $\alpha\in\Omega$.
Then one of the following holds:
		\begin{enumerate}
			\item \(G=\PSL_2(7)\) and \(G_{\alpha}\cong C_2\) or \(G_{\alpha}\cong \Sym_3\).
			\item \(G=\PSL_2(8)\) and \(G_{\alpha}\) is cyclic of order \(2\) or dihedral of order \(6,14\) or
			\(18\).
			\item \(G=\PSL_2(9)\) and \(G_{\alpha}\) is cyclic of order \(2\), dihedral of order \(6\) or
			\(10\), elementary abelian of order \(9\), or the semi-direct product of an elementary abelian group of order \(9\)
			with a cyclic group of order \(2\).
			\item \(G=\PSL_2(11)\) and \(G_{\alpha}\cong C_3\) or \(G_{\alpha}\cong\Alt_4\).
			\item \(G=\PSL_2(13)\) and \(G_{\alpha}\cong C_3\), \(G_{\alpha}\cong C_{13}:C_3\) or \(G_{\alpha}\cong\Alt_4\).
			\item \(q\geq 17\) is odd. If \(q\equiv 1 \mod 4\) then the point stabilizers are either cyclic of order
			\(\frac{q-1}{4}\) or the semi-direct product of an elementary abelian group of order \(q\) with a cyclic group of
			order
			\(\frac{q-1}{4}\).
			If \(q\equiv -1 \mod 4\) then the point stabilizers are cyclic of order \(\frac{q+1}{4}\).
		\end{enumerate}
\end{lemma}

\begin{proof}

    By Remark~\ref{rem4.1} the statement of the lemma holds whenever \(q\leq 41\). Therefore throughout the rest of the proof suppose
    that \(q\geq 43\).
	
	We recall that \(|G|=q\cdot \frac{q^2-1}{(2,q-1)}\), and detailed information about the subgroup structure of \(G\) is stated in Hauptsatz~II\,8.27
	in~\cite{Huppert-I}, which we will use now without further reference.
	
	\medskip
	
	We will proceed by a case distinction on whether \(q\) and \(|G_{\alpha}|\) are coprime or not.
	
	So first additionally suppose that \(q\) and \(|G_{\alpha}|\) have a common prime divisor \(p\).
	Let \(x\in G_{\alpha}\) be of order \(p\).
	Let \(Q\in \Syl_p(G)\) be such that \(x\in Q\) and \(Q_{\alpha}\in \Syl_p(G_{\alpha})\). Then \(Q\) is elementary
	abelian
	and therefore \(Q\leq N_{G}(\langle x \rangle)\). By our usual lemma this implies \(|Q:Q_{\alpha}|\leq
	|N_{G}(\langle x \rangle):N_{G_{\alpha}}(\langle x \rangle)|\leq 4\). Let \(f\) be a positive integer such that
	\(p^f=q\).
	
	If \(p=2\), then \(|G_{\alpha}|\) is divisible by \(\frac{q}{4}=2^{f-2}\), and \(f\geq 5\). Therefore
	\(f-2>
	\frac{f}{2}\).
	
	 If \(p=3\), then \(|G_{\alpha}|\) is divisible by \(\frac{q}{3}=3^{f-1}\), and \(f\geq 3\).
	 Hence
	\(f-1>\frac{f}{2}\).
	
	Finally, if \(p\geq 5\), then \(|G_{\alpha}|\) is divisible by \(q\) by our usual Sylow lemma.
	In particular \(|Q_{\alpha}|\geq \frac{q}{4}\geq \frac{43}{4}>10\).
	
	We now show that \(G_{\alpha}\) is not a \(p\)-group, and the strategy depends on whether \(q\) is even or odd.
	First suppose that \(q\) is odd.
	Assume for a contradiction that \(G_{\alpha}\) is a \(p\)-group.
	Then there exists an element \(a\in
	G_{\alpha}\) of order \(p\)
	with exactly \(4\) fixed points.
	
	Lemma~\ref{numberfpt} says that the number of fixed points of \(a\) is
	\(\frac{|\{\langle a \rangle^g\leq G_{\alpha}\mid g\in G\}|\cdot |N_G(\langle a
		\rangle)|}{|G_{\alpha}|}\), so this number is $4$. As \(|N_G(\langle a
	\rangle)|\) is divisible by \(p^f\) and \(4\) is coprime to \(p\), we deduce that \(|G_{\alpha}|\) is divisible by \(p^f\) and
	hence \(Q_{\alpha
	}=Q\).
	Then our usual lemma gives the contradiction
	\(4\geq |N_G(Q):N_{G_{\alpha}}(Q)|=q\cdot \frac{q-1}{2}/q=\frac{q-1}{2}\geq \frac{42}{2}=21\). Hence
	\(G_{\alpha}\) is not a \(p\)-group in this case.
	
	So now suppose \(q\) is even. As all involutions in \(\SLi_2(q)\cong \PSL_2(q)\) are conjugate,
	the number of involutions in \(Q_{\alpha}\) (i.e. \(|Q_{\alpha}|-1\)) is exactly \(|\{\langle x \rangle^g \leq Q_{\alpha}\mid g\in G\}|\).
	By Lemma~\ref{numberfpt} the number of
	fixed points of \(x\) is
	
	\(\frac{|\{\langle x \rangle^g\leq G_{\alpha}\mid g\in G\}|\cdot |N_G(\langle x
		\rangle)|}{|G_{\alpha}|}\geq \frac{|\{\langle x \rangle^g\leq Q_{\alpha}\mid g\in G\}|\cdot |N_G(\langle x
		\rangle)|}{|G_{\alpha}|}=\frac{(|Q_{\alpha}|-1)\cdot |N_G(\langle x \rangle)|}{|G_{\alpha}|} \geq
	\frac{(|Q_{\alpha}|-1)\cdot q}{|G_{\alpha}|}\).
	
	As \(x\) fixes at most four points, we conclude that \(|G_{\alpha}|\geq \frac{q\cdot (|Q_{\alpha}|-1)}{4}\).
	We %can see, using \texttt{GAP} (see Appendix)
	see by Remark~\ref{rem4.1}, that $\PSL_2(64)$ does not act with fixity 4, so
	from now on we can suppose that  and \(q\geq 128\) in the present case.
	Then \(\frac{\frac{q}{4}-1}{4}\geq \frac{31}{4}>4\), and again $G_ \alpha$ is not a \(p\)-group.

	Together with the information about \(|Q_{\alpha}|\), an inspection of the list of subgroups of \(\PSL_2(q)\)
	reveals that
	\(G_{\alpha}\) is a subgroup of a Frobenius group of order \(q\cdot \frac{q-1}{(2,q-1)}\) that is
	not a \(p\)-group. The Frobenius kernel is \(Q_{\alpha}\) the unique Sylow \(p\)-subgroup of \(G_{\alpha}\).
	Let \(y\in G_{\alpha}\) have prime order dividing \(\frac{q-1}{(2,q-1)}\).
	Then by Lemma~\ref{FrobCyclic}
	the number of fixed points of \(y\) is
	\(\frac{|Q_{\alpha}|\cdot |N_G(\langle y \rangle)|}{|G_{\alpha}|}=\frac{|Q_{\alpha}|\cdot 2\cdot
		\frac{q-1}{(2,q-1)}}{|G_{\alpha}|}\).
	By Lemma~(numberfpt) the number of
	fixed points of \(x\) is \(\frac{|\{\langle x \rangle^g\leq G_{\alpha}\mid g\in G\}|\cdot |N_G(\langle x
		\rangle)|}{|G_{\alpha}|}\). This number coincides with
	\(\frac{|\{\langle x \rangle^g\leq Q_{\alpha}\mid g\in G\}|\cdot |N_G(\langle x
		\rangle)|}{|G_{\alpha}|}\) because \(Q_{\alpha}\) is a characteristic subgroup of \(G_{\alpha}\).
		
		Moreover,  \(|\{\langle x \rangle^g\leq Q_{\alpha}\mid g\in G\}|\) is bounded from above
	by the number of distinct subgroups of order \(p\) of \(Q_{\alpha}\), and \(|N_G(\langle x\rangle)|\) is
	bounded from above by \(|C_G(x)|\cdot |\aut(\langle x \rangle)|=q\cdot (p-1)\). Consequently \(x\) fixes at most
	\(\frac{\frac{|Q_{\alpha}|-1}{p-1}\cdot q\cdot (p-1)}{|G_{\alpha}|}=\frac{q(|Q_{\alpha}|-1)}{|G_{\alpha}|}\) points in
	\(\Omega\).
	Hence $|\FO(x)| \le |\FO(y)|$.
	
	Since the maximal number of fixed points of nontrivial elements in \(G_{\alpha}\) is obtained by an element of prime
	order, we see that
	\(y\) fixes four points. Then the calculation above shows that \(|G_{\alpha}|=|Q_{\alpha}|\cdot \frac{q-1}{2\cdot
		(2,q-1)}\).
	As \(G_{\alpha}\) is a Frobenius group, this is only possible if \(|Q_{\alpha}|=q\) and \(q\) is odd. %TODO By KS 4.1.5
	%und Notizen 17.12.2021
	So \(|G_{\alpha}|=q\cdot
	\frac{q-1}{4}\). Hence \(q-1\) is divisible by \(4\) and so \(q\equiv 1 \mod 4\). This case is listed in (f).
	This finishes the analyzes in the case that \(q\) and \(|G_{\alpha}|\) have a common prime factor.\\
	
	Now we suppose that \(q\) and \(|G_{\alpha}|\) are coprime. Then \(|G_{\alpha}|\) divides \(q^2-1\). Let
	\(x\in G_{\alpha}\) be of prime order \(p\) and such that $|\FO(x)|=4$.
	Then \(p\) divides \(q-1\) or \(q+1\) and we let
	\(\varepsilon \in \{-1,1\}\) be such that \(p\) divides \(q-\varepsilon\). The normalizer of \(\langle x
	\rangle\) is a dihedral group of order \(2\cdot \frac{q-\varepsilon}{(2,q-1)}\).
	Let \(C\) be the cyclic subgroup of order \(\frac{q-\varepsilon}{(2,q-1)}\) of that normalizer.
	With our usual lemma we see that \(|C:C_{\alpha}|\leq |N_G(\langle x \rangle):N_{G_{\alpha}}(\langle x
	\rangle)|\leq 4\).
	Hence \(G_{\alpha}\) contains a cyclic subgroup of order at least \(\frac{q-\varepsilon}{4\cdot (2,q-1)}\geq
	\frac{42}{8}>5\).
	As \(|G_{\alpha}|\) and \(q\) are coprime, an inspection of the list of subgroups of \(G\) shows two possibilities for
	\(G_{\alpha}\):
	a cyclic group of order dividing
	\(\frac{q-\varepsilon}{(2,q-1)}\) or a dihedral group.
	
	For a contradiction we assume that \(G_{\alpha}\) is a dihedral group.
	Then there is an involution \(t\) outside the cyclic group \(C\) of \(G_{\alpha}\). However, \(N_G(t)\) contains a
	cyclic group
	\(D\) that has trivial intersection with \(C\). So by our usual lemma \(|D:D_{\alpha}|\leq
	|N_G(t):N_{G_{\alpha}}(t)|\leq 4\). On the other hand
	\(|D:D_{\alpha}|=\frac{|D|}{2}=\frac{q-\varepsilon}{2(2,q-1)}\geq \frac{42}{4}>4\). This contradiction shows that
	\(G_{\alpha}\) is not a dihedral group.
	
	Therefore \(G_{\alpha}\) is cyclic. Then by Corollary~(cyclicnumberfpt) the number of fixed points
	of a nontrivial element \(c\) of \(G_{\alpha}\) is \(\frac{|N_G(\langle c\rangle)|}{|G_{\alpha}|}=\frac{2\cdot
		\frac{q-\varepsilon}{(2,q-1)}}{|G_{\alpha}|}\). As \(G_{\alpha}\) contains nontrivial elements that fix exactly
	four
	points, it follows that every nontrivial element fixes exactly four points and
	hence \(|G_{\alpha}|=\frac{q-\varepsilon}{2\cdot (2,q-1)}\). As
	\(q-\varepsilon\) is only divisible by \(2\) if \(q\) is odd, we see that \((2,q-1)=2\) and
	\(|G_{\alpha}|=\frac{q-\varepsilon}{4}\). In conclusion, \(q-\varepsilon\) is divisible by \(4\) and so \(q\equiv \varepsilon
	\mod 4\). This case is listed in (f).
\end{proof}

The Suzuki groups appear as potential candidates for fixity 4 actions in Theorem \ref{3and2}, and the next lemma gives the exact possibilities.

\begin{lemma}\label{Suzuki}
Let \(n \in \N, q:=2^{2n+1}\) and
		\(G=\Sz(q)\).
		Suppose that $G$ satisfies Hypothesis \ref{4fix} for some set $\Omega$ and let $\alpha \in \Omega$.
		 Then \(G_{\alpha}\) is cyclic of order \(q+\sqrt{2q}+1\) or \(q-\sqrt{2q}+1\).
\end{lemma}

\begin{proof}
	First we note that \(|G|=q^2(q^2+1)(q-1)=q^2 \cdot (q+\sqrt{2q}+1) \cdot (q-\sqrt{2q}+1) \cdot (q-1)\) and that the numbers \(q^2,
	q+\sqrt{2q}+1, q-\sqrt{2q}+1, q-1\) are pairwise coprime. (For the group order see for example p. 117 in \cite{Wilson}.)
	We also refer to Table \ref{SmallGroupsTable} for $\Sz(8)$, so from now on we suppose that
	$q \ge 32$ and we will several times apply Theorem
	4.1 on p.\,117
	in~\cite{Wilson} for the subgroup structure of $G$.
	With Hypothesis \ref{4fix} we let $x \in G_\alpha$ be of prime order $p$ such that $x$ fixes four points.

	First suppose that $p=2$ and let $x \in S \in \syl_2(G)$.
	Then we see, with Lemma \ref{Inv-lem:TB} (c) and the structure of $S$, that $G_\alpha$
	contains a subgroup of index $2$ or $4$ of $S$ and hence it contains $Z(S)$. Let $N:=N_G(Z(S))$. Then $|N|=q^2(q-1)$ and,
	as $q-1$ is odd and all elements in $Z(S)$ have two or four fixed points,
	it follows that $q-1$ divides $|G_\alpha|$. Then $G_\alpha=N$ contains a Sylow $2$-subgroupof $G$,
	which contradicts Lemma \ref{Inv-lem:TB2}.
	
	Now we know that $p$ is odd, so
	\(p\) divides one of the numbers  \((q+\sqrt{2q}+1), (q-\sqrt{2q}+1)\) or \(q-1\) and
	\(X:=\langle x \rangle\) is contained in a
	maximal subgroup of $G$ of structure \(D_{2(q-1)}\), \(C_{q+\sqrt{2q}+1}:C_4\) or \(C_{q-\sqrt{2q}+1}:C_4\).
	The number of fixed points of \(x\) can be determined with Lemma \ref{FrobCyclic}, so \(4=|\FO(x)|=
	\frac{|N_{G}(X)|}{|G_{\alpha}|}\).
	If \(p\) divides \(q-1\), then \(4=\frac{2(q-1)}{|G_{\alpha}|}\) gives a contradiction because \(q-1\) is odd.
	So now there exists \(\varepsilon\in\{-1,1\}\) such that \(p\) divides \(q+\varepsilon\sqrt{2q}+1\).
	
	Then \(4=\frac{(q+\varepsilon\sqrt{2q}+1)\cdot 4}{|G_{\alpha}|}\) implies that
	\(|G_{\alpha}|=q+\varepsilon\sqrt{2q}+1\). All subgroups of order \(q+\varepsilon\sqrt{2q}+1\) of a group of structure
	\(C_{q+\varepsilon\sqrt{2q}+1}:C_4\) are cyclic, so this yields the statement of the lemma.
\end{proof}

Next we consider the series $\PSL_3(q)$ and $\PSU_3(q)$.

\begin{lemma}\label{PSLPSU3}

Let \(q\geq 3\) be a prime power, let \(\varepsilon\in\{-1,1\}\) and let
		\(G=\PSL^{\varepsilon}_3(q)\). Then \(G\) acts transitively and with
		fixity \(4\) %on a set of size at least \(5\)
		if and only if \(G=\PSU_3(3)\) and the point stabilizers have structure \(((C_3\times C_3):C_3):C_8\).

\end{lemma}

\begin{proof}
    By Remark~\ref{rem4.1} the statement of the lemma holds whenever \(q\leq 9\). Therefore,
    from now on, we let \(q\geq 11\). %TODO überprüfen, dass es auch tatsächlich im Remark passiert

    First we suppose that \(G\) acts transitively and with fixity \(4\) on a set~\(\Omega\) and we let \(\alpha
	\in
	\Omega\).
	
	Assume for a contradiction that \(G_{\alpha}\) has even order. Let \(t\in G_{\alpha}\) be an involution.
	
	The analysis now depends on whether \(q\) is even or odd.
	First suppose that \(q\) is even, which means that $q \ge 16$. Then by Lemma~\ref{3and2simple}, Satz~1
	in~\cite{Ben},
	Theorem~2 in~\cite{GW},
	the Third Main Theorem in~\cite{ABG2} and
	the Main Theorem
	in~\cite{GoHa1974}, we know that
	\(G=\PSU_3(q)\) and that \(G\) has a strongly embedded subgroup.
	Remembering that $q$ is even, we let \(f\in \N\) be such that \(q:=2^f\). Let \(T\in \Syl_2(G)\) be such that \(T_{\alpha}\in
	\Syl_2(G_{\alpha})\) and that \(t\in T_{\alpha}\). Then \(t\in
	Z(T)\) by statement (4) on page~535
	in~\cite{Ben}.
	As \(t\) has at most four fixed points, we see that \(4\geq |C_T(t):C_{T_{\alpha}}(t)|=|T:T_{\alpha}|\). Hence
	\(|T_{\alpha}|\) is one of the following numbers: \(|T|=q^3\), \(\frac{|T|}{2}=\frac{q^3}{2}=2^{3f-1}\) or
	\(\frac{|T|}{4}=\frac{q^3}{4}=2^{3f-2}\). So the order of \(G_{\alpha}\) is
	divisible by \(\frac{q^3}{4}\geq 2^{10}\).
	Inspection of Tables~8.5 and 8.6
	in~\cite{BHRD} reveals that \(G_{\alpha}\)
	lies in a maximal subgroup such that its full pre-image in \(\SUn_3(q)\) has structure \(E_q^{1+2}:(q^2-1)\).
	%	 of a group of order \(q^3\frac{(q+1)}{\ggT(3,q+1)}(q-1)\) because if \(2^{3f-2}\) divided \(2^{3m}\)
	%	where \(\frac{n}{m}\) is an odd prime \(r\) then \(3m \geq 3n-2=3rm-2\) and so \(2\geq 3m(r-1)\geq 6\).
	Let \(I\) be the set of all involutions of \(G\) and  assume that \(I \subseteq N_{G}(T)\). As \(G\) is simple and \(\langle I \rangle\unlhd G\), this implies that \(\langle I \rangle = G\) and in particular \(N_{G}(T)=G\), which is impossible.
	We deduce that \(I \not \subseteq N_{G}(T)\) and we let \(s \in G\) be an involution such that
    \(s \notin N_{G}(T)\). By statement (6) on page~534
	in~\cite{Ben} we know that \(|C_G(s)|\) is
	divisible by \(\frac{q+1}{(3,q+1)}\).
	
	We will need the following observation:
	
	($\ast$) ~~If \(T \le U\leq G\), then \(T \unlhd U\).
	
	This is because, by Tables~8.5 and 8.6
	in~\cite{BHRD}, $U$ lies in a maximal
	subgroup with full pre-image in \(\SU_3(q)\) of structure \(E_q^{1+2}:(q^2-1)\).
	%	by the argumentation above (and
	%	\cite{Hartley1925DeterminationOfTheTernaryCollineationGroupsWhoseCoefficientsLieInTheGF2n}) it is isomorphic to
	%	case~(ii) in~Theorem~3.9 in~\cite{Wilson2009TheFiniteSimpleGroups_Ha4_MI2010-717} and therefore the Sylow
	%	\(2\)-subgroup
	%	is normal in \(U\).
		
%	Assume that \(s\) and \(t\) are not conjugate in \(G\). Then Lemma~2.4
%	in~\cite{Ben} shows that there
%	exists an involution \(u\) in \(C_G(s) \cap C_G(t)\). Then \(T\leq C_G(t)\)
%	implies that \(T \unlhd C_G(t)\) by ($\ast$), and therefore \(C_G(t)\leq N_G(T)\). Now \(u\in C_G(t)\le N_{G}(T)\) and statement (4) on page~534
%	in~\cite{Ben} gives that \(u\in
%	Z(T)\). Once more ($\ast$) applies, so \(T\trianglelefteq C_G(u)\). As \(u\in C_G(s)\), we deduce that \(\langle s, T\rangle
%	\leq C_G(u)\) is a \(2\)-group of order at least \(|T|\), which implies that \(\langle s, T\rangle= T\) and therefore \(s\in T\).
%	This contradicts the choice of \(s\), which proves that
    By Lemma 4.1\,(i) in~\cite{Ben}, $s$ and $t$ are conjugate.
	Let \(h\in G\) be such that \(t^h=s\) and let \(\beta := \alpha^h\).
	Then \(s\in G_{\beta}\), \(s\) has at most four fixed points and therefore \(|C_G(s) :
	C_{G_{\beta}}(s)| \leq 4 \).
	We notice that \(|C_G(s)|\) is divisible by \(\frac{q+1}{(3,q+1)}\), an odd number. This number is divisible by 3 only if \(q+1\) is divisible by \(9\), which means that
	\(|C_{G_{\beta}}(s)|\) is divisible by \(\frac{q+1}{(9,q+1)}>4\).
	
	Assume that \(\frac{q+1}{(9,q+1)}\) is a \(3\)-power. Then there exists some \(l\in \N\) such that
	\(q+1=2^f+1=3^l\).
	If \(l\) is odd, then \(2^f=3^l-1=(3-1)(3^{l-1} + 3^{l-2} + \ldots +3+1)\). The second factor has \(l\) odd summands, so it is odd and it divides \(2^f\). It follows that \(l=1\) and then
	\(f=1\). In conclusion \(q=2\), which contradicts the assumption.
	
	So \(l\) is even. Let \(l_2:=\frac{l}{2}\). Then \(2^f=3^l-1=(3^{l_2}-1)(3^{l_2}+1)\). If \(l_2\)
	is odd, then we argue as in the previous paragraph and it follows that \(l_2=1\), then \(l=2\) and hence \(f=3\). This means that \(q=8<11\),
	whereas we only consider the case $q \ge 11$.
	We deduce that \(l_2\) is even. Then \(3^{l_2} + 1 \equiv (-1)^{l_2}
	+1
	\equiv 1 +1 \equiv 2 \mod 4\). It follows that \(3^{l_2} + 1=2\), but then \(l_2=0\), and this contradicts our assumption
	that \(l \in \N\).
	Therefore \(\frac{q+1}{(9,q+1)}\) is not a \(3\)-power, but it is odd and hence divisible by a prime greater than \(3\).
	
	Let \(r \ge 5\) be a prime dividing \(q+1\) and \(|G_{\alpha}|\). Then \(G_{\alpha}\) contains a
	Sylow \(r\)-subgroup $R$ of \(G\), and we let \(k \in \N\) be such that \(|R|=r^k\). We recall that \(|G_{\alpha}|\) divides \(q^3(q+1)\frac{(q-1)}{(3,q+1)}\). Moreover the numbers \(q\), \(q-1\) and \(3\) are
	coprime to \(r\). Then it follows that \(r^k\) divides \(q+1\).
	As \(|G|=\frac{q^3}{(3,q+1)}(q^3+1)(q^2-1)=\frac{q^3}{(3,q+1)}(q^2-q+1)(q-1)(q+1)^2\), we deduce that
	\(r^{2k}\) divides \(|G|\), contrary to the fact that $R \in \syl_r(G)$ and $|R|=r^k$.
	
	Therefore \(q\) is
	odd.
	
	Now we know that all involutions in \(G\) are conjugate, by Proposition~2.1 in~\cite{GLS6},
	and then Table 8.3 and Table 8.5
	in~\cite{BHRD} show that $t$ is contained in the centre of a
	group %isomorphic to \(\GLi_{\varepsilon}(q)/\ggT(3,q-\varepsilon)\)
	\(U\) such that the full pre-image of $U$ in \(\SLi^{\varepsilon}_3(q)\) is isomorphic to \(\GLi^{\varepsilon}_2(q)\). In particular
	\(|C_G(t)|\) is divisible by \(\frac{q\cdot (q-\varepsilon)^2\cdot (q+\varepsilon)}{(3,q-\varepsilon)}\). As
	\(|C_G(t):C_{G_{\alpha}}(t)|\leq 4\) and \(\frac{q}{4}>2\), we find an element \(x\in
	G_{\alpha}\) of prime order \(p\) such that \(q\) is a power of $p$. Let \(Q\in \Syl_p(G)\) be such that \(x\in Q\). Then \(Z(Q)\leq
	C_G(x)\)
	and therefore \(|Z(Q):Z(Q)\cap G_{\alpha}|\leq |C_{G}(x):C_{G_{\alpha}}(x)|\leq 4\). Next we use the fact that \(|Z(Q)|=q\), whence
	%TODO Referenz
	\(Z(Q)\cap G_{\alpha}\) contains a nontrivial element~\(z\). Therefore \(Q\leq C_{G}(z)\) and we see that
	\(|Q:Q_{\alpha}|\leq |C_{G}(z):C_{G_{\alpha}}(z)|\leq 4\). This implies that \(G_{\alpha}\) is divisible by
	\(\frac{q^3}{(3,q)}\) or \(\frac{q^3}{(4,q)}\). We recall that \(|G_{\alpha}|\) is also divisible by \(\frac{q\cdot
		(q-\varepsilon)^2\cdot (q+\varepsilon)}{(3,q-\varepsilon)}\) and that \(q\geq 11\). Then inspection of the Tables 8.3,
	8.4,
	8.5 and 8.6 in~\cite{BHRD} %TODO include
	%Errata
	shows that \(G=\PSL_3(q)\) and that \(G_{\alpha}\) lies in a maximal subgroup \(M\) such that its full pre-image in
	\(\SLi_3(q)\) has structure \(E_q^2:\GLi_2(q)\).
	Let \(E\) be the normal elementary abelian subgroup of order~\(q^2\).
	Then \(Q_{\alpha}\) contains an element \(a\) in a complement of \(E\) in \(M\). Thus, \(|C_M(a)|=q^2\cdot \frac{q-1}{(3,q-1)}\). On the other hand, \(|C_G(a)=q^3\cdot \frac{q-1}{(3,q-1)}\)
	by Table~2 in~\cite{Simpson}.
%	Then \(|E_{\alpha}|>q\) and therefore there exists an element \(a\in E\setminus Z(Q)\). Let \(A:=\langle a
%	\rangle\).
%	An element \(g\in M\) that normalizes \(A\) maps \(a\) to an element in the same %one-diminsional vector space
%	subgroup of \(E\) of order \(q\). In \(\GLi_2(q)\), the only elements that preserve these
%	subgroups have order dividing \(q-1\), and this implies that \(|N_M(A)|_p=q^2\).
%	As \(\GLi_2(q)\) acts transitively on \(E\), there exists an element in \(M\) that maps \(a\) to an element \(y\in
%	E\cap Z(Q)\). Therefore \(|N_G(\langle a \rangle )|=|N_G(\langle y \rangle)|\) is divisible by \(q^3\). Moreover,
%	\(N_{G_{\alpha}}(A)\leq N_M(A)\), which implies that \(q\) divides \(|N_G(A):N_{G_{\alpha}}(A)|\). But this contradicts the fact that
%	\(|N_G(A):N_{G_{\alpha}}(A)|\leq 4\).
%
%!!!Hier muss die Argumentation anders erfolgen. Man kann da Dinge einfach ausrechnen. %TODO
As \(C_{G_{\alpha}}(a)\leq C_M(a)\), therefore, \(q\) divides \(|C_G(a):C_{G_{\alpha}}(a)|\) but this contradicts the fact that
	\(|C_G(a):N_{G_{\alpha}}(a)|\leq 4\).
	Now all possibilities where \(G_{\alpha}\) has even order are excluded.
	
	We deduce that
	\(|G_{\alpha}|\) is odd.
	
	Moreover, a prime dividing \(|G_{\alpha}|\) %is odd and
	divides at least one of the numbers \(q\), \(q-\varepsilon\), \(q+\varepsilon\)
	and \(q^2+\varepsilon q+1\).
	
	First suppose that there exists an element \(x\) in \(G_{\alpha}\) of prime order \(r\) dividing \(q-\varepsilon\). If
	\(r\geq 5\), then \(G_{\alpha}\) contains a Sylow \(r\)-subgroup of $G$ by Lemma~\ref{Inv-lem:TB}\,(a).
%	, because \(q\), \(q+\varepsilon\) and \(q^2+\varepsilon
%	q +1\) are not divisible by \(r\).
	Then $G_\alpha$ contains an element \(z\) of order $r$ in the centre of a maximal subgroup \(U\) such that the full pre-image of
	$U$ in \(\SLi^{\varepsilon}_3(q)\) is isomorphic to \(\GLi^{\varepsilon}_2(q)\).
	We keep this in mind and look at the case \(r=3\).
%	Let \(P\) be a Sylow \(3\)-subgroup of \(G\) containing \(x\). Then (10-1) in~\cite{GL} gives that \(r(P)=2\). The number \(q+\varepsilon\) is coprime to \(3\), and \(q+\varepsilon q+1\) is divisible by \(3\), but not by \(9\). As \(|G|=q^3\cdot (q^2+q\varepsilon) (q+\varepsilon) \cdot \frac{(q-\varepsilon)^2}{(3,q-\varepsilon)}\), this implies that \(|G|_3=((q-\varepsilon)^2)_3\). Now we apply Theorem~7.6.1 in~\cite{GLS3}: \(G\) does not have a strongly \(3\)-embedded subgroup, which means that the only remaining case of Lemma~\ref{3and2simple} is that \(P\) is elementary abelian of order \(9\). %TODO Diesen Zwischenschritt loswerden
%	Let $M$ be a maximal
%	subgroup whose full pre-image in \(\SLi^{\varepsilon}_3(q)\) has structure \((q-\varepsilon)^2:\Sym_3\). Then $M$ contains a Sylow \(3\)-subgroup of \(G\), so without loss
%	\(P \le M\).
	Let \(P\) be a Sylow \(3\)-subgroup of \(G\) containing \(x\). Then \(|N_M(\langle x\rangle)|\) is divisible by \(9\) and by Table~2 in~\cite{Simpson} also by \(2\).
	This means that %$x\in M$ and that
	\(|N_M(\langle x\rangle)|\) is divisible by \(18\), moreover
	%TODO Der Beweis dafür liegt in der Struktur in SL_{\varepsilon}(3,q) begründet
%	%Then \(N_{N_{\varepsilon}}(Y)\) is divisible by
%	%\(2\cdot 3 \cdot 3\) by Lemma~\ref{SLe3q}\,4.
%	%TODO reicht das so?
%	By Theorem~3.5 and Theorem~3.9 in~\cite{Wilson} we know that \(\SL_3^{\varepsilon}(q)\) has a subgroup of index at most \(q-\varepsilon\) of a subgroup isomorphic to \(\GL_1^\varepsilon(q)\wr \Sym_3\).
%	Then calculations in \(\SL_3^{\varepsilon}(q)\) show that every element of order \(3\) is either inverted or centralized
%	by some involution in \(\SL_3^{\varepsilon}(q)\). Therefore \(N_{N_{\varepsilon}}(Y)\) is divisible by
%	\(2\cdot 3 \cdot 3\).
	\(|N_G(\langle x\rangle):N_{G_{\alpha}}(\langle x\rangle)|\leq 4\). Since this index is not divisible by $6$, it follows that \(P \le G_{\alpha}\). Again \(G_{\alpha}\) contains an
	\(r\)-element \(u\) in the centre of a maximal subgroup \(U\), where the full pre-image of $U$ in \(\SLi^{\varepsilon}_3(q)\) is isomorphic to
	\(\GLi^{\varepsilon}_2(q)\).
	
	Then \(U\leq C_G(u)\) and hence \(|C_G(u)|\) is divisible by \(\frac{q\cdot (q-\varepsilon)^2\cdot
		(q+\varepsilon)}{(3,q-\varepsilon)}\). As \(|C_G(u):C_{G_{\alpha}}(u)|\leq 4\) and
		\(\frac{q\cdot (q-\varepsilon)\cdot
		(q+\varepsilon)}{(3,q-\varepsilon)} \) is divisible by \(8\), this  contradicts the fact that \(|G_{\alpha}|\) is odd.
	Therefore \(|G_{\alpha}|\) and \(q-\varepsilon\) are coprime.
	
	Next suppose that there exists an element \(x\) in \(G_{\alpha}\) of prime order \(r\) such that $r$ divides \(q\). Let \(x\in R\in
	\Syl_r(G)\). Then \(Z(R)\leq C_G(x)\) and consequently \(|Z(R):Z(R)\cap G_{\alpha}|\leq
	|C_{G}(x):C_{G_{\alpha}}(x)|\leq 4\). Therefore \(|Z(R)|=q\), %TODO Referenz
	and it follows that \(Z(R)\cap G_{\alpha}\) contains a nontrivial element \(z\). Now \(|C_{G}(z)|=q^3\cdot
	\frac{q-\varepsilon}{(3,q-\varepsilon)}\). %TODO Das kann an man einfach ausrechnen, weil alle Sylow UG konjugiert
	%sind
	As \(\frac{q-\varepsilon}{(3,q-\varepsilon)}\geq \frac{10}{3}\), point stabilizer have order coprime to \(q-\varepsilon\)
	and \(|C_G(z):C_{G_{\alpha}}(z)|\leq 4\), it follows that \(q-\varepsilon=12\).
	If \(G=\PSL_3(13)\), then \(|N_G(\langle z \rangle)|\) is divisible by \(16\), which gives a contradiction to the fact that \(|G_{\alpha}\) is odd.
	Therefore \(G=\PSU_3(11)\). But then \(|N_G(\langle z \rangle)|\) is divisible by \(8\), which is also impossible.
	%TODO Rechnungen dafür in GAP?
%	But this contradicts the previous paragraph.
	
	Hence \(|G_{\alpha}|\) is coprime to \(q\cdot
	(q-\varepsilon)\).
	
	Now suppose that there exists an element \(x\) in \(G_{\alpha}\) of prime order \(r\) dividing \(q+\varepsilon\). If
	\(r\geq 5\), then \(G_{\alpha}\) contains a Sylow \(r\)-subgroup, because \(q\), \(q+\varepsilon\) and \(q^2+\varepsilon
	q +1\) are not divisible by \(r\).
	In particular there exists a maximal subgroup $U$ of order divisible by $r$ such that
	the full pre-image of $U$ in
	\(\SLi^{\varepsilon}_3(q)\) is isomorphic to  \(\GLi^{\varepsilon}_2(q)\). If \(r=3\), then
	(10-1) in~\cite{GL} shows that \(G\) has cyclic Sylow \(3\)-subgroups, and in particular all subgroups of order \(3\) are conjugate in
	\(G\).
	In both cases there exists an element $b$ of prime order $r$ in a maximal subgroup \(U\) such that the full pre-image of $U$ in \(\SLi^{\varepsilon}_3(q)\)
	is isomorphic to \(\GLi^{\varepsilon}_2(q)\).
	Then \(|C_{G}(b)|\) is divisible by \((q-\varepsilon)\), \(|G_{\alpha}|\) is coprime to \(q-\varepsilon\) and  \(|C_G(b):C_{G_{\alpha}}(b)|\leq 4\). Together all this implies that
	 \(q-\varepsilon=12\).
	If \(G=\PSL_3(13)\), then \(r=7\) and \(|N_G(\langle x \rangle)|\) is divisible by \(16\), which gives a contradiction to the fact that \(|G_{\alpha}\) is odd.
	So \(G=\PSU_3(11)\), but then \(r=5\) and \(|N_G(\langle x \rangle)|\) is divisible by \(16\). This is impossible.
	We deduce that
	\(|G_{\alpha}|\) and \(q\cdot (q^2-\varepsilon)\cdot (q+\varepsilon)\) are coprime. In particular, since one of the numbers \(q\), \(q-\varepsilon\) and
	\(q+\varepsilon\) is divisible by~\(3\), it follows now that \(|G_{\alpha}|\) is coprime to \(6\).
	
	Suppose that \(x\in G_{\alpha}\) has prime order $r$ and fixes exactly four points.
	Suppose further that \(r\) divides \(q^2+\varepsilon
	q+1\) and let \(x \in R\in \Syl_r(G)\).
	By Theorem~6.5.3 in~\cite{GLS3} we know that \(N_G(R)\)	lies in a Frobenius group \(H\) that is a maximal subgroup of \(G\) and has order \(3\cdot \frac{q^2+\varepsilon q+1}{(3,q-\varepsilon)}\), with Frobenius kernel~\(K\) of order
	\(\frac{q^2+\varepsilon q+1}{(3,q-\varepsilon)}\). As \(r \ge 5\), it follows that \(R \le K\). Moreover, $K$ is nilpotent (see for example Theorem~10.3.1
	in~\cite{Gorenstein1980}), so we even have that
	\(R\) is characteristic in \(K\) and \(N_{G}(R)=F\). This means that
	\(G_{\alpha}\leq K\).
	Then Lemma~\ref{numberfpt} gives that  $4=|\FO(x)|$ is
	divisible by \(\frac{|N_{G}(\langle x \rangle)|}{|G_{\alpha}|}\).
	We recall that \(|G_{\alpha}|\) is not divisible by \(3\), but \(|N_{G}(\langle x \rangle)|=|F|\) is, and therefore the number of fixed points of \(x\) is divisible by \(3\). This is false.
\end{proof}

%%%%%%%%%%%%%%%

\section{The case where some involution fixes four points and G has sectional 2-rank at most 4}
%TODO: Ich bin dafür den Titel nur
%\section{The case where G has sectional 2-rank at most 4}
%zu nennen, weil der Rest dann im Fall steht

This section treats the situation from Theorem \ref{3and2}(e) in the special case where $G$ is simple.
Then the main result in \cite{GoHa1974} applies and we see that the alternating groups, except for $\Alt_6$, and also all sporadic groups except for the Lyons group cannot occur, by Remark 4.1.
Some specific groups have already been discussed, see Table \ref{SmallGroupsTable}, as have the series
$\PSL_2(q), \Sz(q), \PSL_3(q)$ and $\PSU_3(q)$
in general.

It remains to consider the series $\PSp_4(q)$, G$_2(q)$, $^2$G$_2(q)$,
$^3$D$_4(q), \PSL_4(q), \PSU_4(q), \PSL_5(q)$ or $\PSU_5(q)$ (where $q$ is
an odd prime power, with further restrictions),
and the sporadic group Ly.

\begin{lemma}\label{Inv-lem:TB2}
Suppose that Hypothesis \ref{4fix} holds, let
$\alpha\in\Omega$ and suppose that $t\in G_\alpha$ is an involution with exactly four fixed points on $\Omega$.
Then $\abs\Omega$ is even, in particular $G_\alpha$ does not contain a Sylow $2$-subgroup of $G$, and all $2$-elements of $G_\alpha$ have an even number of fixed points on $\Omega$.

If $t$ is $2$-central, then $\abs{C_G(t):C_{G_\alpha}(t)}\in\{2,4\}$.
\end{lemma}

\begin{proof}
The first statement follows from the fact that $t$ fixes exactly four points, that $G$ acts with fixity 4 and that, therefore, all
remaining orbits of $t$ on $\Omega$ must have length 2.
This implies that no point stabilizer contains a Sylow $2$-subgroup of $G$, and it also implies that all $2$-elements have an even number of fixed points.

Finally, we suppose that $t$ is $2$-central and we let $T \in \syl_2(G)$ be such that $t \in Z(T)$.

Now $\abs{C_G(t):C_{G_\alpha}(t)} \le 4$ by
Lemma \ref{Inv-lem:TB}~(c) and $T \le C_G(t)$.
Thus, if this index is 1 or 3, then $T \le G_\alpha$, contrary to the first statement.
\end{proof}

\begin{lemma}\label{Inv-lem:2G2}
Suppose that
$q$ is a power of $3$, $q \ge 27$, and that $G=\,^2G_2(q)$.
Then there is no set $\Omega$ such that $(G,\Omega)$ satisfies Hypothesis \ref{4fix} and some involution fixes exactly four points.
\end{lemma}

\begin{proof}
Assume otherwise and let $\alpha \in \Omega$ and $t \in G_\alpha$ be such that
$o(t)=2$ and $t$ fixes exactly four points.
Then Section 4.5.3 in \cite{Wilson} yields that $C_G(t)$ has structure $\langle t \rangle \times \PSL_2(q)$, and we also note that $G$ only has one class of involutions.

Lemma \ref{Inv-lem:TB2} gives that
$|C_G(t) : C_{G_\alpha}(t)| \in \{2,4\}$. As $q \ge 27$, we know that $\PSL_2(q)$ does not have any subgroups of
index 2 or 4. Moreover $t \in C_{G_\alpha}(t)$ and this forces $C_G(t) \le  G_\alpha$, contrary to Lemma \ref{Inv-lem:TB2}.
\end{proof}

Now we need some more specific preparation before we can handle the remaining groups.

\begin{hyp}\label{GoHaHyp}
In addition to Hypothesis \ref{4fix}, suppose that $G$
is a simple Lie type group that
satisfies Theorem \ref{3and2}~(e).
Let $q$ be a power of $3$ and suppose that $G$ is defined over $\GF(q)$ and that

	$G \notin \{\PSL_2(q),\Sz(q), \PSL_3(q),\PSU_3(q), \,^2G_2(q)\}$.

%Moreover, let $\alpha\in\Omega$ and let $t\in G_\alpha$ be an involution such that $\Delta:=\FO(t)$ has size $4$.

%As in parts of Remark \ref{Inv-rem:Ser}, we will now refer to notation from \cite{GLS3} again, specifically Definition 2.3.1 and Remark 2.3.3. We let $\Sigma$ denote a root system for $G$ if $G$ is of untwisted type, with fundamental system $\Pi$, and we let $\tilde{\Sigma}$ and $\tilde{\Pi}$ be the corresponding systems if $G$ is of twisted type. In the twisted case, we denote by $\hat{\Sigma}$ and $\hat{\Pi}$ the reduction modulo the equivalence relation explained in Definition 2.3.1 in \cite{GLS3}. We suppose that $G$ has Lie rank at least 2, we denote the subsets of positive roots of $\Sigma$, $\tilde{\Sigma}$ and $\hat{\Sigma}$ by $\Sigma^+$, $\tilde{\Sigma^+}$ and $\hat{\Sigma^+}$ and we let $S \in \syl_3(G)$. For all $\hat{\alpha} \in \hat{\Sigma}^+$ we let $X_{\hat{\alpha}} \le S$ be the root subgroup that belongs to $\hat{\alpha}$ in $S$. Finally, let $P \le S$ be a subgroup of index 3.
\end{hyp}

\begin{lemma}\label{GoHaq=3}
Suppose that Hypothesis \ref{GoHaHyp} holds and that $q=3$.
Then there is no involution in $G$ that fixes exactly four points.
\end{lemma}

\begin{proof}
The following cases can be treated with \texttt{GAP} (see Table \ref{SmallGroupsTable} and Appendix):
$\PSL_4(3), G_2(3)$, and $\PSp_4(3)$. None of these groups exhibit a fixity four action with some involution fixing four points.
As $q=3$ and $q \equiv 1$ modulo $4$ for the series $\PSU_5(q)$, under our hypothesis,  the group $\PSU_5(3)$ does not need to be considered here.
So we are left with $\PSL_5(3)$ and $\,^3D_4(3)$.

\textbf{Case 1:} $G=\,^3D_4(3)$.

We recall that $G$ has only one class of involutions (Table 4.5.1 in \cite{GLS3}).
So if we assume that $\alpha \in \Omega$ and $t \in G_\alpha$ is an involution with exactly four fixed points, then
Section 4.6.5 in \cite{Wilson} gives that $C_G(t)$ has structure $C_2(\PSL_2(3^3) \times \PSL_2(3)). C_2$.
Also, Lemma \ref{Inv-lem:TB2} tells us that $|G_\alpha|$ is divisible by $|\PSL_2(3^3)| \cdot |\PSL_2(3)|=2^4 \cdot 3^4 \cdot 7\cdot 13$.
Comparing this to $|G| = 2^6 \cdot 3^{12} \cdot 7^2 \cdot 13^2 \cdot 73$ forces $G_\alpha$ to have order divisible by
$2^4 \cdot 3^4 \cdot 7^2 \cdot 13^2$, with Lemma \ref{Inv-lem:TB}.
But then $G$ does not have any maximal subgroup that could contain $G_\alpha$ (see for example Theorem 4.3 in \cite{Wilson}).

\textbf{Case 2:} $G=\PSL_5(3)$.

This time $G$ has two classes of involutions (Table 4.5.1 in \cite{GLS3}), and the corresponding centralizer orders are
$2^8 \cdot 3^4 \cdot
13$ and $2^9 \cdot 3^6 \cdot 5 \cdot 13$.
We will also need that $|G|= 2^9 \cdot 3^{10} \cdot 5 \cdot 11^2 \cdot 13$.
Assume for a contradiction that $\alpha \in \Omega$ and that $t \in G_\alpha$ is an involution with exactly four fixed points.
Let $Q \in \syl_3(C_G(t))$. Then $Q$
fixes $\FO(t)$ or it fixes a unique point in $\FO(t)$. In both cases we deduce that all $3$-elements in $G$ fix exactly four elements or exactly one, which implies that $|\Omega|$ is not divisible by $3$ and that, therefore, $G_\alpha$ contains a Sylow $3$-subgroup of $G$. Without loss $P \le G_\alpha$.

Using Theorem
2.6.5 (d) in \cite{GLS3} and Lemma \ref{Inv-lem:TB} we deduce that
$G_\alpha$ contains a subgroups of index at most 4 of each maximal parabolic subgroup of $G$ that contains $P$, and we know the possible types of subgroups because they are maximal subgroups (see for example Theorem 2.6.7 in \cite{GLS3}).
This information can be calculated with \texttt{GAP}, and it leads to a contradiction because $|G_\alpha|$ is too large.
\end{proof}

\begin{lemma}\label{Inv-lem:Ser}
In addition to Hypothesis \ref{4fix}, suppose that $G$
is a simple Lie type group that
satisfies Theorem \ref{3and2}~(e).
Let $p$ be an odd prime and $n\in\N$ be such that
	$q:=p^n\ge 3$, that $G$ is defined over $\GF(q)$ and that
	$G \notin \{\PSL_2(q),\Sz(q), \PSL_3(q),\PSU_3(q), {^2G}_2(q)\}$.

Moreover, let
$\alpha\in\Omega$ and let $t\in G_\alpha$ be an involution such that $\Delta:=\FO(t)$ has size $4$.

	Then $q \le 3$.
\end{lemma}

\begin{proof}
We assume otherwise and we set $C:=C_G(t)$.
Next we recall the possibilities from the introductory paragraph, including further details from the main result of \cite{GoHa1974}:

$\PSL_4(q)$ ($q\not\equiv 1$ modulo 8),

$\PSU_4(q)$ ($q\not\equiv 7$ modulo 8),

$\PSL_5(q)$ ($q\equiv -1$ modulo 4),

$\PSU_5(q)$ ($q\equiv 1$ modulo 4),

$\PSp_4(q)$, $G_2(q)$ (where $q=3^{2n+1}$) or $^3D_4(q)$.

Here we left out the groups that are excluded by our hypothesis.

%Before we start with our main arguments, we establish that $\abs{C:C_\alpha}\in\{2,4\}$ by Lemma \ref{Inv-lem:TB2}.
%The reason is that all the groups above have only one class of involutions. (See for example Table 4.3.1 and 4.5.1 in \cite{GLS3}.) We will refer to the fact that $\abs{C:C_\alpha}\in\{2,4\}$ as $(\ast)$.

\smallskip
\textbf{Claim 1:} $p \in \pi(G_\alpha)$.

\begin{proof}
Our hypotheses and Lemma \ref{Inv-lem:TB2} yield that $\abs{C:C_\alpha}\in\{1,2,4\}$.
We use Table 4.5.1 in \cite{GLS3}, where we see that $C$ has a section %$K$
isomorphic to a finite group of Lie type or a product of several such groups over $\GF(q)$ and in dimension at least 2. So it follows that
$p \in \pi(C)$. (For more details we refer to Table \ref{Inv-tab:Ser}.) As $p$ is odd, we see that $p$ divides $|C_\alpha|$ and hence $p \in \pi(G_\alpha)$.
\end{proof}

\textbf{Claim 2:}  If $p\ge 5$, then
$G_\alpha$ contains a Sylow $p$-subgroup of $G$ and a subgroup of index at most 4 of a maximal parabolic subgroup.

\begin{proof}
The first statement follows from Claim 1 and Lemma \ref{Inv-lem:TB}. If $P \in \syl_p(G_\alpha)$, then
Theorem 2.6.5 (d) in \cite{GLS3} tells us that there is a subgroup $U \le P$ such that $N_G(U)$
is a maximal parabolic subgroup of $G$. Then
$|N_G(U):N_{G_\alpha}(U)| \le 4$ by Lemma \ref{Inv-lem:TB}, which is the second statement.
\end{proof}

\textbf{Claim 3:}  If $p=3$, then
$G_\alpha$ contains a Sylow $p$-subgroup of $G$ and a subgroup of index 1,2 or 4 of a maximal parabolic subgroup of $G$.

\begin{proof}
 Using Claim 1 we let $P \in \syl_3(G)$ be such that $1 \neq P_\alpha \in \syl_3(G_\alpha)$.
Next we go through the cases of Lemma \ref{syl3} and we see that (a) does not hold. Moreover,
 $\abs{P}\ge q^3\ge 27>9$ by \cite{Wilson} (and by hypothesis) and $P$ does not have maximal class, as we can see in Table \ref{Inv-tab:Ser}. Therefore Lemma \ref{syl3} (d) or (e) holds, and (e) immediately implies our first statement.
We recall that $|C:C_\alpha| \in \{1,2,4\}$, which means that $C_\alpha$ even contains a Sylow $3$-subgroup of $C$. Now every $3$-element in $C$ fixes $\alpha$ and stabilizes $\Delta$, and hence it fixes exactly one or exactly four points on $\Omega$ in total.
This forces $|\Omega| \equiv 1$ modulo 3, which is not compatible with
Case (d) of Lemma \ref{syl3}.
We are left with Case (e) of Lemma \ref{syl3}, where our first assertion holds, and the second statement follows with
Theorem 2.6.5(d) in \cite{GLS3} and Lemma \ref{Inv-lem:TB}.
\end{proof}

We use the claims 2 and 3 and we let
$R\ug G$ be a maximal parabolic subgroup of $G$ such that
$|R:R_\alpha|\le 4$ and in particular $G_\alpha$ contains a Sylow $p$-subgroup $P$ of $G$.
Then Table \ref{Inv-tab:Ser} gives that $G$ has at least two different conjugacy classes of maximal parabolic subgroups.
Using Sylow's Theorem we let $U_1,U_2 \le P$ be such that $N_G(U_1)$ and $N_G(U_2)$ are non-conjugate
maximal parabolic subgroups of $G$. By Theorem 2.6.7 in \cite{GLS3} they are maximal subgroups of $G$.
For both, their intersection with $G_\alpha$ has index at most 4 in $G_\alpha$, and this forces $G_\alpha=G$, which is a contradiction.
\end{proof}

We have now treated all Lie type groups that occur in Case (e) of Theorem \ref{3and2} and that have not been analyzed before.
Finally, we look at the Lyons group.

\begin{lemma}\label{Lyons}
Let $G$ be the sporadic Lyons group.
Then there is no set $\Omega$ such that $(G,\Omega)$ satisfies Hypothesis \ref{4fix} and some involution fixes four points.
\end{lemma}

\begin{proof}
Inspection of the character table of $G$ in the ATLAS \cite{Atlas} shows that $G$ has only one class of involutions,
in particular every involution fixes exactly four points.
Let $\alpha \in \Omega$ be one of the fixed points of some involution $t \in G$.
Then $|C_G(t)|=2 \cdot |\Alt_{11}|$, so
Lemmas 5.3 and \ref{Inv-lem:TB}(a) imply that $G_\alpha$ contains a Sylow subgroup of $G$
for the primes $5, 7$ and $11$.
Now Lemma \ref{Inv-lem:TB}(c) yields that $|G_\alpha|$ is divisible by $2^6 \cdot 3^4 \cdot  5^6 \cdot 7 \cdot 11$
as well.
But then there is no maximal subgroup of $G$ that contains $G_\alpha$, which is a contradiction.
\end{proof}

%%%%%%%%%

\section{The case where point stabilizers have odd order divisible by 3}

Here we consider Case %(b) of Theorem \ref{3and2}.
(4) of Theorem~\ref{3and2simple}.
If $G$ is non-abelian simple, then the CFSG gives
a list of possible groups for the situations where
$G$ has a strongly $3$-embedded subgroup or small Sylow $3$-subgroups of prescribed structure.
The simple groups of Lie type with small Lie rank that occur have already been discussed.
%Some groups have been listed in Table \ref{SmallGroupsTable}, so the only ones left are
% \({}^2\textrm{G}_2(3^{2n+1})\), where $n \in \N$,  \(\textrm{G}_2(q)\), where $q$ is a prime power (with further specifications),
% \({}^2F_4(2)'\) or \({}^2\textrm{F}_4(2^{2n+1})\), where \(n\geq 1\),
% \(\PSp_4(q)\), \(\PSL_4(q)\), \(\PSU_4(q)\), \(\PSL_7(q)\), \(\PSU_7(q)\), \(\PSp_6(q)\), \(\POm_7(q)\) or \(\POm^-_8(q)\), where $q$ is a prime power (with further restrictions).
 In this section we
 %investigate all these groups and
%% that exhibit a fixity 4 action with the given
 %point stabilizer structure. We
 work under the following further hypothesis:

\begin{hyp}\label{odd3}
Suppose that Hypothesis \ref{4fix} holds and that the order of the point stabilizers in $G$ is odd and divisible by 3. Suppose further that
	$G \notin \{\PSL_2(q), \PSL_3(q),\PSU_3(q)\}$. \\
\end{hyp}
%TODO: Wollen wir hier G einfach mit aufnehmen?

%In this section let $H$ be a nontrivial element-wise stabilizer of a set $\Delta \subset \Omega$ of four elements and let $p \in \pi(H)$.

We start with a helpful technical lemma.

%P: Das folgende Lemma hat eigentlich nichts mit Kap 7 zu tun, außer, dass dies der einzige Ort ist, an dem
%wir es benutzen. Deswegen habe ich es in Kapitel 4 verschoben.
%\begin{lemma}\label{helpFrattini}
%Let $A$ be a soluble group and $p,q \in \pi(A)$ be two different prime divisors of the order of $A$. Then
%there is a prime $r \in \{p,q\}$ and a nontrivial $r$-subgroup $R$ of $A$ that is normalized by
%some Sylow $s$-subgroup of $A$, where $\{r,s\} = \{p,q\}$.
%\end{lemma}
%\begin{proof} As $A$ is soluble there is a chief series $1 = A_m \unlhd \cdots \unlhd A_1 = A$
%of normal subgroups $A_i$ of $A$, such that $A_{i}/A_{i+1}$ is an elementary abelian group of prime
%power order. Let $j$ be maximal such that $r$ divides $|A_j|$ for some $r \in  \{p,q\}$, and let
%$R$ be a Sylow $r$-subgroup of $A_j$. Then according to the Frattini argument $A = A_j N_A(R)$.
%Let $s$ be such that $\{r,s\} = \{p,q\}$. Then  $s$ does not divide $|A_j|$ and the assertion holds.
%\end{proof}

\begin{lemma}\label{oddpoint}
Suppose that Hypothesis \ref{odd3} holds for the simple group $G$ and
let $H$ be a nontrivial four point stabilizer. Then the
following hold:
\begin{enumerate}
    \item[(a)] The point stabilizers in $G$ are soluble.
    \item[(b)] One of the following is true:
    \begin{enumerate}
        \item[(i)] $G$ contains a strongly $3$-embedded subgroup or
        \item[(ii)] $3 \in \pi(N_G(H))$.
    \end{enumerate}
\end{enumerate}
\end{lemma}
\begin{proof}
Assertion (a) holds by \cite{FT}, because the point stabilizers have odd order by hypothesis.

For (b) we suppose that (ii) does not hold and we let $\alpha \in \Omega$ be such that $H \le G_\alpha$. Let $p \in \pi(H)$. As $3$ does not divide $|N_G(H)|$ and $H$ has odd order, we see that $p \geq 5$. We recall that $G_\alpha$ is soluble,
and that moreover $3,p \in \pi(G_\alpha)$. Application of  Lemma~\ref{helpFrattini}
yields that $p$ divides $|N_G(X)|$ for some nontrivial $3$-subgroup $X$ of $G_\alpha$ or that \(3\) divides \(|N_G(Y)|\) for some nontrivial \(p\)-subgroup \(Y\) of \(G_{\alpha}\).
In the latter case, Lemma~\ref{Inv-lem:TB}\,(a) implies that we may suppose that \(Y\leq H\). Therefore there is a nontrivial \(3\)-element in \(N_G(Y)\). %stabilizing \(\FO(Y)=\Delta\).
By Lemma~\ref{Inv-lem:TB}\,(b) this element lies in \(N_G(H)\), which contradicts our assumption. Hence $p$ divides $|N_G(X)|$ for some nontrivial $3$-subgroup $X$ of $G_\alpha$.
As $3$ does not divide $|N_G(H)|$, the action of every nontrivial $p$-element on $X$ is fixed point free.
We recall that $p \geq 5$ and then it follows that $|Z(X)| \ge 3^3$.
Now by checking the cases in Theorem \ref{3and2simple} we derive that $G$ contains a strongly $3$-embedded subgroup, which is
(i).
\end{proof}

Next we consider the case that $G$ contains a strongly 3-embedded
subgroup.

\subsection{The strongly 3-embedded case}

\begin{lemma}\label{CFSGstr3emb}
    Let \(G\) be a finite simple group with a strongly \(3\)-embedded subgroup. Then \(G\) is one of the following groups:
    \begin{enumerate}
        \item \(\PSL_2(q)\), where $q$ is a prime power, %and $q=9$ is the only power of $3$ that occurs,
        \item \(\PSL_3(q)\), where $q = 4$ or $q$ is a prime power and $q \equiv -1 \mod{3}$,
        \item \(\PSU_3(q)\), where $q$ is a $3$-power or $q$ is a prime power and  $q\equiv 1 \mod{3}$,
        \item \({}^2\textrm{G}_2(q)\), where $q>3$ is a power of $3$,
        \item \(\Jan_1\) or
        \item \(\Mat_{11}\).
    \end{enumerate}
\end{lemma}

\begin{proof}
      As a consequence of the classification of finite simple groups together with (10-1) in~\cite{GL} and Theorem~3.3.3 in~\cite{GLS3}, the finite non-abelian simple groups with nontrivial cyclic Sylow \(3\)-subgroups are   \(\PSL_2(q)\), where $q$ is not a $3$-power, \(\PSL_3(q)\), where $q \equiv -1 \mod{3}$ and \(\PSU_3(q)\), where $q \equiv 1 \mod{3}$.
Then the statement of the lemma follows from Theorem 7.6.1 in \cite{GLS3}.
\end{proof}

\begin{lemma}\label{2G2}
	
	Let \(n \in \N\), let $q:=3^{2n+1}$ and let \(G={}^2\textrm{G}_2(q)\).
If $\Omega$ is a set such that $(G, \Omega)$ satisfies Hypothesis \ref{odd3}, then
the point stabilizers are Frobenius groups of structure
% more precisely a semidirect product of a group of
\(q^{3} \rtimes C_{\frac{q-1}{2}}\).

%Conversely, if $H \le G$ is a Frobenius group of structure \(3^{3(2n+1)} \ltimes C_{\frac{3^{2n+1}-1}{2}}\), then the action of \(G\) acts transitively, faithfully and with fixity $4$ on \(G/G_{\alpha}\).
\end{lemma}

\begin{proof}
Let $P$ be a Sylow $3$-subgroup of $G$. By part (2) of the Theorem in~\cite{Ward}, $P$ is of class $3$, \(Z(P)\) is elementary abelian of order \(q\) and
\(P'=\Phi(P)\) is an elementary abelian subgroup of order \(q^2\) containing \(Z(P)\).
Let \(\Omega\) be such that $(G, \Omega)$ satisfies Hypothesis \ref{odd3},  and
let $\Delta$ be the union
of the $P$-orbits on $\Omega$ that are of size at most $3$.
As \(|P|=q^3\geq 3^9\) and \(|G_{\alpha}|\) is odd, case (a)\,(ii) of Theorem~\ref{3and2} holds.
Therefore $1 \leq |\Delta| \leq 4$.
As $P$ is of class $3$, it is not of maximal class, and Case (e) of Lemma~\ref{syl3} holds; that is  $P \leq G_\alpha$ for some
\(\alpha \in \Delta\).
Then the element-wise stabilizer  of $\Delta$ in $P$ is a subgroup $Q$
of index at most $3$ in $P$.

According to Theorem~4.2 in~\cite{Wilson}, $N_G(P) = P \rtimes C$ is
a maximal subgroup of $G$, where $C$ is a
cyclic group of order $q-1$.
As \(Z(P)\leq P'\), every subgroup of index at
most $3$ in $P$ contains $Z(P)$. %, which is elementary abelian of order $q$.
This shows that $Z(P) \leq Q$ and therefore \(\FO(Z(P))=\Delta\).

Further notice that $(q-1)/2$ is odd, as $q$ is an odd prime power of $3$.
Since \(C\leq N_G(Z(P))\) this shows that the subgroup of $U$ index $2$ in $C$ also fixes \(\Delta\)
element-wise, and $N_G(P)$ is the stabilizer of $\Delta$. This yields the
assertion.
\end{proof}

Notice that the set $\Delta$ considered in the proof of Lemma~\ref{2G2} has size $2$ and that the subgroup $U \cong C_{(q-1)/2}$ is a nontrivial four point stabilizer.

\begin{lemma}\label{str3emb}
	
	Let \(G\) be a finite simple group with strongly \(3\)-embedded subgroup, suppose that
	$G$ satisfies Hypothesis \ref{odd3}, and let $\alpha \in \Omega$. Then
	 one of the following is true:
	\begin{enumerate}
		\item \(G\cong \Jan_1\) and \(G_{\alpha}\) is cyclic of order \(15\).
		\item \(G \cong {}^2\textrm{G}_2(q)\), where \(n\geq 1\), $q=3^{2n+1}$ and \(G_{\alpha}\) is a Frobenius group
		of structure \(q^{3} \rtimes C_{\frac{q-1}{2}}\).
	\end{enumerate}
\end{lemma}

\begin{proof}
%This is a consequence of \ref{CFSGstr3emb}, where we have already omitted the groups that are excluded in Hypothesis \ref{odd3}.
Omitting the groups that are excluded in Hypothesis~\ref{odd3}, Lemma~\ref{CFSGstr3emb} shows that \(G\) is either \({}^2\textrm{G}_2(q)\), \(\Jan_1\) or \(\Mat_{11}\).
The specific groups can be found in Table \ref{SmallGroupsTable} and the rest follows from Lemma \ref{2G2}.
\end{proof}

As we discussed this case completely we assume in the remainder of the section the following.

\begin{hyp}\label{odd3simple}
Suppose that Hypothesis \ref{odd3} holds, that \(G\) is simple and that \(G\) does not contain a strongly \(3\)-embedded subgroup.
\end{hyp}

We notice that groups with Sylow \(3\)-subgroups of order \(3\) have cyclic Sylow \(3\)-subgroups and
hence a strongly \(3\)-embedded subgroup. Therefore they have already been dealt with in the previous results, and this means that
we have now covered (4)(a) and (b) from Theorem \ref{3and2simple}.
As a next step we show that the cases of extra-special Sylow $3$-subgroups of order $27$ of exponent $3$ or wreath products $C_3 \wr C_3$ cannot happen.

%\subsection{Sylow $3$-subgroups that are either extra-special of order $27$ or wreath products of order $81$}
\subsection{Sylow 3-subgroups that are either extra-special of order 27 or wreath products of order 81}

\begin{lemma}\label{extraspecialWreath}
Let $G$ be a simple group such that a Sylow $3$-subgroup $P$ is extra-special  of order 27 and of exponent $3$ or
a wreath product $C_3 \wr C_3$. Then there does not exist a set $\Omega$ such that $(G,\Omega)$ satisfies Hypothesis \ref{odd3simple}.
\end{lemma}

\begin{proof} Assume otherwise and let $(G, \Omega)$ be such a pair.
Let $P$ be a Sylow $3$-subgroup of $G$ and let
$\Delta$ denote the union of the $P$-orbits of length at most $3$.
Our hypotheses imply that the cases (a), (b) and (d) of Lemma \ref{syl3} do not hold.
If (e) holds, then Proposition \ref{3semireg} gives a contradiction to our hypotheses that point stabilizers have odd order and that $G$ does not have a strongly $3$-embedded subgroup.
The only remaining case is (c), which means that Hypothesis \ref{hyp4} is satisfied and that we can apply Proposition \ref{maxKlasse}. The only possibilities that are compatible with our hypotheses are (e) and (f).
In Case (e) we have that $|\Delta| = 0$, $P$ is extra-special of order $27$ and $|G_\alpha|_3 = 3$ for some $\alpha \in \Omega\setminus{\Delta}$, and in Case (f) we have that
$|\Delta| = 3$, $P \cong C_3 \wr C_3$ and the Sylow $3$-subgroups of point stabilizers are elementary abelian of order $27$.

In both cases, every nontrivial $3$-element of $G$ fixes exactly 0 or 3 points in $\Omega$.

Let $H$ be a nontrivial four point stabilizer and let $\alpha \in \FO(H)$.
Then Hypothesis \ref{odd3} and the previous paragraph imply that $|H|$ is coprime to $6$.
Let $p \in \pi(H)$ and $T \in \syl_p(H)$. We know that $N_G(H)$
is strongly $p$-embedded in $G$ by Lemma \ref{Inv-lem:TB}(h), and then \cite[Theorem 7.3.1]{GLS3} yields that either $T$ is cyclic or  $G$ is one of the groups listed in
that theorem.  As $p > 3$, the only possibilities are $\Alt_9$, $\Alt_{10}$ or $\Alt_{11}$, where the Sylow $3$-structure is not as in our hypothesis, or
$|G|_3 > 3^4$. All this is not possible, and therefore $T$ is cyclic. This implies that $\aut(T)$ is abelian. We keep this prime $p \in \pi(H)$.

Next we make an observation that we will quote several times in the proof:

$(\ast)$~~
If $Y\le
G_\alpha$ is a $3$-group and $1 \neq U \le H$ is a $p$-group, then $C_U(y) =1$ for all $y \in Y^\#$ and
$C_Y(u) =1$ for all $u \in U^\#$.

Assume otherwise and let $1 \neq u \in C_U(y)$.
Then $u$ fixes the three fixed points of $y$ because $p \ge 5$, and conversely $y$ stabilizes the set $\FO(u)$ of size $4$. This is impossible because $y$ cannot have a fourth fixed point.
The other statement holds for the same reason.

By hypothesis $G_\alpha$ has odd order, hence it is a soluble group and Lemma \ref{helpFrattini} applies. We use it for the primes $3$ and $p$, both in $\pi(G_\alpha)$, and we have two cases:

(1) There is a nontrivial $p$-subgroup $U$ of $G_\alpha$ that is normalized by a Sylow $3$-subgroup $Q$ of $G_\alpha$ or

(2) there is a nontrivial $3$-subgroup $Y$ of $G_\alpha$ that is normalized by a Sylow $p$-subgroup $T$ of $G_\alpha$.

For both situations we keep in mind that $p \in \pi(H)$ and $p \ge 5$, so $H$ contains a Sylow $p$-subgroup of $G$ by Lemma \ref{Inv-lem:TB}(a), in particular of $G_\alpha$. This means that we can choose $U \le H$ in (1) and $T \le H$ in (2).

We begin with Case (2), so $Y \le G_\alpha$ is a nontrivial $3$-group normalized by $T \in \syl_p(H)$. Without loss $Y \le P$.

Then $Y\ltimes T$ is a Frobenius group by ($\ast$). As $p \ge 5$, this forces $|Y| \ge 9$ and hence $|G_\alpha|_3 \ge 9$, which excludes Case (e) of Proposition \ref{maxKlasse}.
We conclude that (f) holds, so $P \cong C_3 \wr C_3$, $|\Delta| = 3$, $p = 13$ and $H = T$. The action of $T$ on $Y$ also forces $Y \in \syl_3(G_\alpha)$, and then $|P:Y|=3$ and $P, T \le N_G(Y)$.
As $|\FO(Y)| = 3$, our hypothesis forces $N_G(Y)$ to have odd order. Another application of Lemma \ref{helpFrattini}, this time to the soluble group $N_G(Y)/Y$.
If some nontrivial $3$-subgroup of $N_G(Y)/Y$ normalizes a nontrivial $p$-subgroup, then we consider pre-images in $N_G(Y)$ and the fact that all Hall $\{3,p\}$-subgroups of $N_G(H)$ are conjugate. Together with ($\ast$) this gives a contradiction.

Therefore $T$ normalizes a Sylow $3$-subgroup of $N_G(Y)/Y$, so without loss it normalizes $P/Y$ and hence its full pre-image $P \le N_G(Y)$.
In particular $T$ normalizes $Z(P)$, which has order $3$, and this contradicts ($\ast$) because $Z(P) \le G_\alpha$ by Lemma \ref{tiny}.

In Case (1) we have a nontrivial $p$-subgroup $U$ of $H$ that is normalized by a Sylow $3$-subgroup $Q$ of $G_\alpha$. Then $Q \le N_G(H)$ by Lemma \ref{Inv-lem:TB}(b), and ($\ast$) yields that
$H\rtimes Q$ is a Frobenius group. In particular $H$ is nilpotent. Therefore we may suppose that $U \in \syl_p(H)$, so without loss $U=T$.
As $T$ is cyclic, it follows that $N_H(T) = C_H(T) = H$ and our previous observation yields that $|Q|=3$, that $P$ is extra-special of order $27$ and that $|\Delta| = 0$.
%As a next step we prove that $N_G(U)$  is not divisible by $3$ for each nontrivial subgroup $U$ of $H$.  As $N_G(U) \leq N_G(H)$, it is enough to show that $3$ does not divide $N_G(H)$. Assume to the contrary that   $3$  divides $N_G(H)$.
Lemma \ref{Inv-lem:TB} (e) and (a) gives that $N_G(H)$ acts transitively
on $\FO(H)$. By hypothesis $G_\alpha$ has odd order, so $N_G(H)$ induces a transitive subgroup of $\Alt_4$. Moreover $3 \notin \pi(H)$, so we conclude that $N_G(H)$ contains 3-elements that induce 3-cycles on $\FO(H)$ and hence it
induces the group $\Alt_4$ on $\FO(H)$. Let $M \le N_G(H)$ be the full pre-image of $O_2(N_G(H)/H)$. Then $M \unlhd N_G(H)$ acts transitively on
$\FO(H)$, $T \in \syl_p(M)$ and a Frattini argument gives that
$N_G(H)=M \cdot N_G(T)$, together with Lemma \ref{Inv-lem:TB}(b).
In particular $N_G(T)/H$
contains a subgroup isomorphic to $\Alt_4$.

Let $V \in \Syl_2(N_G(T))$. As $Q \le N_G(T) \le N_G(H)$ and $N_G(H) = (H \cdot V)\cdot Q$ by the previous paragraph, we see that $N_G(T)=(N_H(T) \cdot V) \cdot Q=(H \cdot V)\cdot Q$, as $H$ is nilpotent.
Now $Q$ acts fixed-point-freely on $H \rtimes V$ by the structure of $\Alt_4$.
The group $N_G(H)$ is a Frobenius group with Frobenius kernel $H \cdot V$, so $H \cdot V$ is nilpotent and hence
$[H,V] = 1$.  We recall that $|Q|=3$, so  $N_G(T)=C_G(T) \rtimes Q = (H\cdot V) \rtimes Q$.

Suppose that $X$ is a nontrivial $p$-subgroup of $G$. Without loss $X\le H$, because $H$ contains a Sylow $p$-subgroup of $G$, and then $N_G(X)\le N_G(H)=(H \times V) \rtimes Q$, which implies that $N_G(X) \cap C_G(V)$ has a normal $p$-complement.

So if we consider $C:= C_G(V)$, then
for each nontrivial $p$-subgroup $X \le C$ we have that $N_C(X)$ has a normal $p$-complement, and then
Burnside's $p$-Complement Theorem gives
that $C$ has a normal $p$-complement as well. Let $N:=O_{p'}(C)$ and assume for a contradiction that $3$ divides $|N|$. As $C=N \cdot T$, coprime action gives that $T$ normalizes a Sylow $3$-subgroup $Y$ of $N$. We also recall that $V$ does not centralize
a Sylow $3$-subgroup of $G$, and therefore $|Y| \le 9$ and $[T,Y]=1$. Then $Y$ stabilizes $\FO(T)$, which has size $4$, so $Y$ must fix a point, contrary to ($\ast$). We conclude that $N$ is  a $3'$-group.
We also note that $V$ stabilizes $\FO(T)$ and acts fixed point freely on it. Together with the fact that $N_G(T)/T$ has a subgroup isomorphic to $\Alt_4$, this implies that $V$ is a fours group.

Let $V \le S \in \syl_2(C)$. Then $S \le N$ because $C=N \rtimes T$, so the coprime action of $T$ on $N$ gives that we may choose $S$ to be normalized by $T$.
Moreover, the factorization $N_G(T)=(H \times V) \rtimes Q$ shows that $Q$ normalizes $V$ and hence $Q \le N_G(C) \le N_G(N)$. In the previous paragraph we proved that $N$ is  a$3'$-group, and then the coprime action of $Q$ on $N$ yields that we may choose
$S$ to be $Q$-invariant as well, so that $Q$ normalizes $S \rtimes T$.
%Then $S \cdot T$ is a $\{2,p\}$-Hall subgroup of $C$, and since $C$ is soluble, all $\{2,p\}$-subgroups of $C$ are conjugate to $S \cdot T$. All the subgroups $S\ltimes T$ are conjugate in $C$, which yields that we may assume that $Q$ normalizes this subgroup (the number of conjugates is not divisible by $3$). The fact that $|N_G(Q)|$ is coprime to $2p$
We recall that $|\FO(Q)|=3$. If $s \in S^\#$ is centralized by $Q$, then
$s$ stabilizes $\FO(Q)$ and hence it must fix a point, contrary to our hypothesis that point stabilizers have odd order. Together with ($\ast$)
this shows that $(S \cdot T)\rtimes Q$ is a Frobenius group with Frobenius kernel $S \cdot T$, and in particular $[S,T]=1$.
We conclude that $S \le N_G(T)=N_G(T) = (H \times V)\rtimes Q$ and hence $S=V$.
As $V$ stabilizes $\FO(H)$, a set of size $4$, this yields two cases: $V \in \syl_2(G)$ or $N_G(V)$ contains a subgroup $D \cong D_8$. In the second cases the point stabilizers have even order, which is false. Now the main result in \cite{GW} gives that
$G\cong\PSL_2(q)$ for some prime power $q$ or $G \cong \Alt_7$.
Then a Sylow $3$-subgroup of $G$ is cyclic or elementary abelian, contrary to our hypothesis.
\end{proof}

\subsection{Elementary abelian Sylow 3-subgroups of order 9}

Applying the classification of finite simple groups we get the following result (see Proposition~(1.2) in~\cite{Koshitanielab}).

\begin{lemma}\label{CFSGelab}
    Let \(G\) be a finite simple group with elementary abelian Sylow \(3\)-subgroups of order \(9\). Then \(G\) is one of the following groups:
    \begin{enumerate}
        \item \(\Alt_6, \Alt_7, \Mat_{11}, \Mat_{22}, \Mat_{23}\) or \(\HS\),
        \item \(\PSL^\varepsilon_3(q)\), where \(q\) is congruent to \(3+\varepsilon\) or \(6+\varepsilon\) modulo \(9\),
			except \(\PSL^{-1}_3(2)\cong \PSU_3(2)\), %not simple
			\item \(\PSL^\varepsilon_4(q)\), where \(q\) is congruent to \(3-\varepsilon\) or \(6-\varepsilon\) modulo \(9\),
			\item \(\PSL^\varepsilon_5(q)\), where \(q\) is congruent to \(3-\varepsilon\) or \(6-\varepsilon\) modulo \(9\), or
			\item \(\PSp_4(q)\), where \(q\) is congruent to \(2,4,5\) or \(7\) modulo \(9\), except \(\PSp_4(2)\).
    \end{enumerate}
\end{lemma}

As $G\cong \Mat_{11}$ has a strongly $3$-embedded subgroup, we already discussed this group.
Now we present some easy observation which helps to determine the
groups satisfying Hypothesis~\ref{odd3}.

%\begin{lemma}\label{helpSL2xSL2sub}
%	Suppose that Hypothesis \ref{odd3} holds and that \(G\) is a finite simple group with an elementary abelian
%	Sylow \(3\)-subgroup  of order \(9\) that does not contain a strongly
%	$3$-embedded subgroup.  Then  $3 \in \pi(N_G(H))$, and if $|N_G(X)|$ is of even order for every subgroup $X$ of $G$ of order $3$, then $G_\alpha$ and $N_G(H)$ contain a Sylow $3$-subgroup of $G$.
%\end{lemma}
%\begin{proof}
% The first assertion is a consequence of Lemma~\ref{oddpoint}.
% Suppose that $|N_G(X)|$ is of even order for every subgroup $X$ of $G$ of order $3$. Let $X \leq N_G(H)$ be a subgroup of order $3$.
% Then by assumption $N_G(X)$ is of even order. As the stabilizer $G_\alpha$ for $\alpha \in \Delta$ is of odd order it follows
% that $|\FO(X)| \in \{2,4\}$. Assume $|\FO(X)| = 2$, and let
% $P$ be a Sylow $3$-subgroup of $G$ which contains $X$ as a subgroup.
% Then $P$ acts trivially on $\FO(X)$
% and $|\FO(P)| = 2$. We may assume $P \leq G_\alpha$. It follows that  $H$ is a subgroup whose order is
% coprime to $6$. Recall that  $p \in \pi(H)$.  As $G_\alpha$ is soluble we either get by Lemma~\ref{helpFrattini} that an  element of $H$ of order $p$ acts
% on a nontrivial $3$-subgroup or that $P$ acts on a nontrivial $p$-subgroup. As $\aut(P)$ is a $\{2,3\}$-group, we get in both cases $P$ normalizes a nontrivial $p$-subgroup. This yields $P \leq N_G(H)$.
% \end{proof}

 \begin{lemma}\label{helpSL2xSL2sub}
	Suppose that Hypothesis \ref{odd3simple} holds and that \(G\) has elementary abelian
	Sylow \(3\)-subgroups  of order \(9\).
	%Then  $3 \in \pi(N_G(H))$,
	If $|N_G(X)|$ has even order for every subgroup $X$ of $G$ of order $3$, then every point stabilizer of $G$ contains a Sylow $3$-subgroup of $G$
	and there exists a nontrivial four point stabilizer \(H\)  such that \(N_G(H)\) contains a Sylow \(3\)-subgroup of \(G\).
\end{lemma}
\begin{proof}
Let \(H\) be a nontrivial four point stabilizer in $G$ and let \(\Delta:=\FO(H)\).
 Then \(3\) divides \(|N_G(H)|\) by Lemma~\ref{oddpoint}, as and \(G\) does not contain a strongly \(3\)-embedded subgroup.

 Suppose that $|N_G(X)|$ has even order for every subgroup $X$ of $G$ of order $3$. Let $X \leq N_G(H)$ be a subgroup of order $3$ and let $\alpha \in \Delta$.
 %Then by assumption $N_G(X)$ is of even order.
 As $G_\alpha$ has odd order, every nontrivial \(2\)-element in $N_G(X)$ acts fixed-point-freely on \(\FO(X)\).
 It follows
 that $|\FO(X)| \in \{2,4\}$. Assume that $|\FO(X)| = 2$, and let
 $P$ be a Sylow $3$-subgroup of $G$ such that $X \le P$.
 %Then $P$ acts trivially on $\FO(X)$ and
 As \(P\) is abelian,
the subgroup $P$ acts on $\FO(X)$ and hence fixes its two elements. Therefore
 $|\FO(P)| = 2$,
 and every element in \(P^\#\) fixes exactly two elements in $\Omega$.
 This implies that  $|H|$ is coprime to $6$.
 Without loss we may assume $X\le G_\alpha$, which implies $\alpha \in \FO(P)$ and
 $P \leq G_\alpha$. In particular, if we let $p \in \pi(H)$, then \(p\geq 5\).
 Lemma~\ref{helpFrattini} applied to $G_\alpha$ gives two possibilities:
 Some element \(h\) of $H$ of order $p$ normalizes
a nontrivial $3$-subgroup \(Y\) of $G_\alpha$ or $P$ normalizes a nontrivial $p$-subgroup $U$ of $G_\alpha$.
 In the first case $h \in C_G(Y)$ as \(p\geq 5\). Therefore $h$ acts on $\FO(Y)$
 yielding $\FO(Y) \subseteq \FO(h)$. Then $Y$ has to act trivially on the $2$-set
$\FO(Y)\setminus{\FO(Y)}$, which contradicts $|\FO(Y)| = 2$.
Thus the second case holds true, and by Lemma~\ref{Inv-lem:TB}\,(a) we may suppose that \(U\leq H\). We conclude $P \leq N_G(U) \leq N_G(H)$ with Lemma~\ref{Inv-lem:TB}\,(b).

If $|\FO(X)| = 4$, then $P$ is contained in the element-wise stabilizer of
$\FO(X)$.

 %So there is an element of order \(3\) that stabilizes \(\FO(U)=\Delta\) and that, therefore fixes
%one or four points. But this contradicting the fact that every nontrivial \(3\)-element of \(G\) has exactly two fixed points.
 %As $\aut(P)$ is a $\{2,3\}$-group,
 %So we get in both cases that $P$ normalizes a nontrivial $p$-subgroup. This yields $P \leq N_G(H)$.
% Therefore the assumption was false, and we deduce that \(|\FO(X)|=4\). %Let \(H\) the point-wise stabilizer of \(\FO(X)\).
 %Then \(H\) is nontrivial and
% As \(P\le N_G(X)\), %it acts on \(\FO(X)\). As \(P\) has order \(9\) this implies that \(P\) lies in a point stabilizer.
 %Also, by
% we may apply Lemma~\ref{Inv-lem:TB}\,(b), and this gives that \(P\leq N_G(H)\).
 \end{proof}

 \begin{lemma}\label{SL2xSL2sub}
		Let \(q %\ge 3
		\) be a prime power such that \(q \equiv 2,4,5\) or
	\(7\) modulo \(9\) and let \(G\) be a finite simple group such that $G$  has a subgroup $U = V/Z$, where $V \cong \SL_2(q) \times \SL_2(q)$ and $Z \leq Z(V)$,  and such that the Sylow $3$-subgroups of $G$ are elementary abelian of order $9$.
	Then there does not exist a set \(\Omega\) such that \((G, \Omega)\) satisfies Hypothesis \ref{odd3simple}.
\end{lemma}
\begin{proof}
Assume that \(\Omega\) is such that \((G, \Omega)\) satisfies Hypothesis \ref{odd3simple}.
Since $U$ does contain a Sylow $3$-subgroup of $G$, we
see that $N_G(X)$ has even order for every subgroup $X$ of $G$ of order $3$.
More precisely, for every subgroup of \(G\) of order \(3\) there is an element in \(G\) that inverts this subgroup. %Gilt bereits für U und dann Sylow
Using Lemma~\ref{helpSL2xSL2sub} we find a nontrivial four point stabilizer \(H\) such that $N_G(H)$ contains a Sylow $3$-subgroup $P$ of $G$.

Assume that $P \not\leq H$. Then  $N_G(H)$ induces a subgroup of $\Sym(4)$ on $\Delta:= \FO(H)$, which is divisible by $3$. As $|G_\alpha|$ is odd, we get that
$\Alt(4)$ is induced.
Further, we get that $|P\cap H| = 3$.
%Lemma \ref{Inv-lem:TB}(a) yields that $H$ contains a Sylow \(p\)-subgroup for all $p \ge 5$ in $\pi(H)$.
Due to Frattini's lemma $N_G(P \cap H) \leq N_G(H)$ contains a Sylow-$2$ subgroup of
$N_G(H)$. Therefore
\(N_G(P\cap H)\) has a subgroup isomorphic to $C_3 \times \Alt_4$, which is a $\{2,3\}$-Hall subgroup of $N_G(P \cap H)$.
However, this means that \(N_G(P\cap H)\) does not contain an element inverting \(P\cap H\), which gives a contradiction.
%, which contradicts the structure of $N_U(\tilde{P})$
%for $\tilde{P}$ a Sylow $3$-subgroup of $U$, as in $U$ every $3$-element is inverted by some $2$-element.
%Assume \(P\not\leq H\). Then \(P\cap H\) has order \(3\).
%%Let \(y\in P\setminus H\). Then \(y\in N_G(P\cap H)\), so \(y\) acts on \(\Delta\). Hence \(|\FD(y)|=1\), because \(y\notin H\).
%As noted earlier there exists an element \(s\in N_G(P\cap H)\) of even order. Then \(s\) and \(P\) act on \(\Delta\) and as point stabilizers of \(G\) have odd order
%this implies that \(N_G(H)\) induces the group \(\Alt_4\) on \(\Delta\). %because as \(s\) acts on (\Delta\) by Lemma~\ref{4ptstab} it follow that \(s\in N_G(H)\)

This shows that $P \leq H$ and that Case (e) of Lemma~\ref{syl3} holds. So Proposition~\ref{3semireg} yields that \(G\) contains a strongly \(3\)-embedded subgroup,
contradicting  Hypothesis \ref{odd3simple}.
%Also Fall (e) von Syl3 und damit hat G stark 3-eingebettete UG
%and that $N_G(X) \leq N_G(H)$ for every nontrivial $3$-subgroup of $P$ by Lemma~\ref{4ptstab}. Let $X \leq U$ such that $C_U(X)$ contains a subgroup isomorphic to $\SL_2(q)$ or $\PSL_2(q)$. Then $C_G(X) \leq N_G(X) \leq N_G(H)$ yields that $C_U(X)$ is soluble. Thus $q = 2$ or $q= 3$, which contradicts our assumption on $q$.
\end{proof}

\begin{lemma}\label{odd3el9}
	Suppose that Hypothesis \ref{odd3simple} holds and that \(G\) %is a finite simple group with
	has elementary abelian
	Sylow \(3\)-subgroups of order \(9\).
	Then \(G \cong \Alt_6\) and the point stabilizers are elementary abelian of order~\(9\).
\end{lemma}

\begin{proof}
We use Lemma~\ref{CFSGelab} to determine the possibilities for $G$.
	Table~\ref{SmallGroupsTable} shows that the candidates \(\Alt_7\), \(M_{11}\), \(M_{22}\), \(M_{23}\) and \(\textrm{HS}\) do not exhibit examples for Hypothesis \ref{odd3}.

\(\Alt_6\) is also included in Table \ref{SmallGroupsTable} with exactly the action that we describe in the lemma.
So it remains to prove that the remaining finite simple groups with elementary abelian Sylow $3$-subgroups of order $9$ do not give rise to examples for Hypothesis \ref{odd3simple}.
	
	Let $\varepsilon \in \{1,-1\}$. If \(G =\PSL^{\varepsilon}_3(q)\), then Hypothesis \ref{odd3} is not
	satisfied.
	
	So let \(G\) be one of the remaining groups: \(\PSL^{\varepsilon}_4(q)\), where \(q\equiv
	3-\varepsilon\) or \(6-\varepsilon\) modulo \(9\), \(\PSL^{\varepsilon}_5(q)\), where \(q\equiv 3-\varepsilon\) or \(6-\varepsilon\) modulo \(9\)
	or \(\PSp_4(q)\), where \(q\equiv 2,4,5\) or \(7\) modulo \(9\).
	%Then in all cases \(q\) is a prime power congruent to \(2,4,5\) or \(7\) modulo~\(9\).
	%Let \(H\) be the corresponding matrix group: \(\SL^{\varepsilon}_4(q)\) for
	%\(\PSL^{\varepsilon}_4(q)\), \(\SL^{\varepsilon}_5(q)\) for \(\PSL^{\varepsilon}_5(q)\) and \(\Spl_4(q)\) for
	%\(\PSp_5(q)\).
	%TODO groups and have consistent machen PSL(2,q) ist eine Familie von Gruppen
	Then we consider Theorems~3.5\,(i) and 3.9\,(i) in~\cite{Wilson} and we see that \(\SL_4(q)\), \(\SU_4(q)\),
	\(\SL_5(q)\) and \(\SU_5(q)\) all
	have a subgroup \(V\) isomorphic to \(\SL_2(q) \times \SL_2(q) \cong  \SU_2(q)\times \SU_2(q)\).
	The group \(\Spl_4(q)\) has a subgroup \(V \cong \Spl_2(q) \times \Spl_2(q) \cong \SL_2(q) \times
	\SL_2(q)\).
	In all these cases we see that Lemma \ref{SL2xSL2sub} is applicable, and
	it shows that %, if \(q\geq 3\), then
	none of these groups satisfy Hypothesis \ref{odd3simple}.
\end{proof}

%%%%%%%%%%%%

\section{The case where point stabilizers have order coprime to 6}

\begin{hyp}\label{hyp6}
	In addition to Hypothesis~\ref{4fix}, we let $\alpha \in \Omega$ and we suppose that $2$ and $3$ do not
	divide $|G_\alpha|$. We fix four distinct points $\alpha, \beta, \gamma, \delta \in \Omega$, such that the element-wise stabilizer $H$ of $\Delta := \{\alpha, \beta, \gamma, \delta\}$ is a nontrivial subgroup of $G$.
	Finally, we let $p \in \pi(H)$ and $P \in \syl_p(H)$.
\end{hyp}

We note that Hypothesis \ref{hyp6} implies that $|\Omega| \ge 6$, because in $\Sym_5$ the point stabilizers are $\{2,3\}$-groups.

We begin by collecting the information that we have already generated in earlier sections.
For this we recall that, for all $n \in \N$, $G$ acts as a $(0,n)$-group on $\Omega$ if and only if all elements of $G^\#$ have zero or $n$ fixed points on $\Omega$ (following \cite{PS1}).

\begin{lemma} \label{start}
	Suppose that Hypothesis~\ref{hyp6} holds. Then the following are true:
	\begin{enumerate}
		\item[(a)] $|N_G(H):H| = 4$, $N_G(H)$ is exactly the stabilizer of the set $\Delta$ in $G$ and $N_G(H)$ acts transitively on $\Delta$,
		\item[(b)] $H$ is a TI-subgroup,
		\item[(c)] $G$ acts as a $(0,4)$-group on the set of right cosets of $H$,
		\item[(d)] either $G_\alpha = H$ or  $G_\alpha$ is a Frobenius group with complement $H$.
		\item[(e)] If  \(1 < X\le H\), then  \(|N_{G}(X):N_H(X)|\leq 4\).\\
		In particular, \(P \in \syl_p(G)\)  and therefore
		$(|H|,|\Omega|) = 1$.
		\item[(f)] If $1 < R \le H$, then $|N_G(R)|$ is not divisible by $3$ or $8$.
		\item[(g)] $N_G(P) \leq N_G(H)$.	
	\end{enumerate}
\end{lemma}

\begin{proof}
	As $H$ has odd order, Lemma \ref{Inv-lem:TB} yields that $|N_G(H):H| = 4$. Moreover, $H \neq 1$, so there is a prime $p \ge 5$ in $\pi(H)$ and $H$ contains a Sylow $p$-subgroup of $G$ by Lemma \ref{Inv-lem:TB}(a). Now Part (e) of the same lemma yields that $N_G(H)$ acts transitively on $\Delta$. Of course $N_G(H)$ stabilizes $\Delta$, and now, conversely, suppose that $g \in G$ stabilizes $\Delta$. If $g \in H$, then there is nothing left to prove. Otherwise suppose that $\alpha^g \neq \alpha$ and, using the transitive action of $N_G(H)$ on $\Delta$, let $x\in N_G(H)$ be such that $\alpha^g=\alpha^x$. Then $gx^{-1} \in G_\alpha$, which means that $gx^{-1}$ fixes every point of $\Delta$ (by Hypothesis \ref{hyp6}) and then $g \in N_G(H)$.
 Together this is (a).
	Part (b) is exactly Lemma \ref{Inv-lem:TB}(f).
	Lemma 1.3 combined with Lemma 1.1 of \cite{PS1} proves (c).
	For (d) suppose that \(G_\alpha \neq H\). Assume that there exists \(g\in G_{\alpha}\setminus H\) such that \(H\cap H^g\neq 1\).
	Then (b) yields that \(g\in N_G(H)\), and therefore \(g\) stabilizes the set \(\Delta\setminus \{\alpha\}\), contrary to (a). This proves (d).
	Now we turn to (e)
	and we let $1 \neq X \le H$.
	First we notice that $X$ fixes every point in $\Delta$, but not any more points because of our global fixity 4 hypothesis.
	Thus $N_G(X)$ stabilizes the set $\FO(X)=\Delta$. Now we recall that
	\(|N_{G}(X):N_{G_{\alpha}}(X)|\leq 4\) by Lemma \ref{Inv-lem:TB}(a) and we let
	\(y\in N_{G_{\alpha}}(X)\).
	
	Assume that \(y\) does not fix all
	four points in \(\Delta\). Then it acts as a transposition or as a $3$-cycle on $\{\beta, \gamma, \delta\}$.
	But $o(y)$ has order coprime to $6$, so this is not possible.
	Thus \(y\in H\). We conclude that
	\(|N_{G}(X):N_H(X)|=|N_G(X):N_{G_{\alpha}}(X)|\cdot |N_{G_{\alpha}}(X):N_{H}(X)|\leq 4\), which is the first assertion.
	The second statement follows because, if $p \in \pi(H)$, then $p \ge 5$.
	For (f) we use (e) and the fact that $|H|$ has order coprime to $6$ by hypothesis.
	Assertion (g) holds by (a) because $N_G(P)$ stabilizes $\Delta$.
\end{proof}

%Our strategy now is to first identify the $(0,4)$-actions. This means that first we suppose that $H = G_\alpha$, and then we deduce from these results the actions where $G_\alpha$ is a Frobenius group.
We remark that
the simple group $M_{11}$ is a $(0,4)$-group in its action on the set of cosets of a Sylow $5$-subgroup.
This action is imprimitive, with blocks of size $11$ and with a block stabilizer $F_{55}$ (a Frobenius
group of order $55$). The action of $M_{11}$ on the set of cosets of $F_{55}$ is an example of an
action satisfying Hypothesis \ref{4fix}, and here $G_\alpha$ is a Frobenius group of order coprime to $6$.

\begin{thm}\label{strongcases}
	Suppose that Hypothesis \ref{hyp6} holds. If $G$ is simple, then one of the following holds:
	
	(i) $G$ has cyclic Sylow $p$-subgroups, or
	
	(ii) there is some $n \in \N$ such that $G$ is isomorphic to $\PSL_2(p^n)$ or $\PSU_3(p^n)$
	and $N_G(H)$ is a Borel subgroup.
\end{thm}

\begin{proof}
	$N_G(H)$ is strongly
	$p$-embedded in $G$ by Lemma \ref{Inv-lem:TB}(h), and then
	Proposition 17.11 in \cite{GLS2} yields that
	$N_G(H)$ is also strongly $p$-embedded as defined in \cite{GLS2}. Therefore
	Theorems 7.6.1 and 7.6.2 in \cite{GLS3} apply.
	If $G$ has cyclic Sylow $p$-subgroups, then (i) holds, so it is left to go through the cases of
	Theorem 7.6.1.
	Case (a) of the theorem leads to our Case (ii).
	
	Suppose that Case (b) of the theorem holds, which means that
	$p \ge 5$, $G \cong \A_{2p}$ and $N_G(H)$ is the stabilizer of a partition of $\Omega$ into two subsets of size $p$.
	Then we see that $N_G(H)$ contains elements of order $3$, contrary to our hypothesis and Lemma \ref{start}~(a).
	
	The cases (c) and (d) do not occur because $p \ge 5$ by Hypothesis \ref{hyp6}.
	In Case (e) we have that  $p=5$, $G \cong \,^2 F_4(2)'$ and $N_G(H)$ is the normalizer of a Sylow $5$-subgroup, containing a subgroup of structure $C_4 \ast \SL_2(3)$. This is impossible by Lemma \ref{start}~(a).
	
	In Case (f) we have $p=5$ again, but this time $G \cong$ Mc
	and $N_G(H)$ is the normalizer of a Sylow $5$-subgroup, containing a subgroup of structure $C_3 : C_8$.
	In Case (g) we have that $p=5$ again, and  $G \cong$ Fi$_{22}$, $N_G(H)\cong \aut(D_4(2))$.
	The last case is $p=11$, $G \cong $J$_4$ and $N_G(H)$ is the normalizer of a Sylow $11$-subgroup, containing a subgroup of structure $C_5 \times \GL_2(3)$.
	
	These three cases have in common that they do not occur because of Lemma \ref{start}~(a).
	So this concludes the proof.
\end{proof}

The remainder of this section is devoted to the analysis of all finite simple groups
from the perspective of a possible fixity 4 action satisfying Hypothesis \ref{hyp6}. Further we assume that every Sylow
$p$-subgroup of $H$ is cyclic.
We begin with the alternating groups and the sporadic groups, and then the main work will occur for the groups of Lie type.

\subsection{Alternating and sporadic groups}

\begin{lemma}\label{Altcyclic}
	Suppose that Hypothesis \ref{hyp6} holds, that $n \in \N$, $n \ge 5$ and $G=\A_n$. Let $p \in \pi(H)$ and suppose that $P\in \syl_p(G)$ is cyclic.
	Then $G=\A_{7}$ and $H=G_\alpha$ has order $5$.
\end{lemma}

\begin{proof}
	For the $p$-rank of $G$, we refer to Proposition 5.2.10 in\cite{GLS3}.
	%and we note that $G$ and its central extensions have the same %P: Brauchen wir die central extensions überhaupt?
	%$p$-rank, because $p\ge 5$.
	Thus $1=r(P)=\left[\frac{n}{p}\right]$, which forces $n <2p$.
	Let $h \in H$ be a $p$-cycle.
	If $n-p \ge 3$, then $C_G(h)$ contains a $3$-cycle, which contradicts our hypothesis because
	$C_G(h)$ stabilizes $\FO(h) = \Delta$.
	Hence $n-p \in \{0, 1,2\}$.
	
	So we have that $n \in \{ p, p+1, p+2\}$, and in each case  $N_G(\langle h \rangle)/\langle h \rangle$ is
	cyclic of order $(p-1)$ or $(p-1)/2$. If $a$ divides $(p-1)/2$,
	then $a < n/2$ in all three cases. This shows that $(p-1)/2$ is a
	$2$-power, and an application of Lemma~\ref{start} (f) yields
	that
	$(p-1)/2 \in \{1,2,4\}$. The fact that $p \geq 5$ gives that $p = 5$ and $n \in \{ 5,6,7\}$. Then Lemma~\ref{start} (a)  forces $n = 7$.
\end{proof}

\begin{lemma}\label{spor}
	Suppose that $G$ is a sporadic group satisfying Hypothesis \ref{hyp6}.
	Then $G=\M_{11}$ or $G=\M_{22}$ and in both cases, two examples are possible. The point stabilizers are cyclic of order $5$ or Frobenius groups of order $55$, respectively.
\end{lemma}

\begin{proof}
	%It can be verified with GAP code (indlude code?) that
	By Remark~\ref{rem4.1} the examples that we state are indeed examples where Hypothesis \ref{hyp6} is satisfied. %TODO Die Richtung brauchen wir doch eigentlich nicht.
	So now we prove that these are the only ones, using the notation from Hypothesis \ref{hyp6}. Theorem \ref{strongcases} yields that
	$G$ has cyclic Sylow $p$-subgroups for all $p \in \pi(H)$.
	So we let $p \in \pi(H)$ and $P\in \syl_p(G)$, and we let $X \le P$ be a subgroup of order $p$. We know that $P$ is cyclic and that $|N_G(X)|$ is not divisible by $6$,  $9$ or $8$, using Lemma \ref{Inv-lem:TB}~(c).
	This restricts the possibilities for $p$, as we can see in Tables 5.3a--z in \cite{GLS3}.
	
	\smallskip
	The Mathieu groups:
	
	If $p=5$ or $p=11$, then this leads to the possibilities stated for $\M_{11}$ and $\M_{22}$.
	Otherwise the tables show that the restrictions on $|N_G(X)|$ are not compatible with the
	subgroup structure of $G$.
	
	%P: Der Anfang ist recht kurz und warum kommen jetzt noch einmal alle sporadischen Gruppen?
	\smallskip
	The Janko groups:
	
	Here the tables in \cite{GLS3} show that, again, the restrictions on $|N_G(X)|$ are not satisfied.
	
	\smallskip
	The remaining groups:
	Inspection of the corresponding tables
	shows that, with several applications of Lemma \ref{Inv-lem:TB} (a), eventually $2 \in \pi(H)$ in each case. This is impossible.
\end{proof}

\subsection{Groups of Lie type}

\begin{lemma}\label{Lie1}
	Suppose that $G$ is a simple group of Lie type satisfying Hypothesis \ref{hyp6}, in characteristic~$p$.
	If $G$ has cyclic Sylow $p$-subgroups, then
	$G=\PSL_2(p)$ and $|G_\alpha|=p \cdot \frac{p-1}{4}$.
\end{lemma}

\begin{proof}
	We inspect Theorem 3.3.3 in \cite{GLS3}, and hence Table 3.3.1 there, for the equicharacteristic rank.
    Since \(p\geq 5\), the only simple group there of \(p\)-rank \(1\) is $\PSL_2(p) \cong \PSU_2(p)$.
    Therefore our assertion follows from Theorem \ref{L2q-main} and the fact that $|G_\alpha|$ is coprime to $6$.
%	First we note that $p \ge 5$ and therefore we do not have to consider the exceptions $B_2(2)$, $G_2(2)$, $^2G_2(3)$
%	and $^2F_4(2)$.
%	Let $a \in \N$ be such that $G$ is defined over a field of order $q:=p^a$.
%	If $a \ge 2$, then Table 3.3.1 shows for all types of groups that $r_p(G) \ge 2$, contrary to our hypothesis.
%	If $a=1$, then the table shows that all exceptional types do not occur, where for the Suzuki groups we use that $p\neq 2$.
%	For the classical types, we let $n \in \N$ and we begin with $\PSL_n(p)$, where $n=2$ is the only possibility for cyclic Sylow $p$-subgroups.
%	With $\PSU_n(p)$, we also have $p$-rank 1 only if $n=2$.
%	The odd-dimensional orthogonal groups do not give examples because the contribution of the dimension to the $p$-rank is already too large.
%	For the same reason, the symplectic groups do not occur.
%	Finally, if $G$ is an even-dimensional orthogonal group, then the fact that $p$ is odd gives that $r_p(G)>a=1$, so this is impossible.
%	
%	As $\PSL_2(p) \cong \PSU_2(p)$, our assertion follows from Theorem \ref{L2q-main} and the fact that $|G_\alpha|$ is coprime to $6$.
\end{proof}

Before we enter a detailed analysis of the series of finite simple groups of Lie type,
we prove two preparatory lemmas.

\begin{lemma}\label{help}
	
	(a) Let \(G\) be a finite group with a normal subgroup \(N\) of index \(m \in \N\). Let \(p\in
	\pi(N)\) and \(P\in \Syl_p(N)\). Then \(|N_G(P)|=m\cdot |N_N(P)|\).

	Now suppose that Hypothesis \ref{hyp6} holds and that $G$ is a simple group of classical Lie type over a field with $q$ elements,
	where $q$ is a prime power coprime to $p$. Moreover, suppose that $G$ has cyclic Sylow $p$-subgroups and that $G \notin \{\PSL_2(q),\Sz(q), \PSL_3(q),\PSU_3(q)\}$.
	
	(b) $p$ does not divide $q^2-1$.
	%R: Kommt mir jetzt mit Lemma 8.9 üflü vor, sollten wir einmal prüfen.
	
	(c) $p$ divides $q^{p-1}-1$.
	
\end{lemma}		

\begin{proof}	
	For (a) we use a Frattini argument.
	Then \(G=N\cdot N_{G}(P)\) and \(m\cdot |N|=|G|=|N\cdot N_{G}(P)|=\frac{|N|\cdot |N_{G}(P)|}{|N_G(P)\cap
		N|}=\frac{|N|\cdot |N_{G}(P)|}{|N_N(P)|}\), which implies the statement.
	
	(b) 	If \(p\) divides \(q-1=\Phi_1(q)\) or \(q+1=\Phi_2(q)\), then we use
	(10-1)~in~\cite{GL} and we deduce that $r_p(G) \ge 2$, contrary to
	the hypothesis that $G$ has cyclic Sylow $p$-subgroups.
	
	(c) This is Fermat's Little Theorem.
\end{proof}

In our analysis we will make heavy use of the concept of a \textbf{primitive prime divisor}.
So we recall that, if $q$ is a prime power and $e \in \N$, then a prime divisor $r$ of $q^e-1$ is said to be a primitive prime divisor of $q^e-1$  if and only if $r$ does not
divide any of the numbers $q^i-1$, where $1 \le i < e$.
It was proved by Zsigmondy (see \cite{Zsig}) that, if $e>3$ and  $(q,e) \neq (2,6)$, then there
exists a  primitive prime divisor for
$q^e-1$.
In the following lemma, we use the notation
$d_q(p)$ for the positive integer $e$ such that $p$ is a primitive prime divisor of $q^e-1$, and
$d_{q^2}(p)$ correspondingly.

\begin{lemma}\label{Primesalternativ}
	Let $G$ be a finite simple group of Lie type  that satisfies Hypothesis~\ref{hyp6} and is defined over a field with $q$ elements.
	Further let $p \in \pi(H)$ be such that the Sylow $p$-subgroups of $G$
	are cyclic. Then
% the smallest exponent $k$ for which a particular polynomial (often a divisor of $|G|$), as given in the table, could be divisible by $p$ is given in the following table. We also give
the values for $d_q(p)$ and $d_{q^2}(p)$ are as given in the table.  If we did not calculate one
	of the values, then we put "--" in the respective column.\\
	
	%\begin{table}
	\begin{center}
		\begin{tabular}{|c|c|c|c|}
\hline
%\multicolumn{2}{|c|}{|c|} \\
%\hline
%\hline
$G$ & remarks  & $d_q(p)$ & $d_{q^2}(p)$ \\
			\hline
			$\PSL_n(q)$& $n \geq 4$ &  $n-1,n$ & -- \\
			\hline
			$\PSU_n(q)$& $n \geq 4$, $n$ even  &  $\neq n-1$& $n-1,n/2$ \\
			                       &    \hspace{1cm} $n$ odd &  $ \neq n $  & $(n-1)/2, n$\\
			\hline
			$\PSp_{2n}(q)$ & $n \ge 2$ & --& $n$ \\
			\hline
			$P\Omega_{2n+1}(q)$ & $n \ge 3$ &  --& $n$\\
			\hline
			$P\Omega^+_{2n}(q)$ & $n \ge 4, n$ odd &  $n$ & -- \\
                         & $n$ even, $q \not\equiv 1 \mod 4$ &  $n-1$ & $n-1$ \\
			\hline
			$P\Omega^-_{2n}(q)$ & $n \ge 4$ &  $\neq n $ & $n$\\
                        & & $\neq n-1 $ & $n-1$\\
			\hline
			$^3D_4(q)$ &  & $12$& --\\
			\hline
			$F_4(q)$ & & $8,12$& --\\
   \hline
			$^2F_4(q)$ & &  $6,12$& -- \\
			\hline
			$G_2(q)$ & & $3,6$ & --\\
			\hline
			$^2G_2(q)$ & &  $1,2,6$ & -- \\
			\hline
			$E_6(q)$ & &  $5,8,9,12$  & --\\
			\hline
			$^2E_6(q)$ &  &  $8,10,12,18$  & --\\
			\hline
			$E_7(q)$ & &  $5,7,8,9,10,12,14,18$ & --\\
			\hline
			$E_8(q)$ & &  $7,9,14,15,18,20,24,30$ & --\\
			\hline
			
		\end{tabular}\label{ktablealternativ}
		
	\end{center}
	%\caption{Potential prime divisors of $|G_\alpha|$.}
	%\end{table}
\end{lemma}
\smallskip

%\begin{lemma}\label{Primes} Let $G$ be a finite simple group of Lie type  satisfying Hypothesis~\ref{hyp6}, defined over a field with $q$ elements. Further let $p \in \pi(H)$ be such that the Sylow $p$-subgroups of $G$ are cyclic. Then the smallest $k$ for which $\Phi_k(q)$ could be divisible by $p$ is given in the following table: \end{lemma}

\begin{proof}
	Let $p \in \pi(H)$.
First let $G$ be a classical group of Lie type defined over a vector space $(V, f)$ with bilinear form $f$. If $G = \PSL_n(q)$, then
$f$ is the form that is constant zero. Further let $U$ be a $2$-dimensional subspace of $V$ and $W$ a complement to $U$ in $V$. If $G$ is of unitary or
symplectic type, then let $(U, f_{\mid U})$ be a non-degenerate $2$-dimensional subspace of unitary or symplectic type and let $W$ denote its orthogonal
complement in $V$. Moreover, in all three cases let $K$ be the stabilizer of $U$ in $G$.
Then $K$ contains subgroups $K_U$ and $K_W$ such that, for $\{X,Y\} = \{U,W\}$, the subgroup  $K_X$ acts on the projective space $(P(X), f_{\mid X}$) as $\PSL(X)$, $\PSU(X)$ or  $\PSp(X)$, respectively, and it acts trivially on $Y$ (see \cite[Propositions~4.1.3 and 4.1.4]{KL}).
 If $p \in \pi(K_W)$, then without loss a Sylow $p$-subgroup of $K_W$
 is contained in $H$. Therefore, $K_U \leq N_G(H)$, which implies that $|K_U|$ is not divisible by $3$ and that $|K_U|_2 \leq 4$. If $G$ is linear, unitary or symplectic, then $K_U$ induces $\PSL_2(q)$ on $P(U)$. But the order of this group  is divisible by $3$, contrary to Lemma \ref{start}\,(f). It follows that $p \not\in \pi(K_W)$, which yields the assertion for those groups.

 Now let \(G\) be an orthogonal group. Set \(m\in \{2n,2n+1\}\) and \(\varepsilon\in \{-1,0,1\}\) such that
    \(G=P\Omega_m^{\varepsilon}(q)\), where \(\varepsilon=0\) if \(m\) is odd and \(\varepsilon \in \{-1,1\}\) if \(m\) is
    even. %TODO Muss man die Notation noch weiter erklären?
    Let \(k\) be a positive integer such that \(p\) divides \(q^{2k}-1\). Then there exists \(\varepsilon_2\in \{-1,1\}\)
    such that \(p\) divides \(q^k-\varepsilon_2\).
    Let \(m_1=3\) if \(m\) is odd and let \(m_1=4\) if \(m\) is even. Set \(m_2=m-m_1\) and \(\varepsilon_1=\varepsilon
    \cdot \varepsilon_{2}\).

    Furthermore let $(U, f_{\mid U})$ be a non-degenerate $m_1$-space of $\varepsilon_1$-type. As above, let
    $W$ be the orthogonal complement to $U$ in $V$ and $K$ the stabilizer of $U$ in $G$. Then $K$ contains a subgroup $K_W$
    which acts as  $\Omega_{m_2}^{\varepsilon_2}(q)$ on $W$ and a subgroup \(K_U\) which acts  as \(\Omega_{m_1}^{\varepsilon_1}(q)\) on
    \(U\) (see \cite[Propositions~4.1.6]{KL}).
    If \(p\in \pi(K_W)\), then
    without loss \(K_W\) contains a Sylow \(p\)-subgroup \(P\) of \(H\). Thus \(K_U\leq N_G(P)\).
    Since the order of \(P\Omega_{m_1}^{\varepsilon_1}(q)\) is divisible by \(q(q^2-1)\) and hence by \(3\),
    \(|N_G(P)|\) is divisible by
    \(3\) contrary to Lemma~\ref{start}\,(f). Therefore \(p\) does not divide \(|K_W|\). We conclude that \(k\geq
    n-1\) and, if $n$ is odd, that $k=n$.

Notice also that if $n$ is even and $\varepsilon = 1$, then $G$ contains a subgroup, which is modulo its center isomorphic to $P\Omega_n^+(q) \times P\Omega_n^+(q)$.
Therefore, if $p$ divides the order of the latter group, then the $p$-Sylow subgroups are not cyclic.
 %Then $|P\Omega_n^-(q)|$ is  divisible by $p$.
 It follows that, if $d_q(p) = n$ and $\varepsilon = 1$, then
 $n$ is odd.
 Suppose that $p$ divides $q^{2(n-1)} -1$.
By \cite[Propositions~4.1.6]{KL}
$G$ contains a subgroup isomorphic to $2.(P\Omega_{2}^-(q) \times P\Omega_{2n-2}^-(q))$ or to $\Omega_{2n-2}^-(q) \times \Omega_{2n-2}^-(q)$. We derive from Lemma~\ref{help}(b) that
$p$ divides $q^{n-1} -1$. Further notice that by the same reference $G$ also contains a subgroup isomorphic to $P\Omega_{2(n-1)}^+(q)$ or $P\Omega_{2(n-1)}^+(q)$. By our arguments above, it follows that  $n-1$ is odd and $n$ even.
Furthermore, we derive from  \cite[Propositions~4.1.6]{KL} and Lemma~\ref{start} (f) that either $q$ is even or $q$ is odd and $q-1$ is not divisible by $4$.  The remaining arguments stay the same.
	
Now let $G$ be a group of Lie type of exceptional type. The group $G$ has cyclic Sylow $p$-subgroups by hypothesis and $p \ge 5$. Then we can use Theorem (10-1) and Table 10:2 in \cite{GL},
	where $r_p(G)$ is the exponent of $\Phi_k(q)$.
	This gives exactly the possibilities in our table.
\end{proof}

Next we study the classical groups in detail. There we will repeatedly use the fact that $N_G(P) \leq N_G(H)$ by Lemma~\ref{start}(g).

\newpage
\subsection{Classical groups}

\begin{lemma}\label{tableTX}
Suppose that $G$ is a classical group of Lie type and that $(G,\Omega)$  satisfies Hypothesis~\ref{hyp6}.
Then there is a divisor $D$ of $|N_G(P)|$ as given in
the following table:

	\begin{center}
		\begin{tabular}{|c|c|c|c|c|c|c|}
			\hline
			$G$ & conditions & $d_q(p)$ & $d_{q^2}(p)$ & divisor $D$ & citation & remarks\\
			\hline
			$\PSL_n(q)$& $n \geq 4$ & $n$ & -- & $\frac{(q^n-1) \cdot n}{(q-1) \cdot (q-1,n)}$
			& \cite[7.3]{Huppert-I}& $D= |N_G(P)|$
			\\
			& $ (n,q) \neq (6,2)$ &  & -- &
			& &
			\\
			\hline
			&  & $n-1$ & -- & $\frac{(q^{n-1}-1) \cdot (n-1)}{(q-1,n)}$
			& \cite[7.3]{Huppert-I} &
			\\
		 		 &
			&  &  &
			 &  & \\
			\hline
			$\PSU_n(q)$& $n \geq 4$, $n$ even
			& -- & $\frac{n}{2}$ & $\frac{(q^n-1)n}{(q+1)(q+1,n) }$
			 & \cite[4.2.4]{KL}, \cite[7.3]{Huppert-I} & $\frac{n}{2}$ divides $D$\\
	%		 & $p$ divides $q^{\frac{n}{2}}-\varepsilon$
	%		& -- &  &
	%		 &  & $\varepsilon \in \{1,-1\}$ \\
			 &
			&  &  &
			 &  & \\
			 \hline
			 &
			& $\neq n-1$ & $n-1$ & $\frac{(q^{n-1}+1)(n-1)}{(q+1,n)}$
			 & \cite[4.1.4,~4.3.6]{KL} & $n-1$ divides $D$\\
			 &
			&  &  &
			 &  & \\
			 \hline
			 & $n \geq 4$, $n$ odd
			& $\neq n$ & $n$ & $\frac{(q^{n}+1)n}{(q+1)(q+1,n)}$
			 & \cite[4.3.6]{KL}  & \\
			 		 &
			&  &  &
			 &  & \\
			 \hline
			  &
			& $\neq \frac{n-1}{2}$ & $\frac{n-1}{2}$ & $\frac{(q^{n-1}-1)(n-1)}{(q+1,n)}$
			 & \cite[4.1.4,~4.2.4]{KL}, \cite[7.3]{Huppert-I}  & \\
			   &
			&  &  &
			 &  & \\
			 \hline
		$\PSp_{2n}(q)$ & $p$ divides $q^n-\varepsilon$ & --
& n &		$\frac{(q^n-\varepsilon)2n}{(2,q-1)}$& \cite[4.3.10]{KL}&  $\varepsilon \in \{1,-1\}$ \\
			 		 &
			&  &  &
			 &  & \\
			\hline
					$P\Omega_{2n+1}(q)$ & $n$ even & $\neq n$ & $n$ & $\frac{(q^{n}+1)2n}{(q^n+1,4)}$ & \cite[4.1.6,~4.3.15]{KL}&\\
			 &  &  &  &  & &\\
			 	\hline
			 & $n$ odd & -- & $n$ & $\frac{(q^{n}-\varepsilon)n}{(q^n-\varepsilon,4)}$ & \cite[4.1.6]{KL}&$\varepsilon \in \{1,-1\}$\\
			& $p$ divides $q^n-\varepsilon$ &  &  &  &  & \\
			 \hline
			$P\Omega^+_{2n}(q)$ & $n$ odd & $n$ & -- & $\frac{(q^{n}-1)n}{(q-1,4)}$  & \cite[4.1.20, ~4.2.7]{KL}& $n$ divides $D$ \\
   			 &  &  &  &  &   & \\
			\hline
			 & $n$ even & $n-1$ & $n-1$ & $\frac{(q^{n-1}-1)(n-1)(q-1)}{(q-1,2)}$  & \cite[4.1.6,~4.1.20,~4.3.15]{KL} & $q \not\equiv 1$ mod $4$\\
			 &  &  &  &  &  & \\
			\hline
			 $P\Omega^-_{2n}(q)$ & & $\neq n$ & $n$ & $\frac{(q^{n}+1)n}{(q^n+1,4)}$  & \cite[4.3.15]{KL}& \\
			 		 &
			&  &  &
			 &  & \\
			\hline
			 & & $\neq n-1$ & $n-1$ & $\frac{(q^{n-1}+1)(n-1)}{(q+1,2)}$  & \cite[4.1.6,~4.3.15]{KL}& \\
						 		 &
			&  &  &
			 &  & \\
			\hline

		\end{tabular}\label{TX}
		
	\end{center}

\end{lemma}
\begin{proof}
In the proof we refer to \cite{KL} and we often use the notation there, but we should point out that there is a piece of non-matching notation: The letter $n$ in \cite{KL}
does not have the same meaning as our $n$ in the table if $G$ is a symplectic or an orthogonal group.
Also, the notation that is used in the specific results that we refer to is often explained earlier, e.g. on pages 57, 60 and in 80 onward (in \cite{KL}).

Now we consider the cases in the table one by one,
and we set $N:= |N_G(P)|$.

\smallskip
\underline{$G = \PSL_n(q)$:}

If $d_q(p) = n$, then the claim follows with \cite[Satz 7.2]{Huppert-I}.
Let $d_q(p) = n-1$.  Then according to \cite[Lemma~4.1.4]{KL}, $G$ contains a subgroup $Z\PSL_{n-1}(q)B$, where $Z$ commutes with $\PSL_{n-1}(q)$ and has order
$(q-1)/(q-1,n)$ and $B$ has order $(q-1,n-1)$.
With Lemma~\ref{help}(a) and \cite[Satz 7.2]{Huppert-I} we obtain that $D$, as given in the table, divides $N$.
\\
\underline{$G = \PSU_n(q)$, $n \geq 4$:}

Let $n$ be even.
We start with the case where $d_{q^2}(p) = n/2$.
 According to \cite[4.2.4]{KL}, $G$ has a subgroup
of structure $A.\PSL_{n/2}(q^2).B.2$, where $A$ has order  $(q-1)\cdot (q+1,\frac{n}{2})/(q+1, n)$
and is central in $A.\PSL_{n/2}(q^2)$ and $B$ has order $(q^2-1, \frac{n}{2})/(q+1,\frac{n}{2})$. Therefore, the
%Let $\varepsilon \in
%\{1,-1\}$ be
%such that $p$ divides $q^{n/2}- \varepsilon$. By \cite[Lemma~4.5.6]{KL}, %$G$ contains a subgroup
%isomorphic to $\PSp_n(q)$ extended by a group of order $(2,q-1)(q+1, %n/2)/(q+1,n)$.
 assertion follows by application of Lemma~\ref{help} (a) and (b) and the results for $\PSL_n(q)$.

Next suppose that $d_{q^2}(p) = n-1$. Then $p$ divides $q^{n-1} +1$ because by the table
in Lemma \ref{Primesalternativ} it does not divide $q^{n-1} -1$.
Here \cite[Lemma~4.1.4]{KL} tells us that $G$ has a subgroup
of structure $C_a\PSU_{n-1}(q)C_b.C_b$, where $a = q+1$ and $b = (q+1, n-1)$. Then the subgroup isomorphic to $\PSU_{n-1}(q)$ acts
trivially on $C_a$. Moreover $\PSU_{n-1}(q)$ has a subgroup of structure $C_c C_{n-1}$, where
$c = (q^{n-1}+1)/(q+1)(q+1,n-1)$ by \cite[Lemma~4.3.6]{KL}.
Then, using Lemma~\ref{help}(a), we obtain that $N$ is divisible by $D$ as given in the table.

Now suppose that $n$ is odd. If $d_{q^2}(p) = n$, then $p$  divides $q^n+1$.
In this case the claim follows as explained for $n$ even and
$d_{q^2}(p) = n-1$ by quoting \cite[Lemma~4.3.6]{KL}.
Finally suppose that $d_{q^2}(p) = (n-1)/2$. By \cite[Lemma~4.1.4]{KL}, there is a subgroup in $G$ of structure
$C_a. \PSU_{n-1}(q).C_b$, where $a = (q+1)/(q+1,n)$ and $b = (q+1,n-1)$. Now we apply our result that we obtained  for $G = \PSU_n(q)$, $n$ even and $d_{q^2}(p) = n/2$, to obtain the divisor as presented in the table.
\\
\underline{$G = \PSp_n(q)$:}

Let $\varepsilon \in \{ 1,-1\}$ be
such that $p$ divides $q^n- \varepsilon$.
According to \cite[Lemma~4.3.10]{KL}, $G$ contains a subgroup isomorphic to $\PSp_2(q^n).n \cong
\PSL_2(q^n).n$. The assertion follows in this case by Lemma~\ref{help}(a) and because $\PSL_2(p^n)$
contains a dihedral subgroup of order $2\cdot n \cdot (q^n-\varepsilon)/(q-1,2)$.
\\

We deviate from the ordering of the groups in the table now, because it makes our arguments easier.

\underline{$G = \POm_{2n}^+(q)$:}

First suppose that $d_q(p) = n$. Then $n$ is odd. We apply
\cite[Lemma~4.2.7]{KL} and \cite[Satz 7.2]{Huppert-I}:
If $q$ is even, then $G$ contains a subgroup of structure $\GL_n(q)C_{(n,2)}$, which provides us with the required divisor $D$.
If $q$ is odd, then we need to distinguish two cases
depending on whether $(q-1)/2$ is even or odd. In these cases
$G$ has a subgroup of structure $C_{(q-1)/a} \PSL_n(q)C_{(q-1, n)}$, where $a=4$ in the first case and $a=2$ in the second case. From this and the fact
that $(q^n-1,2) = (q-1, 2)$, we deduce the required
divisor.

Now suppose that $d_{q^2}(p) = n-1$. Then $n$ is even.  and $p$ divides $q^{n-1}-1$ according to Lemma~\ref{Primesalternativ}.
By \cite[Proposition~4.1.6]{KL}, $G$ contains a subgroup isomorphic to
$(\Omega^+_{2}(q) \times \Omega^+_{2n-2}(q))C_2$, and we
 find our divisor $D$ by applying the results for $d_q(p) = n-1$ to $\Omega^+_{2n-2}(q)$ and by recalling that $|\Omega^+_{2}(q) | = (q-1)/2$.
\\
\underline{$G = \POm_{2n}^-(q)$:}

If $d_{q^2}(p) = n$, then the normalizer of $T$ in $\Omega_{2n}^-(q)$ has order $(q^n+1)n/(2,q-1)$ by \cite[Proposition~4.3.15]{KL}. In this case, we derive the assertion from the fact that $|\Omega_{2n}^-(q): P\Omega_{2n}^-(q)| = (q^n+1,4)/(q-1,2)$.
It remains to consider the case where $d_{q^2}(p) = n-1$. According to \cite[Proposition~4.1.6]{KL}, and since $N$ is not divisible by $8$, the group $G$ has a subgroup isomorphic to $\Omega_2^-(q) \times
\Omega_{2n-2}^-(q)$. Now we apply \cite[Proposition~4.3.15]{KL} to the subgroup of $G$ isomorphic to $\Omega_{2n-2}^-(q)$, and  obtain the divisor $D$ of $N$ as written in the table.
\\
\underline{$G = \POm_{2n+1}(q)$:}

Then $q$ is odd and
$d_{q^2}(p) = n$.
By \cite[Proposition~4.1.6]{KL}, $G$ has subgroups of structure $\Omega_{2n}^\varepsilon (q). C_2$, where $\varepsilon \in \{1,-1\}$. If $n$ is even, then $p$ divides $q^n+1$ and we find the required divisor by
applying our results for $P \Omega_{2n}^-(q)$  and $d_{q^2}(p) = n$. If $n$ is odd, then either $p$ divides $q^n+1$ and we apply the results for $P \Omega_{2n}^-(q)$  and $d_{q^2}(p) = n$, once more, or $p$ divides $q^n-1$ and we apply the results for $P \Omega_{2n}^+(q)$  and $d_{q}(p) = n$.
\end{proof}

\begin{lemma}\label{IsPrimeLinearUnitaer}
	 Suppose that $(G,\Omega)$  satisfies Hypothesis~\ref{hyp6}.  If $n \geq 4$ and if $G = \PSL_n(q)$ or $G=\PSU_n(q)$,
	then $d_q(p)$ and $d_{q^2}(p)$ are prime numbers, respectively.
\end{lemma}
\begin{proof} First we note: If $a$ is a divisor of the natural number $l$, then $q^a -1$ divides $q^l-1$, because
	$$q^l-1 = (q^a -1)((q^{a})^{b-1} + \cdots + q^a + 1),~\mbox{where}~ b = l/a.$$
If $l$ is odd, then $q^a +1$ divides $q^l+1$, because
$$q^l+1 = (q^a+1)((q^{a})^{b-1} -  \cdots - q^a + 1).$$

Assume for a contradiction that  $G = \PSL_n(q)$, that $d_q(p)$ is not a prime and that $1 < a <  d_q(p)$ is a divisor of $d_q(p)$.
 Moreover, suppose that $d_q(p) = n$. We observe that $a \mid n$  and $n \geq 4$, which together yields that $a < n-1$.
%If $p$ divides $q^n-1$, then $q^a-1$ must divide $(q-1)(n,q-1)$ by Lemma~\ref{Primesalternativ}.
If $(q,a) = (2,6)$, then  we chose $a = 3$ instead of $a = 6$.
Therefore, there exists a primitive prime divisor $r$ of  $q^a-1$. Since $d_q(r) =  a >1$, we also have that $1 < d_q(r) <n-1 $. According to the table in Lemma \ref{tableTX} it follows that $|N_G(P)|$ is divisible by $r$. Since $r \neq 2$, we obtain that $r >3$ and that $r \in \pi(H)$ (see Lemma \ref{start}), contrary to Lemma~\ref{Primesalternativ}, or that $r = 3$,
contrary to Lemma~\ref{start}.

We conclude that $d_q(p) = n-1$.  Then we use the same argument as in the previous paragraph, just replacing $n$ by $n-1$, and  we obtain a similar contradiction. Thus $d_q(p)$ is a prime.

Now suppose that $G=\PSU_n(q)$. Then $d_{q^2}(p) = (n-1)/2$ or $n$, if $n$ is odd, and $d_{q^2}(p) =n-1$ or $n/2$ if $n$ is even.
Assume that  $d_{q^2}(p)$ is not a prime and that it is divisible by $1 < a <  d_{q^2}(p)$.
Suppose first that $n$ is odd. If $d_{q^2}(p) = n$, then $a$ is odd as well, $p$ divides $q^n+1$, and  $q^a+1$ also divides $q^n+1$. Then, as in the case where $G = \PSL_n(q)$, we arrive at a contradiction by choosing a primitive prime  divisor of $q^{2a} -1$ if $(q,a) \neq (2,3)$. If $(q,a) = (2,3)$, then $2^3 +1 = 9$ and  $(q+1)(4,q+1) = 3$, which  yields a contradiction to Lemma~\ref{start}(f).
Thus our assumption implies that $d_{q^2}(p) = (n-1)/2$, by Table \ref{ktablealternativ}  Then $p$   and $q^{2a}-1$ divide $q^{n-1} - 1$.
 If $(q,a) = (2,3)$ ,  then we notice that $n-1$ divides $|N_G(P)|$ by Lemma \ref{tableTX}, so we arrive at the contradiction that $3$ divides $|N_G(P)|$ (see Lemma~\ref{start}(f)).
 Thus $(q,a) \neq (2,3)$ and we can chose a primitive prime divisor $r$ for $q^{2a}- 1$.  Since $2a \geq 4$, we see that $r \in \pi(N_G(P))$ by the table in Lemma \ref{tableTX}. Consequently $r \in \pi(H)$, which yields a contradiction to Lemma~\ref{Primesalternativ}, because $d_{p^2}(r) = a< (n-1)/2$.

We conclude that $n$ is even.  Then $d_{q^2}(p) = n-1$ or $n/2$.  In the first case $p$ divides $q^{n-1 } + 1$ and in the second case it divides $q^n -1$.
We then derive a contradiction in a similar way as for the case where $n$ is odd.
\end{proof}

\begin{lemma}\label{Auch1weniger}
	Let $G$ be a simple classical group of Lie type, but none of the groups $\PSL_n(q)$ or $\PSU_n(q)$, and suppose that $(G, \Omega)$ satisfies Hypothesis~\ref{hyp6}.	
%If $d_q(p)$  or $d_{q^2}(p)$ are given as in Table..,
Then $n \in \{2,4\}$ or
there is a prime $r \in \pi(H)$ such that
$d_q(r) $  or $d_{q^2}(r)$ equals $n-1$.
\end{lemma}
\begin{proof}
If $d_q(p) = n-1$ or $d_{q^2}(p) = n-1$, then our assertion holds. Therefore assume that $d_q(p) = n$ or $d_{q^2}(p) = n$.
 According to Table~\ref{TX}, the order of $N_G(P)$ is
divisible by the odd part of $n$.
Assume there is an  odd prime divisor $r$ of $n$.  Then  $r \in \pi(H)$ by Lemma~\ref{start}(a). Since $n$ divides $q^{n-1} -1$,  it follows that
$d_q(r) < n$ or $d_{q^2}(p) < n$ and we read from Table \ref{TX} that $d_q(r) = n-1$ or $d_{q^2}(p) = n-1$.

Now assume that $n$ is a $2$-power and that there does not exist a prime $r \in \pi(H)$ such that $d_q(r) = n-1$ or $d_{q^2}(r) = n-1$.
Then according to Lemma~\ref{Primesalternativ} $G = \PSp_{2n}(q)$ or $G=P\Omega_{2n+1}(q)$.  In both cases $n$ divides  $|N_G(P)|$. Therefore, an application of Lemma~\ref{start}(a) yields the assertion.
\end{proof}

\begin{lemma}\label{symplectic}
	Suppose that $G = \PSp_{2n}(q)$, where $ n \geq 2$, or that $G=P\Omega_{2n+1}(q)$, where $n \geq 3$. Suppose further that $(G, \Omega)$ satisfies Hypothesis~\ref{hyp6}.	
	Then  $G = \PSp_{4}(q)$ and $G_\alpha = H $  is a cyclic group of order $(q^2+1)/(2,q+1)$.
\end{lemma}

\begin{proof}
Lemmas~\ref{Auch1weniger} and \ref{Primesalternativ} imply that $n \in \{2,4\}$.
First  assume that $G =  P\Omega_{2n+1}(q)$.  Then $n \geq 3$ yields that $n =4$. Further recall that $q$ is an odd prime power.
Since $4$ is even and $d_{q^2}(p) = 4$, it follows that  $p$ divides $q^4+1$. Then we turn to the table in Lemma \ref{tableTX} and we see that
$|N_G(P)|$ is divisible by $(q^4+1)4$.  Since $q$ is an odd prime power, it follows that $8$ divides $|N_G(P)|$, contradicting Lemma~\ref{start}(f).
This shows that $(P\Omega_{2n+1}(q), \Omega)$ does not satisfy Hypothesis~\ref{hyp6} for any set $\Omega$.

It remains to consider the situation where $G=\PSp_{2n}(q)$. Since $n$ is even, we see that $p$ divides $q^n +1$. Then
$|N_G(P)|$ is  divisible by $(q^n +1 )2n/(2,q-1)$, which is divisible by $2n$. This forces $n = 2$ (see Lemma~\ref{start}(f)), and  $G_\alpha= H = O(N_G(P))$ is as stated.
\end{proof}

We remark that  $(G, G/U)$, where $G = \PSp_{2n}(q)$,  $U$  is cyclic of order $(q^2+1)/(2,q+1)$ and $G$ acts on $G/U$ by right multiplication satisfies Hypothesis~\ref{hyp6}, by Lemma \ref{tableTX}.

\begin{prop}\label{ClassicalNotUnitaryCo6}
Let $G$ be a simple classical group of Lie type, and suppose that $(G, \Omega)$ satisfies Hypothesis~\ref{hyp6}.	
Then one of the following holds:
\begin{itemize}
	\item[(a)]  $G= \PSL_2(q)$ is as in Lemma~\ref{L2q-main}(f).
		\item[(b)] $G= \PSU_4(q)$, where $q \in \{2,3\}$, and $G_\alpha$ is cyclic of order $5$.
	\item[(c)] $G= \PSp_4(q)$, and $G_\alpha$ is cyclic of order $(q^2 +1) /(2, q-1)$.
	\item[(d)] $G= \POm^-_8(q)$, and  $G_\alpha$ is cyclic of order
	$(q^4+1)/(2,q+1)$.
\end{itemize}	
\end{prop}
\begin{proof} The symplectic and the odd-dimensional orthogonal groups have already been discussed in Lemma~\ref{symplectic}, and small Lie type groups have been treated in Lemma \ref{L2q-main}, so we obtain (a) and (c). Now we may suppose that $G \in\{ \PSL_n(q), \PSU_n(q)\}$ or
$G = \POm^{\varepsilon}_{2n}(q)$, where $\varepsilon \in \{1,-1\}$, and $n \geq 4$ in all cases.
\medskip
\\		
Suppose that  $G =  \PSL_n(q)$ or $G = \PSU_n(q)$, where $n \geq 4$.  Then by  Lemma~\ref{IsPrimeLinearUnitaer}
	$d_q(p)$ or $d_{q^2}(p)$, respectively, is a prime number.  Moreover,  Lemma~\ref{Primesalternativ} tells us that $d_q(p)$ and $d_{q^2}(p)$ can have two different values.
	
	Assume first that $G = \PSL_n(q)$.  If $d_q(p) = n$, then according to \cite[Kap. II, Theorem~(7.3)(c)]{Huppert-I} we have that
		$$|N_G(P)| = n\cdot [(q^n-1)/m(q-1)]~ \mbox{, where}~ m = ((q^n-1)/(q-1), n).$$
		As $n \geq 4$ and $n=d_q(p)$ is a prime, it is an odd prime. This implies that $(q^n-1)/(q-1) =  q^{n-1} + \cdots + q +1$ and $|N_G(P)|$  are odd numbers, which contradicts Lemma~\ref{start}(a).   Then we conclude that $d_q(p) = n-1$ and that $n-1 \geq 3$ is a prime.  As $|N_G(P)|$  is divisible by $n-1$ , see
		 \cite[Kap. II, Theorem~(7.3)(c)]{Huppert-I}, we deduce from Lemma~\ref{Primesalternativ} that $d_q(n-1) = n-1$, contrary to Lemma~\ref{help}(c).
\medskip
\\		
Next suppose that $G = \PSU_n(q)$ and assume, in addition, that  $n \geq 4$ is odd. If $d_{q^2}(p) = n$, then by \cite[Corollary~5 and the proof of  Theorem~5.3]{Hest}
		we have that $|N_G(P)| = (q^n+1)n/[(q+1)(n,q+1)]$.  This shows that $|N_G(P)|$  is odd and and hence $|N_G(H)|$ is odd, contrary to Lemma~\ref{help}(a). This yields that
		$d_{q^2}(p) = (n-1)/2$.
		
		We derive from Lemma~\ref{IsPrimeLinearUnitaer}
		that $(n-1)/2 \geq 2$ is a prime dividing $|N_H(P)|$. Then Lemmas~\ref{start} and \ref{help} imply that $(n-1)/2 = 2$ and $n=5$.
According to Lemma~\ref{tableTX} $|N_G(P)|$ is divisible
  by $(q^4-1)\cdot  4/(q+1,5)$.
  Since $N_G(P)$ is a subgroup of the stabilizer of $\Delta$, Lemma~\ref{help} (b) yields that  $q+1$ is a divisor of $24$. Therefore
  $q+1 = 4$ by Lemma~\ref{start} (f), which yields that $8$ divides $D$ in contradiction to Lemma~\ref{start} (f).
 % Now $p$ divides $q^2+1$
%		and $|N_G(P)|$ is divisible by $(q^4 - 1)4/(5,q+1)$ by Lemma~\ref{sgpPSU}. Therefore, $q$ is a $2$-power, and $|N_G(P)|$ as well as $|N_G(H)|$ are divisible by $2^2-1 = 3$, which contradicts
%		Hypothesis~\ref{hyp6}.
		This shows that $n$ is even.

		Suppose next that  $d_{q^2}(p) = n -1$. Then $n-1$ is a prime by Lemma~\ref{IsPrimeLinearUnitaer}, which divides $|N_G(P)|$ by Lemma~\ref{tableTX}. Then Lemma~\ref{help}(c)
		yields that $d_{q^2}(n-1) = n/2$. By Lemma~\ref{IsPrimeLinearUnitaer}  $d_{q^2}(p) = n/2 \geq 2$ is a prime,
  and it divides $|N_G(P)|$ by Lemma~\ref{tableTX}. If $n/2 > 2$, then $n/2 \neq 3$ by Lemma~\ref{start} (f) and $d_{q^2}(n/2) < n/2$ by Lemma~\ref{help} (c), contrary to Lemma~\ref{Primesalternativ}.
  This shows that  $n/2 = 2$  and $n = 4$.
  Then application of Lemma~\ref{TX} gives that $|N_H(P)|$ is divisible by
$$(q^4-1)/(q+1)(q+1,4)= (q^2+1)(q-1)/(q+1,4).$$ By Lemma~\ref{start} (g) also $|H|$ is divisible by that number.
Therefore, Lemma~\ref{start} (b) yields that $(q-1)$ divides $(q+1,4)$. Since $(q-1,q+1) \in \{1,2\}$, we conclude $q-1 \in \{1,2\}$ and $q \in \{2,3\}$.
%		Therefore $q-1 \in \{1,2,4\}$ and $q+1 \in \{3,4,6\}$.  If $n/2 \geq 3$, then Lemma \ref{tableTX} implies that $n/2$ divides $q+1$. We conclude that $(q,n) = (2,6) $
%		and we obtain the contradiction  that $|N_G(T)|$ is divisible by $3$. This shows that $n/2 = 2$. If $q= 5$, then again $|N_G(T)|$ is divisible by $3$ by  Table \ref{TX}.
	Therefore statement (b) holds.
\medskip
\\		
Now suppose that $G = \POm^{\varepsilon}_{2n}(q)$ .
		If $\varepsilon = +$ and $d_q(p) = n$ or $\varepsilon = -$ and $d_{q^2}(p) = n$, then
		Lemma~\ref{Auch1weniger} implies that
	 $n \in \{2,4\}$ or that there is $r \in \pi(H)$ such that $d_q(r) = n-1$ or $d_{r^2}(p) = n-1$. If $n \in \{2,4\}$, then $n \geq 4$ gives that $n = 4$.

	 Assume first that $G =P\Omega_{2n}^+(q)$.
	If $d_{q^2}(r) = n-1$ for some $r \in \pi(H)$, then $n$ is even by Lemma~\ref{TX}. This shows that in both cases, $n \in \{2,4\}$ or  $d_{q^2}(p) = d_q(p) = n-1$,
	Lemma \ref{tableTX} implies that  $|N_G(P)|$ is divisible by $n-1$. Now we argue as before: $n-1$ is
 odd, and therefore divisible by some prime $t$. According to Lemma~\ref{help} (c) $d_q(t) \leq t-1 < n-1$ in contradiction to Lemma~\ref{Primesalternativ}.
This shows that $G \neq \POm_{2n}^+(q)$.

Thus $G = \POm_{2n}^-(q) $ and $d_{q^2}(p) \in \{n-1, n\}$. Assume that $n$ is not a $2$-power. Then, by Lemma~\ref{Auch1weniger}, there is a prime $r \in \pi(H)$ such that $d_{q^2}(r) = n-1$. Therefore, if  $T \in \Syl_r(H)$, then $(n-1)$ divides $|N_G(P)|$ by Lemma \ref{tableTX}. Hence we see, just as in the proof of
Lemma~\ref{Auch1weniger}, that $n-1$ is in $\{2,4\}$, which yields again that $n-1 = 4$ and $n=5$. Then $|N_G(P)|$ is divisible by $3 = 2^2+1$, if $q$ is even, and by $8$, if $q$ is odd, by \cite[Table~8.68]{BHRD}.
	
This shows that $n \in \{2,4\}$ and therefore $n = 4$. Then $p$ divides $q^4+1$, and $G_\alpha$ is cyclic of order  $(q^4+1)/(4,q^4+1) = (q^4+1)/(2,q+1)$, which is (d).
\end{proof}

%%%%%

\subsection{Exceptional groups}

\begin{hyp}\label{exc}
In addition to Hypothesis \ref{hyp6}, we suppose that
$q$ is a prime power and $G$ is one of the following groups:
$^3D_4(q)$, $F_4(q)$, $^2F_4(q)$, $G_2(q)$, $^2G_2(q)$, $E_6(q)$, $^2E_6(q)$, $E_7(q)$ or $E_8(q)$.
\end{hyp}

\begin{lemma}\label{3D4q}
Suppose that Hypothesis \ref{exc} holds and that $G={}^3D_4(q)$.
Then $G_\alpha$ is cyclic of order $q^4-q^2+1$.
\end{lemma}

\begin{proof}
Lemma \ref{Primesalternativ} tells us that $p$ divides $\Phi_{12}(q)=q^4-q^2+1$.
Inspection of Theorem~4.3 in \cite{Wilson} shows that $G$ has a maximal
	subgroup \(M\) of structure \(C_{q^4-q^2+1}:4\). Let \(C\leq M\) be cyclic of order
	\(q^4-q^2+1\). Then without loss \(P \le C\), in fact $P$ is characteristic in $C$ and hence normal in $M$.
	Then $M=N_G(P)$, because $G$ is simple.
	Together with Lemma \ref{start}(e) this implies that
	\(|M:N_{H}(P)|\leq 4\) and therefore \(|H|\) is divisible by \(q^4-q^2+1\). The list in Theorem~4.3
	in~\cite{Wilson} only leaves the possibility that $H$ is contained in a maximal subgroup of $G$ that is isomorphic to $M$, and then it follows that
	\(G_{\alpha}\) is cyclic of order \(q^4-q^2+1\).
\end{proof}

\begin{lemma}\label{F4q}
Suppose that $G=F_4(q)$.
Then there is no set $\Omega$ such that Hypothesis \ref{exc} holds.

\end{lemma}

\begin{proof}
Assume for a contradiction that Hypothesis \ref{exc} holds for $G$.
Then
Lemma \ref{Primesalternativ}  implies that $p$ divides $\Phi_8(q)=q^4+1$ or $\Phi_{12}(q)=q^4-q^2+1$.

In the second case, where \(p\) divides \(q^4-q^2+1\), we consider Table~5.1
	in~\cite{LSS} and we see that \(G\) has a subgroup
	\(M\) of structure \({}^3\textrm{D}_4(q).3\). Let \(D\leq M\) be isomorphic to \({}^3\textrm{D}_4(q)\).
	Then Theorem~4.3
	in~\cite{Wilson} yields that \(D\) has a cyclic subgroup \(C\) of order
	\(q^4-q^2+1\), and without loss $P \le C$. Then
	\(|N_D(P)|=4(q^4-q^2+1)\) and, since $p\ge 5$, Lemma \ref{help}(a) yields that \(|N_M(P)|=3\cdot
	|N_D(P)|=12\cdot (q^4-q^2+1)\), contrary to Lemma \ref{start}(f).
	
	Therefore, the first case above holds, i.e. \(p\) divides \(q^4+1\).
	
	Then we use Table~5.1 in~\cite{LSS} and we see that \(G\) has
	a subgroup \(M\) of structure \((2,q-1).\Omega_9(q)\).
	Let $Z \unlhd M$ be such that $|Z|=(2,q-1)$ and let $Q \in \syl_p(M)$.  As $p$ and $|Z|$ are coprime,
	the coprime action of $Q$ on $Z$ implies that  \(N_{M/Z}(QZ/Z)=N_M(Q)Z/Z=N_M(Q)/Z\).
	We deduce that
	 \(|N_G(Q)|\) is divisible by \(|N_{M/Z}(QZ/Z)|\cdot |Z|\).
	Now \(M/Z \cong\Omega_9(q)\), and we can refer to the corresponding lemmas.
	If \(q\) is even, then \(\Omega_9(q)\cong \PSp_8(q)\), where the normalizer of a Sylow $p$-subgroup has order
	divisible by $8$. This contradicts Lemma \ref{start}(f).
	
	If $q$ is odd, then we turn to
 Table \ref{TX} for the group $\POm_{2n+1}(q)$, where
	the normalizer of a Sylow $p$-subgroup also has order
	divisible by $8$, again contrary to Lemma \ref{start}(f). This final contradiction proves the lemma.
\end{proof}

\begin{lemma}\label{2F4q}
Suppose that  $G=\,^2F_4(q)$ or that $q=2$ and
$G=\,^2F_4(q)'$.
Then there is no set $\Omega$ such that Hypothesis \ref{exc} holds for $G$ and $\Omega$.
\end{lemma}

\begin{proof}
If $q=2$ and
$G=\,^2F_4(q)'$, then our claim follows
from Remark~\ref{rem4.1}.

Next assume that $G=\,^2F_4(q)$ and that Hypothesis \ref{exc} holds.
Lemma \ref{Primesalternativ}  gives that $p$ divides
 \(\Phi_6(q)=q^2-q+1\) or \(\Phi_{12}(q)=q^4-q^2+1\).
	
In the first case, where \(p\) divides \(q^2-q+1\), the Main Theorem in~\cite{Malle2F4}
tells us that $G$ has a maximal subgroup \(M\) of structure \(\SUn_3(q):2\). Let \(S\leq M\) be isomorphic to
	\(\SUn_3(q)\) and, without loss, let  \(P\le S\).
	Since \(q +1 \equiv 2^{2n+1} +1 \equiv (-1)^{2n+1} +1 \equiv -1 +1 \equiv 0 \mod 3\), we see that \((3,q+1)=3\) and hence
	\(Z(S)\) is cyclic of order \(3\). Therefore \(|N_S(P)|\) is divisible by \(3\). Moreover $p \ge 5$ by hypothesis and
	Lemma \ref{help}(a) yields that \(|N_M(P)|=2\cdot |N_S(P)|\), whence \(|N_G(P)|\)
	has order divisible by \(6\). This contradicts Lemma \ref{start}(f).
	
	In the second case, where \(p\) divides \(q^4-q^2+1\), we first observe that
	\(q^4-q^2+1=(q^2+q+1+\sqrt{2q}(q+1))(q^2+q+1-\sqrt{2q}(q+1))\). Let \(\varepsilon \in \{-1,1\}\) be such that
	\(p\)
	divides \((q^2+q+1+\varepsilon\sqrt{2q}(q+1))\). Then the Main Theorem in~\cite{Malle2F4}
	shows that \(G\) has a maximal subgroup \(M\) of structure \(C_{(q^2+q+1+\varepsilon\sqrt{2q}(q+1))}:12\).
	Without loss $P \le M$ and then $|N_G(P)|$ is divisible by \(12\), contrary to Lemma \ref{start}(f).
\end{proof}

\begin{lemma}\label{G2q}
There is no set $\Omega$ such that Hypothesis \ref{exc} holds in the case where $G=G_2(q)$.
\end{lemma}

\begin{proof}
We assume otherwise. Then Lemma \ref{Primesalternativ}  yields that every
\(p\in \pi(H)\) divides \(\Phi_{3}(q)=q^2+q+1\) or \(\Phi_6(q)=q^2-q+1\). Let
	\(\varepsilon\in \{-1,1\}\) be such that \(p\) divides \(q^2+\varepsilon q+1\). Then \(H\) contains a Sylow
	\(p\)-subgroup $P$
	of \(G\) by Lemma \ref{start}(e) and \(P\) lies in a maximal subgroup \(M\) of structure
	\(\SLi^{\varepsilon}_3(q):2\) by Table~4.1
	in~\cite{Wilson}.
	Let \(S\leq M\) be a subgroup of index \(2\) in \(M\) such that \(S \cong \SLi^{\varepsilon}_3(q)\).
	Then \(P\leq S\) and Lemma \ref{help}(a) implies that \(|N_M(P)|=2\cdot |N_S(P)|\).
	If \(Z(S)\neq 1\), then \(|Z(S)|=(3,q-\varepsilon)=3\) and hence \(|N_S(P)|\) is divisible by \(3\).
	Hence \(|N_G(P)|\) is divisible by \(6\).
	If \(Z(S)=1\), then \(S\cong \PSL^{\varepsilon}_3(q)\). Then Theorem~6.5.3
	in~\cite{GLS3} tells us that \(S\) has a subgroup \(F\) that is a Frobenius group with Frobenius
	kernel
	\(K\) of order \(q^2+\varepsilon q+1\) and Frobenius complement of order \(3\). Since \(p\geq 5\), it follows that
	\(P\leq K\), in fact \(P=O_p(K)\), because Frobenius kernels are nilpotent. Hence \(N_F(P)=F\) and
	we see that $3$ divides  \(|N_F(P)|\), which means that
	\(|N_G(P)|\) is divisible by \(6\) in this case as well.
	In both cases Lemma \ref{start}(f) gives a contradiction.
\end{proof}

\begin{lemma}\label{2G2q}
Suppose that Hypothesis \ref{exc} holds and that $G=\,^2G_2(q)$.
Then $G_\alpha$ is cyclic of order $\frac{q-1}{2}$.
\end{lemma}

\begin{proof}
Lemma \ref{Primesalternativ} gives that $p$ divides
$\Phi_1(q)=q-1$, $\Phi_2(q)=q+1$ or
 \(\Phi_6(q)=q^2-q+1\).

We begin with the last case and assume that \(p\) divides \(q^2-q+1=(q+\sqrt{3q}+1)(q-\sqrt{3q}+1)\).
Then there exists \(\varepsilon\in \{-1,1\}\)
	such
	that \(p\) divides \(q-\varepsilon\sqrt{3q}+1\).
	Theorem~4.2 in \cite{Wilson} yields that $G$ has a
	subgroup \(M\) of structure \(C_{q-\varepsilon\sqrt{3q}+1}:6\).
	Then without loss $P \le M$ and
	 \(|N_G(P)|\) is divisible by~\(6\), which is impossible by Lemma \ref{start}(f).
	
	Next we assume that \(p\) divides \(q+1\). Then, by Theorem~4.2
	in~\cite{Wilson}, $G$ has a subgroup \(A\) of structure
	\(2^2\times D_{(q+1/2)}\). Let \(E\) be elementary abelian of order~\(4\) and let \(D\) be  a dihedral group of order
	\(\frac{q+1}{2}\) such that \(E\times D\leq A\). Moreover, let \(Y \le D\) be a subgroup of order \(p\). Then
	\(Y \unlhd D\) and consequently \(N_A(Y)=E\times D=A\). As \(|D|\) is even, we see that \(|N_A(Y)|\) is divisible
	by \(8\). Since we can choose $Y$ as a subgroup of $P$, this contradicts Lemma \ref{start}(f) as well.
	We conclude that \(p\) divides \(q-1\).
	
	Theorem~4.2 in~\cite{Wilson} shows that $G$ has a maximal subgroup \(M\) of structure \(2\times\PSL_2(q)\). Let \(C\) be cyclic of order \(2\) and let \(L \cong \PSL_2(q)\) be such that \(C\times L\leq M\). Let \(Y\le L\) be a subgroup of order \(p\).
	Then Satz~II\,8.3
	in~\cite{Huppert-I} shows that  \(N_L(Y)\) is a dihedral group of order \(2\cdot
	\frac{q-1}{2}\), whence \(N_M(Y)=C\times N_L(Y)\) has order \(2\cdot 2\cdot
	\frac{q-1}{2}\). As before we may choose $Y$ as a subgroup of $P$ and we deduce with Lemma \ref{start}(e) that
	 \(|H|\) is divisible by
	\(\frac{q-1}{2}\). Inspection of Hauptsatz~II\,8.27
	in~\cite{Huppert-I} shows that  the only subgroups of \(\PSL_2(q)\) of order coprime to \(6\) and divisible
	by \(\frac{q-1}{2}\) are cyclic of order \(\frac{q-1}{2}\). Combining this with
	the list of maximal subgroups of \(G\) in
	Theorem~4.2 in~\cite{Wilson} reveals that \(G_{\alpha}\) is cyclic of order
	\(\frac{q-1}{2}\).
\end{proof}

\begin{lemma}\label{PaulasubgroupsE6}

Suppose that \(q\) is a prime power, that
 $G=E_6(q)$ (and then $\varepsilon=1$) or $G={^2}E_6(q)$ (and then $\varepsilon=-1$) and that
		\(p\geq 5\) is a prime that divides neither \(q^6-1\) nor \(q^4-1\).
		 If \(p\) divides one of the following numbers (left colomn),
		then in each
		case there exists a \(p\)-subgroup \(Q\) of \(G\)
		such that \(|N_G(Q)|\) is divisible by the corresponding number in the right colomn of the table.\\
		
		\begin{center}
			\begin{tabular}{l|c}
				\(p\) divides & \(|N_G(Q)|\) is divisible by \\
				\hline
				\(q^4-q^2+1\)&\(3\)\\
				\(q^6+\varepsilon q^3+1\)& \(3\)\\
				\(q^4+1\)& \(8\)\\
				\(q^4+\varepsilon q^3+q^2+\varepsilon q+1\) %
				&\(5\)
			\end{tabular}
		\end{center}
		\end{lemma}
		
	\begin{proof}
	If \(p\) divides \(q^4-q^2+1\), then we use Table~5.1
	in~\cite{LSS}. We see that \(G\) has a
	subgroup~\(M\) of structure \(({^3}D_4(q)\times C_{q^2+\varepsilon q+1}).C_3\) and we choose a normal subgroup
	\(N\) of~\(M\) of
	index~\(3\).
	If \(p\) divides \(q^6+\varepsilon q^3+1\), then we set \(e:=(3,q-1)\) and we use Table~5.1
	in~\cite{LSS} again. Here \(G\) has a subgroup
	\(M\) of structure \(\PSL^{\varepsilon}_3(q^3).(C_e\times C_3)\). In this case, let \(N\leq M\) be isomorphic to
	\(\PSL^{\varepsilon}_3(q^3)\).
	Then in both cases, \(p\) divides \(|N|\) and \(|M:N|\) is divisible by \(3\). Let \(Q\in \Syl_p(N)\).
	By Lemma~\ref{help}~(a), it follows that $3$ divides \(|N_M(Q)|\), as stated in the table.
		
	Next suppose that \(p\) divides \(q^4+1\) or \(q^4+\varepsilon q^3+q^2+\varepsilon q+1\) and
	let \(h:=(4,q-1)\). Then Table~5.1
	in~\cite{LSS} gives that \(G\) has a subgroup
	\(M\) of structure \(C_h.(\POm^{\varepsilon}_{10}(q)\times C_{(q-\varepsilon)/h})\).
	Let \(Y\) be an arbitrary \(p\)-subgroup of \(M\) and let \(Z\) be a normal subgroup of \(M\) of
	order~\(h\). Since \(p\geq 5\), we see that \(Y\) acts coprimely on \(Z\) and then
	\(N_{M/Z}(YZ/Z)=N_M(Y)Z/Z\).
	Therefore \(|N_M(Y)|\) is divisible by \(|N_{M/Z}(YZ/Z)|\). In particular, since \(M/Z\) has a subgroup isomorphic to
	\(\POm^{\varepsilon}_{10}(q)\), we deduce that, for every
	\(p\)-subgroup \(R\) of \(\POm^{\varepsilon}{10}(q)\), there exists a
	\(p\)-subgroup \(Y\) of \(M\) such that \(|\N_G(Y)|\) is divisible by
	\(|\N_{\POm^{\varepsilon}_{10}(q)}(R)|\).
	
	To specify our analysis further, we suppose that
	 \(p\) divides \(q^4+1\). Then \(p\) divides \(q^8-1\) and, if \(l<4\), then our hypothesis guarantees that \(p\)
	does not divide
	\(q^{2l}-1\). Now the internal structure of \(\POm^{\varepsilon}_{10}(q)\)
	gives a \(p\)-subgroup~\(R\) of \(\POm^{\varepsilon}_{10}(q)\) such that
	\(|N_{\POm^{\varepsilon}_{10}(q)}(R)|\) is divisible by \(\frac{(q^4+1)\cdot 2 \cdot 4 \cdot
		(q+\varepsilon)}{(4,q^5-\varepsilon)}\). If \(q\) is even, then \(\frac{(q^4+1)\cdot 2 \cdot 4 \cdot
		(q+\varepsilon)}{(4,q^5-\varepsilon)}=(q^4+1)\cdot 2 \cdot 4 \cdot (q+\varepsilon)\), which is a number divisible
	by \(8\). If
	\(q\) is odd, then \(q^4+1\) and \(q+\varepsilon\) are both divisible by \(2\), and hence \(\frac{(q^4+1)\cdot 2 \cdot
		4
		\cdot
		(q+\varepsilon)}{(4,q^5+\varepsilon)}\) is divisible by \(\frac{2\cdot 2 \cdot 4 \cdot 2}{4}=8\).
	Consequently, if \(p\)
	divides \(q^4+1\), then the information in the table is correct.
	
	Next we suppose that \(p\) divides \(q^4+\varepsilon q^3+q^2+\varepsilon q+1\). In
	particular, \(p\) divides
	\(q^5-\varepsilon\), which divides \(q^{10}-1\). If \(p\) divides \(q^4+1\), then \(p\) also divides \(q^8-1\), which divides
	\((q^{10}-1)-(q^8-1)=q^8\cdot (q^2-1)\). %
	By hypothesis \(p\) divides neither \(q^6-1\) nor \(q^4-1\), and this
	implies,
	for all \(l<5\), that \(p\) does not divide \(q^{2l}-1\). Again we use the subgroup structure of \(\POm^{\varepsilon}_{10}(q)\), this time we let \(R\in \Syl_p(\POm^{\varepsilon}_{10}(q))\). Then \(|N_{\POm^{\varepsilon}_{10}(q)}(R)|\) is divisible by \(\frac{5\cdot
		(q^5-\varepsilon)}{(4,q^5-\varepsilon)}\), and this number is divisible by \(5\). Again this is what we state in the table.
	
\end{proof}

\begin{lemma}\label{Paulacopr6E6}
Suppose that  $G=E_6(q)$ or $G={^2}E_6(q)$. Then there is no set $\Omega$ such that
Hypothesis \ref{exc} holds.

%Let \(q\) be a prime power, \(\varepsilon\in\{-1,1\}\), and		\(G=\textrm{E}^{\varepsilon}_6(q)\). Then \(G\) does not act transitively, %		with fixity~\(4\), and such that point stabilizers have order coprime to \(6\) on any set.
		\end{lemma}
		
	\begin{proof}
	We will prove the lemma by contradiction. Therefore, assume that there exists a set \(\Omega\) such that
	Hypothesis  \ref{exc} holds, with all its notation.
	Next we look at the values in the table in Lemma \ref{Primesalternativ} for $E_6(q)$ and ${^2}E_6(q)$.
We notice that \(\Phi_{10}(q)=\Phi_5(-q)\) and \(\Phi_{18}(q)=\Phi_9(-q)\), and then the table gives that
there is $\varepsilon \in \{1,-1\}$ such that
	\(p\) divides \(\Phi_5(\varepsilon q)=q^4+\varepsilon q^3+q^2+\varepsilon
	q+1\), \(\Phi_{8}(q)=q^4+1\),
	\(\Phi_9(\varepsilon q)=q^6+\varepsilon q^3+1\), or \(\Phi_{12}(q)=q^4-q^2+1\). In particular, $p$ and $q$ are coprime.
	
	We assume, for a contradiction, that a prime divisor \(r\) of \(|H|\) divides \(q^6-1\) or \(q^4-1\).
	Then the smallest positive integer \(k\) such
	that \(r\) divides \(\Phi_k(q)\) is in \(\{1,2,3,4,6\}\).
	Thus,
	\mbox{(10-2)}~in~\cite{GL} implies that
	the \(r\)-rank of
	\(G\) is at least \(2\). This contradicts Theorem~\ref{strongcases}.
	Therefore, \(|H|\) and \(q^6-1\) are coprime, as are \(|H|\) and \(q^4-1\), and we can apply
	Lemma~\ref{PaulasubgroupsE6}.
	
	If \(p\) divides \(q^4- q^2+1\) or \(q^6+\varepsilon q^3+1\), then the lemma gives a
	\(p\)-subgroup \(R\)
	such that \(|N_G(R)|\) is divisible by
	\(3\), contrary to Lemma~\ref{start}\,(f).
	
	If \(p\) divides \(q^4+1\), then Lemma~\ref{PaulasubgroupsE6} gives a \(p\)-subgroup~\(R\) such that
	\(|N_G(R)|\) is divisible by \(8\), which again contradicts Lemma~\ref{start}\,(f).
	
	As a consequence, \(p\) divides \(q^4+\varepsilon q^3+q^2+\varepsilon q+1\).
	Then Lemma~\ref{PaulasubgroupsE6} gives a \(p\)-subgroup \(R\) such that \(|N_G(R)|\) is divisible by
	\(5\) and, by Lemma~\ref{start}\,(f), it follows that $5 \in \pi(H)$.
	Now
	\(q^4-1\) is divisible by \(5\) by Lemma \ref{help}\,(c), which contradicts the fact that \(|H|\) and \(q^4-1\) are coprime.
	
\end{proof}

\begin{lemma}\label{E7q}
There is no set $\Omega$ such that Hypothesis \ref{exc} holds in the case where $G=E_7(q)$.
\end{lemma}

\begin{proof}
We assume otherwise and
consider Table~5.1
	in~\cite{LSS}.
	If $p$ divides $|\POm^{\varepsilon}_{12}(q)|$, then (using Lemma \ref{start}\,(e)) a Sylow $p$-subgroup of $H$ contains a Sylow $p$-subgroup of
	a subgroup of $G$ of structure $d.(\PSL_2(q) \times \POm^{\varepsilon}_{12}(q))$, where $d=(2,q-1)$. Now we notice that $d \in \{1,2\}$ and that $3$ divides $|\PSL_2(q)|$, and this contradicts Lemma \ref{start}\,(f).
This excludes the cases where $p$ divides $q^{2i}-1$, where $i \in \{1,\ldots, 5\}$.  	

In a similar way, we rule out more cases from the table in Lemma \ref{Primesalternativ}
by using Lemma \ref{Paulacopr6E6} and the possibility that
$p$ divides $|E_6(q)|$ or $|{^2}E_6(q)|$. Then, as above, a Sylow $p$-subgroup of $H$ contains a
Sylow $p$-subgroup of a subgroup $M$ of $G$ of structure
$(3,q-1).(~E_6(q) \times (q-(1/(3,q-1))~).(3,q-1)$ or
$(3,q+1).(~{^2}E_6(q) \times (q+(1/(3,q+1)) ~).(3,q+1)$.

Then we find a $p$-subgroup $Q$ of $M$ such that its normalizer has order divisible by
\(\frac{5\cdot (q-1)}{(4,q^5-1)}\) or
	\(\frac{5\cdot (q+1)}{(4,q^5+1)}\), respectively.
Calculations show that these numbers are always at least 5, which by Lemma \ref{start}(e) means that
there is a prime $r \in \pi(H)$ dividing them.
	
	First assume that \(r=5\). Since \(q\neq 5\) by Lemma \ref{Lie1}, Fermat's Little Theorem
	implies that \(q^4-1\) is divisible by \(5\). Then Theorem~(10-1)
	in~\cite{GL} shows that \(G\) has
	\(5\)-rank at least~\(2\), contradicting Theorem \ref{strongcases}.	
We recall that $r \notin \{2,3\}$, and hence $r \ge 7$ and $r$
must divide $q-1$ or $q+1$. Then Theorem~(10-1)
	in~\cite{GL} and Theorem \ref{strongcases} give another contradiction.
	
	The previous paragraphs exclude amost all possibilities from the table in Lemma \ref{Primesalternativ}, only
	\(\Phi_7(q)=q^6+q^5+q^4+q^3+q^2+q+1\) and \(\Phi_{14}(q)=q^6-q^5+q^4-q^3+q^2-q+1\) remain.
	We let \(\varepsilon\in \{-1,1\}\) be such that \(p\)
	divides \(q^6+\varepsilon q^5+q^4+ \varepsilon q^3 +q^2+ \varepsilon q+ 1\),we set
	 \(f:=\frac{(4,q-\varepsilon)}{(2,q-1)}\) and inspect Table~5.1
	in~\cite{LSS}. Then $G$ has a subgroup \(M\) of structure \(f.\PSL^{\varepsilon}_8(q)\).
	Without loss $P \le M$ and we choose \(Z\unlhd M\) of
	order~\(f\). The coprime action of $P$ on $Z$ yields that \(N_{M/Z}(PZ/Z)=N_M(P)Z/Z\), whence
	 \(|N_M(P)|\) is divisible by \(|N_{M/Z}(PZ/Z)|\).

	As \(M/Z\cong \PSL^{\varepsilon}_8(q)\), we can refer to the arguments for
	\(\PSL_n(q)\) and \(\PSU_n(q)\), and then we find a $p$-subgroup
	\(Q\) such that
	\(|N_{\textrm{E}_6(q)}(Q)|\) is divisible by \(\frac{7\cdot (q-\varepsilon)}{(8,q-\varepsilon)}\).
	Then \(|N_G(P)|\) is divisible by \(7\), and Lemma~\ref{start}(e) leaves two possibilities:
	 \(q=7\), which contradicts Lemma~\ref{Lie1}, %~\ref{strongcases} %9.1,?
  or $q$ and $7$ are coprime, whence Fermat's Little Theorem gives that $7$ divides
	\(q^6-1\). Then Theorem~(10-1)
	in~\cite{GL} shows that \(r_7(G)\ge 3\), contrary to Theorem \ref{strongcases}.
\end{proof}

\begin{lemma}\label{E8q}
There is no set $\Omega$ such that Hypothesis \ref{exc} holds in the case where $G=E_8(q)$.
\end{lemma}

\begin{proof}
We assume otherwise and follow a strategy similar to that in the previous proof, for $E_7(q)$.
First we assume that $p$ divides $|E_7(q)|$.
Then Lemma \ref{start}\,(e) gives that a Sylow $p$-subgroup of $H$ contains a
Sylow $p$-subgroup of a subgroup $M$ of $G$ of structure
$(2,q-1).(\PSL_2(q) \times E_7(q))$. Since $3 \in \pi(\PSL_2(q))$, this contradicts Lemma \ref{start}\,(f).

Together with the table in Lemma \ref{Primesalternativ}, this excludes the cases where $p$ divides
 \(\Phi_7(q)\),
 %=q^6+q^5+q^4+q^3+q^2+q+1\),
 \(\Phi_9(q)\),
 %=q^6+q^3+1\),
	\(\Phi_{14}(q)\), and
	%=q^6-q^5+q^4-q^3+q^2-q+1\), \(\Phi_{15}(q)=q^8-q^7+q^5-q^4+q^3-q+1\),
	\(\Phi_{18}(q)\).
	%=q^6-q^3+1\), \(\Phi_{20}(q)=q^8-q^6+q^4-q^2+1\), \(\Phi_{24}(q)=q^8-q^4+1\) or \(\Phi_{30}(q)=q^8+q^7-q^5-q^4-q^3-q+1\).
	
Now we assume that \(p\) divides \(\Phi_{15}(q)\) or \(\Phi_{30}(q)\). We recall that $P\in \syl_p(H)$ is a Sylow subgroup of $G$ and we let \(T\leq G\) be as
	in~\cite{LSS}. Then without loss \(P\le T\), and
	Table~5.2
	in~\cite{LSS} implies that  \(|N_G(P)|\) is divisible by \(30\), contrary to Lemma \ref{start}(f).

Next we assume that \(p\) divides \(\Phi_{24}(q)=q^8-q^4+1\). Then by Table~5.1
	in~\cite{LSS} we find a subgroup \(M\) of $G$ of structure \(\PSU_3(q^4).8\). Let \(U\leq M\) be isomorphic to \(\PSU_3(q^4)\). Then without loss $P \le U$ and Lemma \ref{help}(a) implies that
	\(|N_M(P)|=8\cdot |N_U(P)|\).
	In particular $|N_G(P)|$ is divisible by $8$, contrary to Lemma \ref{start}(f).

	The final case is that \(p\) divides \(\Phi_{20}(q)=q^8-q^6+q^4-q^2+1\). Then, with Table~5.1
	in~\cite{LSS}, we find  a subgroup \(M\) of $G$ of structure
	\(\SUn_5(q^2).4\). Let \(U\leq M\) be isomorphic to \(\SUn(5,q^2)\).
	
	\smallskip
	\textbf{Case 1:} \(Z(U)\neq 1\).
	
	Then \(|Z(U)|=(5,q^2-1)=5\). Without loss \(P\le U\) and in particular \(Z(U)\leq
	N_G(P)\). Moreover \(|N_M(P)|=4\cdot |N_U(P)|\) by Lemma \ref{help}(a), so \(|N_G(P)|\) is divisible by 20.
	Then Lemma \ref{start}(e) implies that $5 \in \pi(G_\alpha)$.
	We know that \(q\neq 5\) by Lemma 9.1, so
	\(q^4-1\) is divisible by \(5\). But then Theorem~(10-1)
	in~\cite{GL} shows that \(r_5(G) \ge 4\), contrary to Theorem \ref{strongcases}.
	
	\smallskip
	\textbf{Case 2:} \(Z(U)=1\).
	
	Then \(U\cong \PSU_5(q^2)\) and we refer to our results about
	(\(\PSU_n(q)\). Without loss \(P\leq U\) and consequently \(|N_G(P)|\) is divisible by \(\frac{(q^5+1)\cdot
		5}{(q+1)\cdot (5,q+1)}=\frac{q^5+1}{q+1}\cdot 5\).
		Lemma \ref{start}(e) implies that $5 \in \pi(G_\alpha)$, which gives two possibilties:
	
	\(q=5\), which contradicts Lemma 9.1, or \(q^4-1\) is divisible by \(5\).
	In the second case we use Theorem~(10-1)
	in~\cite{GL}, and then we see that \(r_5(G) \ge 4\), contrary to
	Theorem \ref{strongcases}.

\end{proof}

%%%%%

\begin{thm}\label{Liecomplete}
Suppose that $G$ is a simple group of Lie type satisfying Hypothesis \ref{hyp6}.
Then one of the following holds:

$G=\PSL_2(q)$, where the defining characteristic $r$ is at least $5$, $n \in \N$, and $G_\alpha$ has index 2 in the normalizer of a
Sylow $r$-subgroup.

$G=\PSU_4(2)$ or $G=\PSU_4(3)$, and $|G_\alpha|=5$.

$G=\PSp_4(q)$ and $G_\alpha$ is cyclic of order $\frac{q^2+1}{(2,q+1)}$.

$G=\POm^-_8(q)$ and $G_\alpha$ is cyclic of order $\frac{q^4+1}{(2,q+1)}$.

$G=\,^3D_4(q)$ and $G_\alpha$ is cyclic of order $q^4-q^2+1$.

$G=\,^2G_2(q)$ and $G_\alpha$ is cyclic of order $\frac{q-1}{2}$.
\end{thm}

\begin{proof}
Theorem \ref{strongcases} gives two cases:
In Case (ii) we inspect Section 7 for the groups that are mentioned.
In most cases of Theorem \ref{L2q-main}, we see immediately that $G_\alpha$
has order divisible by $2$ or $3$.
The only remaining possibility is that $G_\alpha$ has index 2 in the normalizer of a
Sylow $p$-subgroup of $G$, and for this to be compatible with Hypothesis \ref{hyp6}, we must have that $p \ge 5$.
Lemma \ref{PSLPSU3} gives that $\PSU_3(p^n)$ only occurs with point stabilizers of order divisible by $3$.

In Case (i) of Theorem \ref{strongcases}, we first consider the situation where the defining characteristic
is in $\pi(H)$. Then we use Lemma \ref{Lie1}, which gives a special case of the statement about $\PSL_2(q)$.
In the cross-characteristic case, we go through the results for the individual series of groups of Lie type.
The only possibilities come from Lemmas \ref{ClassicalNotUnitaryCo6}, \ref{3D4q} and
\ref{2G2q}. All other series of groups do not give examples satisfying Hypothesis \ref{hyp6}.
\end{proof}

%%%%%%%%%%%%%%

\section{Proof of Theorem \ref{main4fpsimple}}

%After all this work we can prove our main result. For this we go through the main results of the previous sections and deduce that the possibilities in Table \ref{MainTable} are the only ones where a finite simple group acts with fixity 4.

The hypothesis of Theorem \ref{main4fpsimple} tells us that $G$ is a finite simple non-abelian group that acts transitively and with fixity 4 on a set $\Omega$.
Let $f$ denote the maximum number of fixed points of involutions in $G$.

If Case (1) of Theorem \ref{3and2simple} holds, then $1 \le f \le 4$ and
$G$ is isomorphic to $\PSL_2(q)$, $\Sz(q)$ or $\PSU_3(q)$ for some power $q$ of $2$. Here we use the main result from \cite{Ben}.
Then the corresponding results from Section 5 apply, more precisely Lemmas \ref{L2q-main}, \ref{Suzuki} and \ref{PSLPSU3}.
All these cases are contained in Table \ref{MainTable}.

If Case (2) of Theorem \ref{3and2simple} holds, then $f \le 3$ and
$G$ is isomorphic to $\PSL_2(q)$,  to $\PSL_3(q)$ or to $\PSU_3(q)$ for some odd prime power $q$.
Again we turn to Section 7, this time to Theorem \ref{L2q-main} and Lemma \ref{PSLPSU3}.
All the possibilities there are covered by Table \ref{MainTable}.

If Case (3) of Theorem \ref{3and2simple} holds, then $f=4$,
 $G$ has sectional $2$-rank at most $4$ and
there is a list of possibilities for $G$. We have already seen that all possibilities for
$2$-dimensional Lie type groups are contained in Table \ref{MainTable}, then we refer to
Lemma \ref{Inv-lem:Ser} for the Lie type groups of dimension at least 3,
and the only remaining group is the sporadic group Ly.
But this was discussed in Lemma \ref{Lyons}.

In Case (4) of Theorem \ref{3and2simple} we have a number of sub-cases, all of which
have been treated in Section 8. For the strongly 3-embedded case we use Lemma \ref{str3emb}, and all possibilities listed there are contained in Table \ref{MainTable}. This covers Cases (4) (a) and (b).
For (c) we use Lemma \ref{odd3el9}, which is also a case from our table.
Cases (d) and (e) are covered by Lemma \ref{extraspecialWreath} and do not give any examples.

Finally, we turn to Case (5), and hence to Section 9.
Given Hypothesis \ref{hyp6}, we have already seen the groups that are listed in Theorem \ref{strongcases}(ii), so we only need to consider those with cyclic Sylow $p$-subgroups.
The only example coming from an alternating group here is $\Alt_7$, as explained in Lemma \ref{Altcyclic}, and this is included in Table \ref{MainTable}.
The only sporadic examples are $\M_{11}$ and $\M_{22}$, as can be seen in Lemma \ref{spor}, and these are also in our table.
A summary of the results about Lie type groups is given in
Theorem \ref{Liecomplete}, and this is where the remaining cases from table \ref{MainTable} come from.

%Ausblick auf die nächsten Schritte?

%%%%%%%%%%%%%%%%%%

%%%%%%%%%%%%%%%%%%

%%%%%%%%%%%%%%%%%%%

\section{Appendix: GAP code and more fixed point information}

%Explain input and output.
%
%\vspace{1cm}
%
%
%\texttt{TestTomF4:=function(tom)\\
%local marks, g, fin;\\
%marks := MarksTom(tom);;\\
%fin := [ ];;\\
%for g in [1..Length(marks)] do
%if Set( [ 2 .. Length( marks[g] ) ],\\
%i -> marks[g][i] < 5 )
%= [ true ]
%and ( 4 in marks[g] )\\
%then Add( fin, [ g, StructureDescription( RepresentativeTom(
%tom, g) ), marks[g][1]] );\\
%fi;\\
%od;\\
%return fin;\\
%end;\\}

For a Table of Marks \texttt{t} of a group \(G\) and for a positive integer \texttt{k}
the following function \texttt{TestTom(t,k)} determines if there are transitive actions of \(G\) with fixity \(k\).
If there are none, then an empty list is returned.
Otherwise the output includes the isomorphism type of the point stabilizers, along with some other information, for every type of fixity $k$ action.

\begin{verbatim}
TestTom:=function(t,k)
  #t Table of Marks of a group
  #k fixity to test
  local marks,g,fin;
  marks:=MarksTom(t);;
  fin:=[];;
  for g in [1..Length(marks)] do
    if Set([2..Length(marks[g])],i->marks[g][i]<k+1)=[ true ] and
      (k in marks[g]) and marks[g][1]>k
    then Add(fin,[StructureDescription(RepresentativeTom(t,g)), marks[g]]);
    fi;
  od;
  return fin;
end;
\end{verbatim}

For many simple groups the Table of Marks is already pre-computed and accessible via the GAP-package TomLib~\cite{TomLib}.
%Examples for two of the groups from Remark 4.1, showing that they do not give
%any examples for fixity 4 actions:\\
%
%\begin{verbatim}
%gap> t:=TableOfMarks(PSL(2,5));
%TableOfMarks( Group([ (3,5)(4,6), (1,2,5)(3,4,6) ]) )
%gap> TestTom(t,4);
%[ ]
%\end{verbatim}
As an example, we look at the sporadic Held group  %the function used for a group from Remark~\ref{rem4.1}, showing that
\(\He\) and we see, with the code above, that it does not give
any example for fixity \(4\) actions:

\begin{verbatim}
gap> t:=TableOfMarks("He");
TableOfMarks( "He" )
gap> TestTom(t,4);
[  ]
gap>
\end{verbatim}

%\vspace{1cm}
%In Lemma ... we also need to exclude $\PSL_2(64)$, which is done here:\\
%
%
%\begin{verbatim}
%gap> t:=TableOfMarks("L2(64)");
%TableOfMarks( "L2(64)" )
%gap> TestTom(t,4);
%[  ]
%gap>
%\end{verbatim}

The next example
%for a group from Remark~\ref{rem4.1}
is \(\Alt_6\cong \PSL_2(9)\).
The function \texttt{TestTom} could be used in the following way:

\begin{verbatim}
gap> t:=TableOfMarks("L2(9)");
TableOfMarks( "A6" )
gap> TestTom(t,4);
[ [ "C2", [ 180, 4 ] ], [ "S3", [ 60, 4, 3, 1 ] ], [ "S3", [ 60, 4, 3, 1 ] ],
  [ "C3 x C3", [ 40, 4, 4, 4 ] ], [ "D10", [ 36, 4, 1, 1 ] ],
  [ "(C3 x C3) : C2", [ 20, 4, 2, 2, 2, 2, 2, 2 ] ] ]
\end{verbatim}

We see that there are six different actions of \(\PSL_2(9)\) on sets of cosets with fixity $4$. Two with a point stabilizer isomorphic to \(\Sym_3\)
and one for \(C_2\), \(C_3\times C_3\), \(D_{10}\) and \((C_3\times C_3):C_2\) each.
In each case the size of the set \(\PSL_2(9)\) is acting on is given (for example \(36\), when the point stabilizer are dihedral of order \(10\))
and some additional information about the number of fixed points of elements.

%More fixed point information for the group that we have checked with \texttt{GAP}.
%	We give an example for how to read this table:
%	For the group $G=...$ there are two different ways to act with fixity 4.
%	In the first case the point stabilizers have structure ... and in the second case they have structure ...
%	For both cases, the little box in the right column shows the number of fixed points for each conjugacy class.

%TODO: Die GAP-Funktion ist noch nicht optimal für das, von dem ich glaube, wofür wir sie brauchen, weil sie noch mehr Informationen ausgibt,
%von denen ich nicht sehe, wofür wir sie benötigen. Ich habe jetzt trotzdem ersteinmal alles erklärt, so wie es im Code steht.

More fixed point information for the groups can be gained by the following GAP-function:

\begin{verbatim}
FixedPointsTom:=function(tom, num)
  local G, indC, mat, CN, l;
  G := UnderlyingGroup(tom);;
  indC := Filtered ( [2..Length(MarksTom(tom))],
    i -> IsCyclic( RepresentativeTom( tom, i ) ) );;
  mat := MatTom(tom)[num];;
  CN := OrdersTom(tom);;
  l:=List(indC,i->[ CN[i],Order(Centralizer(G, RepresentativeTom(tom,i))),mat[i]]);;
  return l;
end;
\end{verbatim}

\vspace{2cm}
The input is again a Table of Marks \texttt{tom} of a group and \texttt{num} identifies the specific action we are looking at.
If we return to the example \(G=\PSL_2(9)\), then the table of marks looks like this:

\newpage
\begin{verbatim}
gap> tom:=TableOfMarks("L2(9)");
TableOfMarks( "A6" )
gap> Display(tom);
 1:  360
 2:  180 4
 3:  120 . 6
 4:  120 . . 6
 5:   72 . . . 2
 6:   90 2 . . . 2
 7:   90 6 . . . . 6
 8:   90 6 . . . . . 6
 9:   60 4 3 . . . . . 1
10:   60 4 . 3 . . . . . 1
11:   45 5 . . . 1 3 3 . . 1
12:   40 . 4 4 . . . . . . . 4
13:   36 4 . . 1 . . . . . . . 1
14:   30 2 6 . . . . 2 . . . . . 2
15:   30 2 . 6 . . 2 . . . . . . . 2
16:   20 4 2 2 . . . . 2 2 . 2 . . . 2
17:   15 3 3 . . 1 3 1 1 . 1 . . 1 . . 1
18:   15 3 . 3 . 1 1 3 . 1 1 . . . 1 . . 1
19:   10 2 1 1 . 2 . . 1 1 . 1 . . . 1 . . 1
20:    6 2 3 . 1 . . 2 1 . . . 1 2 . . . . . 1
21:    6 2 . 3 1 . 2 . . 1 . . 1 . 2 . . . . . 1
22:    1 1 1 1 1 1 1 1 1 1 1 1 1 1 1 1 1 1 1 1 1 1

\end{verbatim}

Now \texttt{num} refers to the specific line of this table of marks.
We have seen earlier that \(G=\PSL_2(9)\) acts on a set of size \(20\) and with a point stabilizer \(U\) of structure \((C_3\times C_3):C_2\).
This refers to line \texttt{16} of \texttt{tom}. If we want to know the number of fixed points for every element of \(U\) in this action (which we call the fixed point profile), then we use the following GAP-code:

\begin{verbatim}
gap> FixedPointsTom(tom,16);
[ [ 2, 8, 4 ], [ 3, 9, 2 ], [ 3, 9, 2 ], [ 5, 5, 0 ], [ 4, 4, 0 ] ]
\end{verbatim}

In general the output of \texttt{FixedPointsTom} is a list of lists, where
the first entry in each of the lists refers to the order of an element \(g\in G\), the second to the order of \(C_G(g)\)
and the last one is the number of fixed points of \(g\) in its action on \(G/U\).
So in our case, where \(G=\PSL_2(9)\) and \(U\cong (C_3\times C_3):C_2\), we can see that elements of order \(2\) fix four points, elements of order \(3\)
fix two points independent of their conjugacy class, and non-trivial elements of other orders
have no fixed points on \(G/U\).

Collecting this information for all fixity \(4\) actions we obtain a fixed point profile, see Table~\ref{SGTdetails}, for most of the groups from
Remark~\ref{rem4.1}:

\begin{longtable}{lll}	%HEAD
	$Group$ & point stabilizer structure & corresponding fixed point profile\\[0.5ex]\hline\hline
	\endhead \\[-2ex]
	%MAIN
	\(\Alt_6 \cong \PSL_2(9)\)            & \(C_2\) & \tiny{
\begin{tabular}{|c|c|c|c|c|} \hline
      2A & 3A & 3B & 4A & 5A \\ \hline
      4   &  0  &  0  &  0  &  0 \\ \hline
    \end{tabular}}\\\\
	& \(\Sym_3\)&\tiny{
\begin{tabular}{|c|c|c|c|c|} \hline
      2A & 3A & 3B & 4A & 5A \\ \hline
      4   &  3  &  0  &  0  &  0 \\ \hline
    \end{tabular}}\\\\
    	& \(\Sym_3\)&\tiny{
\begin{tabular}{|c|c|c|c|c|} \hline
      2A & 3A & 3B & 4A & 5A \\ \hline
      4   &  0  &  3  &  0  &  0 \\ \hline
    \end{tabular}}\\\\
	& \(C_3 \times C_3\)&\tiny{
\begin{tabular}{|c|c|c|c|c|} \hline
      2A & 3A & 3B & 4A & 5A \\ \hline
      0   &  4  &  4  &  0  &  0 \\ \hline
    \end{tabular}}\\ \\
	& \(D_{10}\) &\tiny{
\begin{tabular}{|c|c|c|c|c|} \hline
      2A & 3A & 3B & 4A & 5A \\ \hline
      4   &  0  &  0  &  0  &  1 \\ \hline
    \end{tabular}} \\\\
	&  \((C_3 \times C_3) : C_2\)   &  \tiny{
\begin{tabular}{|c|c|c|c|c|} \hline
      2A & 3A & 3B & 4A & 5A \\ \hline
      4   &  2  &  2  &  0  &  0 \\ \hline
    \end{tabular}}\\[5ex]
	\hline
	\\
	\(\Alt_7\)                            & \(C_5\) & \tiny{
\begin{tabular}{|c|c|c|c|c|c|c|} \hline
      2A & 3A & 3B & 4A & 5A & 6A & 7A \\ \hline
      0   &  0  &  0  &  0  &  4  &  0  &  0 \\ \hline
    \end{tabular}}\\\\
	&  \(\Alt_6\) & \tiny{
\begin{tabular}{|c|c|c|c|c|c|c|} \hline
      2A & 3A & 3B & 4A & 5A & 6A & 7A \\ \hline
      3   &  4  &  1  &  1  &  2  &  0  &  0 \\ \hline
    \end{tabular}} \\\\  \hline\\
	
	\(\PSL_2(7) \cong \PSL_3(2)\)         & \(C_2\) & \tiny{
\begin{tabular}{|c|c|c|c|} \hline
      2A & 3A & 4A & 7A \\ \hline
      4   &  0  &  0  &  0 \\ \hline
    \end{tabular}}\\\\
	& \(\Sym_3\)                                                & \tiny{
\begin{tabular}{|c|c|c|c|} \hline
      2A & 3A & 4A & 7A \\ \hline
      4   &  1  &  0  &  0 \\ \hline
    \end{tabular}} \\ \\ \hline\\
	
	\(\PSL_2(8)\)                         & \(C_2\) & \tiny{
\begin{tabular}{|c|c|c|c|} \hline
      2A & 3A & 7A & 9A \\ \hline
      4   &  0  &  0  &  0 \\ \hline
    \end{tabular}}\\\\
	&  \(\Sym_3\) &\tiny{
\begin{tabular}{|c|c|c|c|} \hline
      2A & 3A & 7A & 9A \\ \hline
      4   &  3  &  0  &  0 \\ \hline
    \end{tabular}} \\\\
	& \(D_{14}\) & \tiny{
\begin{tabular}{|c|c|c|c|} \hline
      2A & 3A & 7A & 9A \\ \hline
      4   &  0  &  1  &  0 \\ \hline
    \end{tabular}}\\\\
	& \(D_{18}\) & \tiny{
\begin{tabular}{|c|c|c|c|} \hline
      2A & 3A & 7A & 9A \\ \hline
      4   &  1  &  0  &  1 \\ \hline
    \end{tabular}}\\\\  \hline \\
	
	\(\PSL_2(11)\)                        & \(C_3\) & \tiny{
\begin{tabular}{|c|c|c|c|c|} \hline
      2A & 3A & 5A & 6A & 11A\\ \hline
      0   &  4  &  0  &  0 & 0 \\ \hline
    \end{tabular}}\\\\	
	& \(\Alt_4\)                                            &      \tiny{
\begin{tabular}{|c|c|c|c|c|} \hline
      2A & 3A & 5A & 6A & 11A\\ \hline
      3   &  4  &  0  &  0 & 0 \\ \hline
    \end{tabular}}\\   \\ \hline\\

	\(\PSL_2(13)\)                        & \(C_3\) & \tiny{
\begin{tabular}{|c|c|c|c|c|} \hline
      2A & 3A & 6A & 7A & 13A\\ \hline
      0   &  4  &  0  &  0 & 0 \\ \hline
    \end{tabular}}\\\\
	&  \(\Alt_4\) &  \tiny{
\begin{tabular}{|c|c|c|c|c|} \hline
      2A & 3A & 6A & 7A & 13A\\ \hline
      3   &  4  &  0  &  0 & 0 \\ \hline
    \end{tabular}}\\\\
	& \(C_{13} : C_3\)                                             &
	 \tiny{
\begin{tabular}{|c|c|c|c|c|} \hline
      2A & 3A & 6A & 7A & 13A\\ \hline
      0   &  4  &  0  &  0 & 2 \\ \hline
    \end{tabular}}\\                                                 \\\hline\\

    	\(\PSL_2(17)\)                        & \(C_4\) & \tiny{
\begin{tabular}{|c|c|c|c|c|c|} \hline
      2A & 3A & 4A & 8A & 9A & 17A\\ \hline
      4   &  0  &  4  &  0 & 0 &0 \\ \hline
    \end{tabular}}\\\\
	&  \(C_{17}:C_4\) & \tiny{
	\begin{tabular}{|c|c|c|c|c|c|} \hline
      2A & 3A & 4A & 8A & 9A & 17A\\ \hline
      4   &  0  &  4  &  0 & 0 & 2 \\ \hline
    \end{tabular}}\\
    \\\hline\\

        	\(\PSL_2(19)\)                        & \(C_5\) & \tiny{
	\begin{tabular}{|c|c|c|c|c|c|} \hline
      2A & 3A & 5A & 9A & 10A & 19A\\ \hline
      0   &  0  &  4  &  0 & 0 & 0 \\ \hline
    \end{tabular}}\\
    \\\hline\\

            	\(\PSL_2(23)\)                        & \(C_6\) & \tiny{
	\begin{tabular}{|c|c|c|c|c|c|c|} \hline
      2A & 3A & 4A & 6A & 11A & 12A & 23A\\ \hline
      4   &  4  &  0  &  4 & 0 & 0 & 0\\ \hline
    \end{tabular}}\\
    \\\hline\\

        	\(\PSL_2(25)\)                        & \(C_6\) & \tiny{
\begin{tabular}{|c|c|c|c|c|c|c|c|} \hline
      2A & 3A & 4A & 5A & 5B & 6A & 12A & 13A\\ \hline
      4   &  4  &  0  &  0  & 0   &  4  & 0     &  0\\ \hline
    \end{tabular}}\\\\
	&  \((C_5 \times C_5):C_6\) & \tiny{
\begin{tabular}{|c|c|c|c|c|c|c|c|} \hline
      2A & 3A & 4A & 5A & 5B & 6A & 12A & 13A\\ \hline
      4   &  4  &  0  &  2  & 2   &  4  & 0     &  0\\ \hline
    \end{tabular}}\\
    \\\hline\\

                	\(\PSL_2(27)\)                        & \(C_7\) & \tiny{
	\begin{tabular}{|c|c|c|c|c|} \hline
      2A & 3A & 7A & 13A & 14A \\ \hline
      0   &  0  &  4  &  0 & 0\\ \hline
    \end{tabular}}\\
    \\\hline\\

        	\(\PSL_2(29)\)                        & \(C_7\) & \tiny{
\begin{tabular}{|c|c|c|c|c|c|c|} \hline
      2A & 3A & 5A & 7A & 14A & 15A & 29A\\ \hline
      0   &  0  &  0  &  4  & 0   &  0  & 0 \\ \hline
    \end{tabular}}\\\\
	&  \((C_{29}):C_7\) & \tiny{
\begin{tabular}{|c|c|c|c|c|c|c|} \hline
      2A & 3A & 5A & 7A & 14A & 15A & 29A\\ \hline
      0   &  0  &  0  &  4  & 0   &  0  & 2 \\ \hline
    \end{tabular}}\\
    \\\hline\\

                    	\(\PSL_2(31)\)                        & \(C_8\) & \tiny{
	\begin{tabular}{|c|c|c|c|c|c|c|c|} \hline
      2A & 3A &4A &  5A & 8A & 15A & 16A & 31A\\ \hline
      4   &  0  &  4  &  0  & 4   &  0   &   0    & 0 \\ \hline
    \end{tabular}}\\
    \\\hline\\

        	\(\PSL_2(37)\)                        & \(C_9\) & \tiny{
\begin{tabular}{|c|c|c|c|c|c|c|} \hline
      2A & 3A & 6A & 9A & 18A & 19A & 37A\\ \hline
      0   &  4  &  0  &  4  & 0   &  0  & 0 \\ \hline
    \end{tabular}}\\\\
	&  \((C_{37}):C_9\) & \tiny{
\begin{tabular}{|c|c|c|c|c|c|c|} \hline
      2A & 3A & 6A & 9A & 18A & 19A & 37A\\ \hline
      0   &  4  &  0  &  4  & 0   &  0  & 2 \\ \hline
    \end{tabular}}\\
    \\\hline\\

        	\(\PSL_2(41)\)                        & \(C_{10}\) & \tiny{
\begin{tabular}{|c|c|c|c|c|c|c|c|c|} \hline
      2A & 3A & 4A & 5A  & 7A  & 10A & 20A & 21A & 41A\\ \hline
      4   &  4  &  0  &  4  & 0     &  4   & 0      & 0     & 0   \\ \hline
    \end{tabular}}\\\\
	&  \((C_{41}):C_{10}\) & \tiny{
\begin{tabular}{|c|c|c|c|c|c|c|c|c|} \hline
      2A & 3A & 4A & 5A  & 7A  & 10A & 20A & 21A & 41A\\ \hline
      4   &  4  &  0  &  4  & 0     &  4   & 0      & 0     & 2   \\ \hline
    \end{tabular}}\\
    \\\hline\\

	\(\PSU_3(3)\)                         & \(((C_3 \times C_3) :C_3) :C_8\)                                                  &
	\tiny{
\begin{tabular}{|c|c|c|c|c|c|c|c|c|} \hline
      2A & 3A & 3B & 4A & 4B & 6A & 7A & 8A & 12A\\ \hline
      4   &  1  &  1  &  4 & 0   & 1   & 0    & 2   & 1\\ \hline
    \end{tabular}}\\
	\\\hline\\

	\(\PSU_4(2)\cong \PSp_4(3)\)          & \(C_5\)                                                                           &
	 \tiny{
\begin{tabular}{|c|c|c|c|c|c|c|c|c|c|c|c|c|c|} \hline
      2A & 2B & 3A & 3B & 3C & 4A & 4B & 5A  & 6A & 6B & 6C & 6D & 9A & 12A\\ \hline
      0   &  0  &  0  &  0 & 0    &  0  &   0  & 4  &   0  &  0 &  0  &  0  &  0  &  0  \\ \hline
    \end{tabular}}\\ \\\hline\\

	\(\Sz(8)\)                            & \(C_5\) &  \tiny{
\begin{tabular}{|c|c|c|c|c|} \hline
      2A & 4A & 5A & 7A & 13A\\ \hline
      0   &  0  &  4  &  0 & 0 \\ \hline
    \end{tabular}}\\\\
	& \(C_{13}\)                                                               &
	 \tiny{
\begin{tabular}{|c|c|c|c|c|} \hline
      2A & 3A & 6A & 7A & 13A\\ \hline
      0   &  0  &  0  &  0 &  4 \\ \hline
    \end{tabular}}\\                                                       \\\hline\\
	
	\(M_{11}\)                            & \(C_5\) &  \tiny{
\begin{tabular}{|c|c|c|c|c|c|c|} \hline
      2A & 3A & 4A & 5A & 6A & 8A & 11A\\ \hline
      0   &  0  &  0  &  4  &  0  &  0  &   0\\ \hline
    \end{tabular}}\\\\
	 & \(C_{11} :C_5\) &
	 \tiny{
\begin{tabular}{|c|c|c|c|c|c|c|} \hline
      2A & 3A & 4A & 5A & 6A & 8A & 11A\\ \hline
      0   &  0  &  0  &  4  &  0  &  0  &   1\\ \hline
    \end{tabular}}\\\\
	 & \(\PSL_2(11)\)                                          &
	 \tiny{
\begin{tabular}{|c|c|c|c|c|c|c|} \hline
      2A & 3A & 4A & 5A & 6A & 8A & 11A\\ \hline
      4   &  3  &  0  &  2  &  1  &  0  &   1\\ \hline
    \end{tabular}}\\\\\hline\\
	\(M_{12}\)                            & \(M_{11}\)                                                                        &
	 \tiny{
\begin{tabular}{|c|c|c|c|c|c|c|c|c|c|c|c|c|} \hline
      2A & 2B & 3A & 3B & 4A & 4B & 5A  & 6A & 6B & 8A & 8B & 10A & 11A\\ \hline
      0   &  4  &  3  &  0 & 0    &  4  &  2  & 0   &   1 &  2  &  0  &  0    &  1  \\ \hline
    \end{tabular}}\\  \\
    	                        & \(M_{11}\)                                                                        &
	 \tiny{
\begin{tabular}{|c|c|c|c|c|c|c|c|c|c|c|c|c|} \hline
      2A & 2B & 3A & 3B & 4A & 4B & 5A  & 6A & 6B & 8A & 8B & 10A & 11A\\ \hline
      0   &  4  &  3  &  0 & 4    &  0  &  2  & 0  &   1  &  0 &  2  &  0  &    1  \\ \hline
    \end{tabular}}\\   \\
    \hline\\
	\(M_{22}\)                            & \(C_5\) &  \tiny{
\begin{tabular}{|c|c|c|c|c|c|c|c|c|} \hline
      2A & 3A & 4A & 4B & 5A & 6A & 7A & 8A & 11A\\ \hline
      0   &  0  &  0  &  0  &  4  &  0  &   0 &   0 &   0\\ \hline
    \end{tabular}}\\\\
	&  \(C_{11}:C_5\)                                                           &    \tiny{
\begin{tabular}{|c|c|c|c|c|c|c|c|c|} \hline
      2A & 3A & 4A & 4B & 5A & 6A & 7A & 8A & 11A\\ \hline
      0   &  0  &  0  &  0  &  4  &  0  &   0 &   0 &   1\\ \hline
    \end{tabular}} \\                                                       \\\hline\\

	\(J_1\)                               & \(C_{15}\)                                                                        &
	\tiny{
\begin{tabular}{|c|c|c|c|c|c|c|c|c|} \hline
      2A & 3A & 5A & 6A & 7A & 10A & 11A & 15A & 19A\\ \hline
      0   &  4  &  4  &  0  &  0  &  0  &   0 &   4 &   0\\ \hline
    \end{tabular}} \\
     \\[0.5ex]\hline\hline\\[-1ex]
	\caption[Small Groups Table for fixity 4]{For each group, we give the possible point stabilizer structures for fixity 4 actions, and also the fixed point profile.}
	\label{SGTdetails}

\end{longtable}

%%%%%%%%%%%%%%%%%%

%%%%%%%%%%%%%%%%%%%

%%%%%%%%%%%%%%%%%%

\end{document}